\documentclass[titlepage,11pt]{article}

\usepackage{latexsym}
\usepackage{amsfonts}
\usepackage{amsmath}
\usepackage{amssymb}
\usepackage{comment} 
\usepackage{mathtools}
\usepackage{tikz}
\usepackage{xcolor}
\usetikzlibrary{arrows}
\usetikzlibrary{arrows.meta}
\usetikzlibrary{decorations.markings}
\usetikzlibrary{decorations.pathreplacing,calligraphy}
\usepackage{float}
\newcommand\blackslug{\hbox{\hskip 1pt \vrule width 4pt height 8pt depth 1.5pt
		\hskip 1pt}}
\newcommand\bbox{\hfill \quad \blackslug \bigbreak}

\def\LL{,\ldots,}

\newcommand{\dist}{\operatorname{dist}}

\newcommand{\voro}{\operatorname{\Delta}}

\newcommand{\rk}{\operatorname{rk}}
\newcommand{\gtype}{\operatorname{rk}}
\newcommand{\cupcup}{\cup \cdots\cup}

\newcommand{\bd}{\operatorname{bd}}

\newcommand{\mac}{\mathcal}

\title{Asymptotic structure. III\@. Excluding a fat tree}
\author{
	Tung Nguyen\thanks{Part of this work was done while the first author was at Princeton University and was supported by a Porter Ogden Jacobus Fellowship, AFOSR grant FA9550-22-1-0234, and NSF grant DMS-2154169.
		Currently supported by a Titchmarsh Research Fellowship and a Christ Church Research Centre Grant.}\\
	University of Oxford,\\ Oxford, UK
	\and
	Alex Scott\thanks{Supported by EPSRC grant EP/X013642/1.  AS would like to thank the Department of Mathematics at the University of Arizona, where part of this work was completed.}\\
	University of Oxford, \\
	Oxford, UK
	\and
	Paul Seymour\thanks{Supported by AFOSR grants
		FA9550-19-1-0187 and FA9550-22-1-0234, and by NSF grants DMS-1800053 and DMS-2154169.}\\
	Princeton University,\\ Princeton, NJ 08544, USA}

\date{
}

\newtheorem{thm}{}[section]

\newcommand{\Proof}{\noindent{\bf Proof.}\ \ }

\begin{document}
	\maketitle
	\begin{abstract}
		Robertson and Seymour proved that for every finite tree $H$, there exists $k$ such that every finite graph $G$ with no $H$ minor has path-width at most $k$; and 
		conversely, for every integer $k$, there is a finite tree $H$ such that every finite graph $G$ with an $H$ minor has path-width more than $k$. If we (twice) replace 
		``path-width'' by ``line-width'', the same is true for infinite graphs $G$.
		
		We prove an analogue in coarse graph theory: for every finite tree $H$, graphs without fat $H$-minors are quasi-isometric to graphs with bounded line-width. More precisely, or every finite tree $H$, there exists $k$ such that for every $c\ge 0$, there exist $L,C$ such that every graph 
		that does not contain $H$ as a $c$-fat minor admits an $(L,C)$-quasi-isometry to a graph with line-width at most $k$; and conversely,
		for all $k$ there is a finite tree $H$ such that for all $L,C$ there exists $c$ such that no graph that contains $H$ as a $c$-fat minor admits an $(L,C)$-quasi-isometry to a graph 
		with line-width at most $k$.
	\end{abstract}
	
	\section{Introduction}

	Graphs in this paper may be infinite, and have no loops or parallel edges.  For vertices $x,y$ in a graph $G$ we write $\dist_G(x,y)$ for the distance between $x$ and $y$, that is, the number of edges in the shortest path of $G$ between $x$ and  $y$ (if no path exists we set $\dist_G(x,y) = \infty$).  If $G$ is connected then the distance function is a metric on the vertex set of $G$.  More generally, if each of $X$ and $Y$ is a vertex, or a set of vertices, or a subgraph of $G$, then $\dist_G(X,Y)$ is the minimum of $\dist_G(x,y)$ over vertices $x$ in $X$ and $y$ in $Y$.
	
	Coarse graph theory is concerned with the large-scale geometric structure of graphs, and in particular with properties that are resistant to quasi-isometry.	
	Let $G,H$ be graphs, and let $\phi:V(G)\to V(H)$ be a map.
	Let $L,C\ge 0$; we say that $\phi$ is an {\em $(L,C)$-quasi-isometry}
	if:
	\begin{itemize}
		\item for all $u,v$ in $V(G)$, if $\dist_G(u,v)$ is finite then $\dist_H(\phi(u),\phi(v))\le L \dist_G(u,v)+C$;
		\item for all $u,v$ in $V(G)$, if $\dist_H(\phi(u),\phi(v))$ is finite then $\dist_G(u,v)\le L \dist_H(\phi(u),\phi(v))+C$;
		and
		\item for every $y\in V(H)$ there exists $v\in V(G)$ such that $\dist_H(\phi(v), y)\le  C$.
	\end{itemize}
	
	Many notions from classical graph theory have coarse graph theoretic analogues: in this paper we will be concerned with an analogue of graph minors.	
	Let us first recall the definition of a minor of a graph $H$.  
	Let $G$ be a graph, and for each vertex $x$ of $H$ let $\eta(x)$ be a non-null connected subgraph of 
	$G$, with the $\eta(x)$ pairwise vertex-disjoint.  If, for each $uv\in E(H)$, there is an edge of $G$ between $\eta(u)$ and $\eta(v)$ then we say that $G$ contains $H$ as a {\em minor}, and $\eta$ {\em exhibits} $H$ as a minor 
	of $G$. 
	
	In the coarse graph-theoretic setting, an analogous role is played by {\em fat minors}: these were introduced by 
	Bonamy, Bousquet, Esperet, Groenland, Liu, Pirot and Scott \cite{bonamy}
	(an equivalent notion was considered earlier by Chepoi, Dragan, Estellon, Habib, and Vax\`es \cite{chepoi}).
	If $G$ is a graph, we write $U(G)$ for $V(G)\cup E(G)$. 
	Let $G,H$ be graphs, and let $c\ge 0$ be an integer. For each $x\in U(H)$, let $\eta(x)$ be a non-null connected subgraph of 
	$G$, all pairwise vertex-disjoint, such that 
	\begin{itemize}
		\item for each $uv\in E(H)$, there is an edge of $G$ between $\eta(u)$ and $\eta(uv)$ ($=\eta(vu)$), and an edge between $\eta(v)$ and $\eta(uv)$;
		\item $\dist_G(\eta(x),\eta(y))>c$ for all distinct $x,y\in U(H)$, except when one of $x,y$ is in $V(H)$, the other is 
		in $E(H)$,
		and the edge is incident in $H$ with the vertex. (See Figure \ref{fig:fatminor}.)
	\end{itemize}
	\begin{figure}[H]
		\begin{center}
			\begin{tikzpicture}[]
				
				\tikzstyle{every node}=[inner sep=1.5pt, fill=black,circle,draw]
				\def\r{1}
				\node (v0) at (0,0) {};
				\node (v1) at ({\r* cos(90)}, {\r*sin(90)}) {};
				\node (v2) at ({\r* cos(162)}, {\r*sin(162)}) {};
				\node (v3) at ({\r* cos(234)}, {\r*sin(234)}) {};
				\node (v4) at ({\r* cos(306)}, {\r*sin(306)}) {};
				\node (v5) at ({\r* cos(18)}, {\r*sin(18)}) {};
				\foreach \from/\to in {v0/v1,v0/v2,v0/v3,v0/v4,v0/v5,v1/v2,v2/v3,v3/v4,v4/v5,v5/v1}
				\draw [-] (\from) -- (\to);
				
			\end{tikzpicture}
		\end{center}
		
		\begin{center}
			\begin{tikzpicture}[]
				
				\tikzstyle{every node}=[inner sep=1.5pt, fill=black,circle,draw]
				\def\r{.5}
				\def\s{2}
				\node (v0) at (0,0) {};
				\node (v1) at ({\s* cos(90)}, {\s*sin(90)}) {};
				\node (v2) at ({\s* cos(162)}, {\s*sin(162)}) {};
				\node (v3) at ({\s* cos(234)}, {\s*sin(234)}) {};
				\node (v4) at ({\s* cos(306)}, {\s*sin(306)}) {};
				\node (v5) at ({\s* cos(18)}, {\s*sin(18)}) {};
				\foreach \from/\to in {v0/v1,v0/v2,v0/v3,v0/v4,v0/v5,v1/v2,v2/v3,v3/v4,v4/v5,v5/v1}
				\draw [-] (\from) -- (\to);
				
				\draw[fill= gray!20] (0,0) circle (\r) {};
				\draw[fill= gray!20] ({\s* cos(90)}, {\s*sin(90)}) circle (\r) {};
				\draw[fill= gray!20] ({\s* cos(162)}, {\s*sin(162)}) circle (\r) {};
				\draw[fill= gray!20]  ({\s* cos(234)}, {\s*sin(234)}) circle (\r) {};
				\draw[fill= gray!20] ({\s* cos(306)}, {\s*sin(306)}) circle (\r) {};
				\draw[fill= gray!20] ({\s* cos(18)}, {\s*sin(18)}) circle (\r) {};
				
				\def\t{1}
				\def\a{.2}
				\def\b{.4}
				\draw[fill= white] ({\t* cos(90)}, {\t*sin(90)}) ellipse ({\a} and {\b});
				\draw[rotate around= {72:({\t* cos(162)}, {\t*sin(162)})}, fill= white] ({\t* cos(162)}, {\t*sin(162)}) ellipse ({\a} and {\b}) {};
				\draw[rotate around= {144:({\t* cos(234)}, {\t*sin(234)})}, fill= white]  ({\t* cos(234)}, {\t*sin(234)}) ellipse ({\a} and {\b}) {};
				
				\draw[rotate around= {216:({\t* cos(306)}, {\t*sin(306)})}, fill= white] ({\t* cos(306)}, {\t*sin(306)}) ellipse ({\a} and {\b}) {};
				\draw[rotate around= {288:({\t* cos(18)}, {\t*sin(18)})}, fill= white] ({\t* cos(18)}, {\t*sin(18)}) ellipse ({\a} and {\b}) {};
				
				\def\a{.2}
				\def\b{.55}
				\def\t{1.6}
				\draw[rotate around= {90:({\t* cos(270)}, {\t*sin(270)})}, fill= white] ({\t* cos(270)}, {\t*sin(270)}) ellipse ({\a} and {\b});
				\draw[rotate around= {162:({\t* cos(342)}, {\t*sin(342)})}, fill= white] ({\t* cos(342)}, {\t*sin(342)}) ellipse ({\a} and {\b});
				\draw[rotate around= {234:({\t* cos(54)}, {\t*sin(54)})}, fill= white] ({\t* cos(54)}, {\t*sin(54)}) ellipse ({\a} and {\b});
				\draw[rotate around= {306:({\t* cos(126)}, {\t*sin(126)})}, fill= white] ({\t* cos(126)}, {\t*sin(126)}) ellipse ({\a} and {\b});
				\draw[rotate around= {18:({\t* cos(198)}, {\t*sin(198)})}, fill= white] ({\t* cos(198)}, {\t*sin(198)}) ellipse ({\a} and {\b});

			\end{tikzpicture}
		\end{center}
		\caption{A fat minor.} \label{fig:fatminor}
	\end{figure}
	In these circumstances, we say that $G$ contains $H$ as a {\em $c$-fat minor}, and $\eta$ {\em exhibits} $H$ as a $c$-fat minor 
	of $G$. 
	
	If $G$ is a graph and $X\subseteq V(G)$, then $G[X]$ denotes the subgraph of $G$ induced on $X$.
	Sometimes, we are given a subset $X\subseteq V(G)$, and the function $\eta$ satisfies that $\eta(x)\subseteq G[X]$ for all $x\in U(H)$.
	In this case, we say that {\em $X$ includes a $c$-fat $H$-minor of $G$}. Note that this is not the same thing as saying that $G[X]$
	contains $H$ as a $c$-fat minor, because being $c$-fat depends on a distance function, and we are using the distance function of $G$ rather than that
	of $G[X]$.  If $\eta$ exhibits $H$ as a $c$-fat minor of $G$, then
	$\eta(H)$ denotes the subgraph $\bigcup_{x\in U(H)}\eta(x)$.
	
	In the classical graph-theoretic setting, the structure of graphs that do not contain an $H$-minor is well understood: the graph minor structure theorem of Robertson and Seymour~\cite{GM16} says (roughly) that $H$-minor-free graphs can be built by gluing together graphs that can (almost) be embedded on a surface of bounded genus.  But what about graphs that do not contain a {\em fat} $H$-minor?   Do they have a similar structure to graphs with no $H$-minor, at least up to quasi-isometry?
	
	It is easy to see that if $G$ is quasi-isometric to a graph with no $H$-minor then $G$ does not contain a (suitably) fat $H$-minor.  More precisely:
	\begin{thm}\label{easyfat}
		Let $L,C,c\ge 0$, and suppose that $G$ contains $H$ as a $c$-fat minor, and there is an $(L,C)$-quasi-isometry from $G$ to some graph $G'$. 
		If $c\ge L(L+C)+C$ then 
		$G'$ contains $H$ as a minor.
	\end{thm}
	
	Georgakopoulos and Papasoglu~\cite{agelos} conjectured that a converse also holds:
	\begin{thm}\label{fatminorconj}
		{\bf False conjecture:} For every graph $H$ and all $c\ge 0$, there exists $L,C\ge 0$, such that if a graph $G$ does not contain $H$ as a $c$-fat minor, then $G$ admits an $(L,C)$-quasi-isometry to a graph with no $H$ minor.
	\end{thm}
	This is known to be true when:
	\begin{itemize}
		\item $H=K_3$, by Manning's theorem~\cite{manning} (and see~\cite{agelos}); 
		\item $H=K_{1,m}$, by Georgakopoulos and Papasoglu~\cite{agelos};
		\item $H=K_{2,3}$, by Chepoi, Dragan, Newman, Rabinovich, and Vax\`es~\cite{chepoi};
		\item $H=K_4^-$ (that is, $K_4$ with one edge deleted) by Fujiwara and Papasoglu~\cite{fuji};
		\item $H=K_{2,t}$, by Albrechtsen, Distel and Georgakopoulos~\cite{agelos2}.
	\end{itemize}
	However, the conjecture has recently been shown to be false in general, by Davies, Hickingbotham,
	Illingworth and McCarty~\cite{davies}, building on the construction of \cite{coarsemenger}; and even more recently, Albrechtsen, Distel and Georgakopoulos~\cite{agelos3} showed that it is false when
	$H$ is the graph of the octahedron.
	There remains hope that 
	\ref{fatminorconj} is true for some nontrivial classes of graphs $H$ (for instance, it may be true for all trees $H$ as far as we know). But, for general $H$, the large-scale structure of graphs without a fat $H$-minor is still mysterious.

	Our result is not exactly of this type, but has the same flavour. 
	A {\em path-decomposition} of a graph $G$ is a pair $(T, (B_t:t\in V(T)))$, where $T$ is a path (possibly infinite)
	and each $B_t$ is a subset of $V(G)$ (called a {\em bag}), such that:
	\begin{itemize}
		\item $\bigcup_{t\in V(T)}B_t = V(G)$;
		\item for every edge $e=uv$ of $G$, there exists $t\in V(T)$ with $u,v\in B_t$; and
		\item for all $t_1,t_2,t_3\in V(T)$, if $t_2$ lies between $t_1,t_3$ in $T$, then
		$B_{t_1}\cap B_{t_3}\subseteq B_{t_2}$.
	\end{itemize}
	The {\em width} of a path-decomposition $(T,(B_t:t\in V(T)))$ is the maximum of the numbers $|B_t|-1$ for $t\in V(T)$,
	or $\infty$ if there is no finite maximum;
	and the {\em path-width} of $G$ is the minimum width of a path-decomposition of $G$.
	
	The path-width theorem of Robertson and Seymour~\cite{GM1} shows that the path-width of a graph is closely related to the tree minors that it contains:
	\begin{thm}\label{GM1} 
		For every finite tree $H$, there exists $k$ such that every finite graph $G$ with no $H$ minor has path-width at most $k$; and conversely, for every integer $k$, there is a finite tree $H$ such that every finite graph with an $H$ minor has path-width more than $k$.
	\end{thm}
	It is important here that $G$ is finite: for instance, if $G$ is the disjoint union of countably many infinite paths, then $G$ does not contain
	the claw $K_{1,3}$ as a minor, and yet $G$ does not have finite path-width. This can be fixed. Say a {\em line-decomposition} of $G$ is a family
	$(B_t:t\in T)$ of subsets of $V(G)$, where $T$ is a linearly-ordered set, satisfying the same three bullets above (with $V(T)$ replaced by $T$ everywhere).
	The {\em width} of a line-decomposition $(B_t:t\in T)$ is the maximum of the numbers $|B_t|-1$ for $t\in T$,
	or $\infty$ if there is no finite maximum;
	and the {\em line-width} of $G$ is the minimum width of a line-decomposition of $G$.  For finite graphs, line-width and path-width are the same, but line-width works better with infinite graphs: for example, it is easy to see that the line-width of an infinite graph $G$ is the supremum of the line-widths of its components.   We proved in \cite{subtrees} that $G$ has line-width at most $k$
	if and only if every finite subgraph has path-width at most $k$. Consequently, if we replace ``path-width'' by ``line-width'' in \ref{GM1}, then the statement
	of \ref{GM1} is true even for infinite graphs $G$.

	Our aim in this paper is to give a coarse graph theory analogue of this ``line-width'' version of the path-width theorem \ref{GM1}.
	Our main result says:
	\begin{thm}\label{mainthm}
		For every finite tree $H$ there exists $k\ge 0$ such that for every $c\ge 0$, there exist $L,C\ge 0$ such that every graph
		that does not contain $H$ as a $c$-fat minor admits an $(L,C)$-quasi-isometry to a graph with line-width at most $k$. 
		
		Conversely,
		for every $k$ there is a finite tree $H$ such that for all $L,C$ there exists $c$ such that no graph that contains $H$ as a $c$-fat minor admits an
		$(L,C)$-quasi-isometry to a graph with line-width at most $k$.
	\end{thm}
	The second half is easy: choose $c> L(L+C)+C$, and let $H$ be a finite tree with line-width more than $k$. If $G$ contains $H$ as a $c$-fat minor,
	and there is an $(L,C)$-quasi-isometry from $G$ to $G'$, then by \ref{easyfat}, $G'$ contains $H$ as a minor and hence has line-width more than $k$. 
	The first half is much more difficult and will occupy the whole paper, except for the final section which concerns graph searching.
	
	As a first step, let us eliminate the quasi-isometry. If $X\subseteq V(G)$, we say $X$ has {\em quasi-size} at most $(k,r)$
	if there is a set $Y\subseteq V(G)$ with $|Y|\le k$, such that $\dist_G(x,Y)\le r$ for each $x\in X$. 
	If $X\subseteq V(G)$, let us say a line-decomposition $(B_t:t\in T)$ of $G[X]$ has {\em quasi-width} at most $(k,r)$ in $G$ if 
	$B_t$ has quasi-size at most $(k,r)$ for each $t\in T$. Let us say $X$ has {\em quasi-line-width} at most $(k,r)$
	if $G[X]$ admits a line-decomposition with quasi-width at most $(k,r)$.

	It is important that the distance function used in the definition of 
	quasi-size and quasi-line-width is that defined by $G$, not that defined by $G[X]$. This will be the case throughout the paper. 
	Even when speaking of subgraphs of $G$, we will never use their distance functions: the distance function in use will always be that of $G$.
	We sometimes write $X$ for $G[X]$ when $X\subseteq V(G)$. This should cause no confusion since there is always only one graph $G$ under consideration.

	It was proved in~\cite{coarsepw}, strengthening results of \cite{hick} and an earlier version of \cite{coarsetw}, that:
	\begin{thm}\label{quasitw}
		For all $k,r$, there exist $L,C\ge 1$ such that if $G$ has quasi-line-width at most $(k,r)$,
		then
		$G$ admits an $(L,C)$-quasi-isometry to a graph with line-width at most $k$.
	\end{thm}
	Thus, to complete the proof of \ref{mainthm}, it suffices to prove the following:
	\begin{thm}\label{mainthm2}
		For every finite tree $H$ there exists $k\ge 0$ such that for every $c\ge 0$, there exists $r= O(c)$ such that every graph
		that does not contain $H$ as a $c$-fat minor has quasi-line-width at most $(k,r)$.
	\end{thm}
	
	$H_0$ denotes the tree with one vertex. 
	For $\ell\ge 1$ an integer, let $H_\ell$ be the finite tree such that every vertex has degree one or three, and for some vertex $r$ (called the {\em root})
	every path from $r$ to a vertex of degree one has length exactly $\ell$. Every tree $H$ is a minor of $H_\ell$  for some choice of
	$\ell$, and then, if $G$ does not contain $H$ as a $c$-fat minor then it does not contain $H_\ell$ as a $c$-fat minor. Consequently
	it suffices to prove \ref{mainthm2} when $H=H_\ell$, that is:
	\begin{thm}\label{mainthm3}
		For every $\ell\ge 1$, there exists $k\ge 0$ such that for every $c\ge 0$, there exists $r=O(c)$ such that every graph
		that does not contain $H_\ell$ as a $c$-fat minor has quasi-line-width at most $(k,r)$.
	\end{thm}
	\section{Buildings and tie-breakers}
	
	If $X\subseteq V(G)$, $\bd(X)$ denotes
	the set of vertices in $X$ that have a neighbour in $V(G)\setminus X$, and is called the {\em boundary} of $X$.
	A key idea of the proof is that we work with subsets $X\subseteq V(G)$ such that, simultaneously, $G[X]$ has bounded quasi-line-width and $\bd(X)$ has bounded quasi-size, and it turns out that we can make the two bounds the same with little loss.
	Let us say $X\subseteq V(G)$ has {\em quasi-bound} at most $(a,b)$ if $G[X]$ has quasi-line-width at most $(a,b)$ and $\bd(X)$ has quasi-size at most $(a,b)$.

	A {\em tie-breaker} in $G$ is
	a well-order $\Lambda$ of the set of all edges of $G$ (this exists, by the well-ordering theorem). 
	If $P,Q$ are distinct finite paths of $G$, we say $P$ is {\em $\Lambda$-shorter} than $Q$
	if either
	\begin{itemize}
		\item $|E(P)|<|E(Q)|$; or
		\item $|E(P)|=|E(Q)|$, and the first element (under $\Lambda$) of $(E(P)\setminus E(Q))\cup (E(Q)\setminus E(P))$ belongs to $P$.
	\end{itemize}
	This defines a total order on the set of all finite paths of $G$.
	A {\em $\Lambda$-geodesic} means a finite path $P$ such that no other path joining its ends is $\Lambda$-shorter than $P$.
	Every $\Lambda$-geodesic of $G$ is a geodesic of $G$, but the converse is false.
	(The point of the tie-breaker is that there is only
	one $\Lambda$-geodesic between any two vertices, while this is not true for geodesics.)
	It is easy to check that
	if $P$ is a $\Lambda$-geodesic then so are all subpaths of $P$. Given a graph $G$, we will keep some tie-breaker $\Lambda$ fixed, 
	and usually suppress the dependence of other objects on the choice of $\Lambda$.
	
	A {\em building} in a graph $G$ is a nonempty subset $X\subseteq V(G)$ such that $G[X]$ is connected.
	If $\mathcal{T}$ is a set of pairwise vertex-disjoint buildings in a graph $G$, we define $V(\mac T)=\bigcup_{X\in \mac T}X$. 
	Again, let $\mathcal{T}$ be a set of pairwise vertex-disjoint buildings, in a graph $G$, now with a tie-breaker $\Lambda$. 
	We say that $X\in \mathcal{T}$ is {\em $\Lambda$-closest} to $v$ if there is a path of $G$ between $v,X$ that is $\Lambda$-shorter
	than any other path between $v$ and $Y$ for $Y\in \mathcal{T}$.
	For each $X\in \mathcal{T}$, let $\voro_{\mathcal{T}}(X)$ (or just $\voro(X)$, when $\mathcal{T}$ is clear) be the set of all $v\in V(G)$
	such that $X$ is $\Lambda$-closest to $v$. We call $\voro(X)$ the {\em Voronoi cell} of $X$, and  the collection of subsets $\voro(X)\;(X\in \mathcal{T})$ is called the
	{\em Voronoi partition} defined by $\mathcal{T}$. (It is indeed a partition, provided that each component of
	$G$ includes a member of $\mathcal{T}$; and this is true for us since we will always assume that $G$ is connected and $\mathcal{T}$ is non-null.)
	The main purpose of the tie-breaker is to make the Voronoi partition defined by $\mathcal{T}$ well-defined.
	We see that:
	\begin{itemize}
		\item $X\subseteq \voro(X)$ for each $X\in \mathcal{T}$;
		\item the sets $\voro(X)\;(X\in \mathcal{T})$ are pairwise disjoint and have union $V(G)$;
		\item for each $X\in \mathcal{T}$ and each $v\in \voro(X)$, $\dist_G(v,X)\le \dist_G(v,Y)$ for each $Y\in \mathcal{T}$, and there is
		a path of $G[\voro(X)]$ between $v$ and $X$ of length $\dist_G(v,X)$.
	\end{itemize}
	
	If $C$ is a set of buildings, we define $\voro_{\mac T} (C) = \bigcup_{X\in C}\voro_{\mac T} (X)$. 
	If $\mac T$ is a set of pairwise vertex-disjoint buildings, and $X,Y\in \mac T$ are distinct,  we say that $X$ {\em adjoins} $Y$ 
	(in $\mac T$) or {\em $\mac T$-adjoins} $Y$ if there is an edge between $\voro_{\mac T}(X)$ and $\voro_{\mac T}(Y)$. A set 
	$C\subseteq \mac T$ is {\em adjoin-connected} (or {\em $\mac T$-adjoin-connected} in case of ambiguity) if for every partition of $C$ into two nonempty sets $A,B$, some member of $A$ 
	adjoins some member of $B$. 
	
	Here is an easy lemma.
	\begin{thm}\label{increaserad}
		Let $X\subseteq V(G)$ have quasi-line-width at most $(a,b)$, and let $\bd(X)$ have quasi-size at most $(a',b')$. Let 
		$Y\subseteq V(G)$ with $X\subseteq Y$, such that every vertex in $Y\setminus X$ has distance at most $r$ from $\bd(X)$. Then $\bd(Y)$ has quasi-size at
		most $(a',b'+r)$, and $Y$ has quasi-line-width at most $(a+a',\max(b,b'+r))$.
	\end{thm}
	\Proof
	Let $Z=Y\setminus X$. Every vertex in $Z$ has distance at most $r$ from $\bd (X)$, and the latter has quasi-size at most $(a',b')$;
	so $Z\cup \bd (X)$ has quasi-size at most $(a',b'+r)$. Since $\bd (Y)\subseteq Z\cup \bd (X)$ (since $X\subseteq Y$), it follows that $\bd(Y)$ has quasi-size at
	most $(a',b'+r)$. By adding $Z$ to each bag of a line-decomposition of $G[X]$, we deduce that $Y$ has quasi-line-width at most 
	$(a+a',\max(b,b'+r))$. This proves \ref{increaserad}.~\bbox
	\section{Societies}
	
	For inductive purposes, it is sometimes helpful to strengthen ``fat'' to ``superfat''.
	With notation as before, let $\ell\ge 1$, let $r$ be the root of $H_\ell$,
	let $B_1,B_2,B_3$ be the three components of $H_\ell\setminus \{r\}$, and for $i = 1,2,3$, let $e_i$ be the edge of $H_\ell$ between
	$r$ and $V(B_i)$. We say that $\eta$ exhibits
	$H_\ell$ as a {\em $c$-superfat} minor of a graph $G$ if $\eta$ exhibits $H_\ell$ as a $c$-fat minor of $G$,
	and $\dist_G(x,y)> 3c$ for all $i,j$ with $1\le i<j\le 3$, and all
	$x\in V(\eta(B_i))\cup V(\eta(e_i))$, and all $y\in  V(\eta(B_j))\cup V(\eta(e_j))$.
	When $\ell=0$, we say that $\eta$ exhibits $H_0$ as a $c$-superfat minor of a graph $G$ if it exhibits $H_0$ as a $c$-fat minor of $G$.
	
	We need some more definitions. 
	Let $c,\ell,d_0\ge 2$ be integers, fixed throughout the paper. We will be concerned with graphs that do not contain $H_\ell$ as a $c$-fat
	minor. ($d_0$ is a large number, much larger than $c$ and $\ell$, that we will specify later; in fact we will define $d_0=5\cdot 4^{5\ell^2}c$, but its exact value will not matter until the end of the paper.)
	A {\em century} is an integer $k$ with  $0\le k\le \ell$. 
	Let $k$ be a century, and let $\tau = (d, \alpha,\beta)$ be a triple of
	three positive integers; we call $\tau$ a {\em canon}.
	Let $\Lambda$ be a tie-breaker in a connected graph $G$, and use it to define Voronoi cells as usual.
	A {\em $k$th-century $\tau$-society}
	in $G$ is a set $\mathcal{T}$ of pairwise vertex-disjoint buildings in $G$, where each member of $\mac T$ is assigned to be a ``house'' or ``fort'' of $\mac T$ (and not a house if $k=0$),
	satisfying the following (where $\rk(X)=k-1$ if $X$ is a house and $\rk(X)=k$ if $X$ is a fort):
	\begin{itemize}
		\item For all $v\in V(G)$, there exists $X\in \mathcal{T}$ such that $\dist_G(v,X)\le d_0$;
		\item Every two distinct members $X,Y$ of $\mac T$ have distance more than $d$.
		\item Every fort of $\mac T$ includes a $c$-superfat $H_k$-minor of $G$.
		\item Each fort of $\mac T$ has quasi-bound at most $(\alpha, \beta)$.
		\item $\voro_{\mac T}(C)$ has quasi-bound at most $(\alpha,\beta)$, for every $\mac T$-village $C$.
		(A {\em $\mac T$-community} is an adjoin-connected set of houses, and a {\em $\mac T$-village} is a maximal $\mac T$-community.)
	\end{itemize}
	
	We will not need to work with arbitrary canons; we will mostly just be interested in canons
	$\tau = (d, \alpha,\beta)$ where $d$ is large compared with $c$. Let us call them ``suitable'' for now.
	The strategy for the proof of the main theorem \ref{mainthm3} is to prove that:
	\\
	\\
	{\bf (*) }{\em For every century $k$, and every suitable canon $\tau=(d, \alpha,\beta)$, 
		if $G$ has no $c$-fat $H_\ell$ minor, and admits a $k$th-century $\tau$-society, then $G$ has bounded quasi-line-width (where the bound depends on 
		$k,\tau$).}
	\\
	\\
	If we could prove (*) when $k=0$, then \ref{mainthm3}
	follows, since every graph admits a $0$th-century $(d,1,1)$-society (just pick a maximal set of vertices 
	pairwise at distance more than $d$, and call them all forts). On the other hand, (*) is easy to prove when $k=\ell$. So we will work
	by induction on $\ell-k$: assume that (*) is true for $k+1$, and try to prove that it is also true for $k$.
	
	For the inductive step, we can assume that $G$ admits a $k$th-century $\tau$-society (for some $k<\ell$ and some suitable 
	canon $\tau$). We need to bound the quasi-line-width of $G$, and we do so by showing that $G$
	also contains some $(k+1)$st-century $\tau'$-society (for some suitable $\tau'$ depending on $\tau$)
	and apply the inductive hypothesis. Getting from a $k$th-century society to a $(k+1)$st-century one is thus the main part of the paper, 
	and will be accomplished by first converting our society to a ``realm'', and optimizing the realm, and then extracting from that 
	(via ``governments'' for the realm) the $(k+1)$st-century society we want. But in preparation for converting the given society to a 
	realm, we first need to ``civilize'' it. 
	
	With $\tau$ as before, let $\mac T$ be a $k$th-century $\tau$-society in a graph $G$. We say $\mac T$ is {\em civilized} if
	for every $\mac T$-village $C$, $V(C)$ has quasi-line-width at most $(\alpha,\beta)$ and
	$\bigcup_{X\in C}\bd (X)$ has quasi-size at most $(\alpha,\beta)$ (and hence the same holds for every $\mac T$-community $C$, since both these 
	properties are inherited under taking subsets). This is a great strengthening of the final ``society'' axiom, and we will carry it out
	in the next two sections.

	\section{The grouping lemma}
	
	This section proves a key result 
	that we call the ``grouping lemma''. 
	Suppose we have some disjoint connected subgraphs of a graph $G$, and we would like to join some of them together with short paths, such
	that
	at the end, different components are not very close together, and each vertex in the short paths we used for joining
	is close to one of the original subgraphs. This is not always possible, for instance, if $G$ is a uniform binary tree of large depth,
	and the subgraphs we begin with are the leaves of the tree $G$. But the grouping lemma says that it is possible, provided that
	$G$ does not contain some tree as a fat minor. The argument is similar to that in section 4 of~\cite{AS5}.
	
	For each $i \ge 1$, let $B_i$ be the rooted binary tree of depth $i$. Thus, $B_i$ is the finite rooted tree such that every vertex has degree one or three except for the root $r$ which has degree two, and
	every path from $r$ to a vertex of degree one has length exactly $i$. (This differs from the tree $H_i$ previously defined in that
	the root now has degree two, and in $H_i$ the root has degree three.) Let $B_0$ be the tree with one vertex.
	
	A {\em $(\ge d)$-subdivision} of a graph $H$ is a graph $H'$ obtained by subdividing each edge at least $d$ times. Thus,
	$V(H)\subseteq V(H')$.
	Each edge $e$ of $H$ is replaced by a path joining the ends of $e$, and we call this path the {\em subdivided $e$}. For each subgraph $X$ of $H$, the {\em subdivided $X$} means the subgraph of $H'$ consisting of the vertices in $V(X)$ together with the union of the subdivided $e$ over all edges $e$ of $X$.
	
	Suppose that $G$ is a graph, with a subgraph $J$ that is a $(\ge d)$-subdivision
	of $B_{\ell}$, where $\ell\ge 1$. Suppose, moreover, that for each vertex $v\in B_\ell\subseteq V(J)$ with degree at least two in $B_\ell$, if $X, Y$
	are the two components of $B_\ell\setminus v$ that
	do not contain the root of $B_\ell$, and $e,f$ are the edges of $B_\ell$ between $v$ and $X,Y$ respectively, then the union
	of the subdivided $e$ and the subdivided $f$ is a geodesic in $G$ between the subdivided $X$ and the subdivided $Y$. In this case we
	call $J$ a {\em geodesic $(\ge d)$-subdivision} of $B_\ell$.
	
	In order to prove the grouping lemma, we first need the following:
	\begin{thm}\label{gettree}
		Let $\ell\ge 1$ and $d\ge 3c\ge 6$ be integers. Suppose that $G$ is a graph, with a subgraph $J$ that is a geodesic $(\ge d)$-subdivision
		of $B_{\ell}$.
		Then $G$ contains $B_\ell$ as a $c$-fat minor.
	\end{thm}
	\Proof From the definition of subdivision, $V(B_\ell)\subseteq V(G)$. For each edge or subgraph $X$ of $B_\ell$, let $\phi(X)$ be the subdivided $X$.
	Let $r$ be the root of $B_\ell$. For each $v\in V(B_\ell)$ of degree at least two, let $X_v$, $Y_v$ be the two components of
	$B_\ell \setminus v$ that do not contain $r$, and let $e_v, f_v$ be the edges of $B_\ell$ between $v$ and $X_v, Y_v$ respectively.
	For each $v\in V(B_\ell)$ of degree at least two, let $\eta(v)$ be the subpath of the geodesic $\phi(e_v)\cup \phi(f_v)$
	with vertex set all vertices of $\phi(e_v)\cup \phi(f_v)$ with distance at most $2c$ from $v$. For each $v\in V(B_\ell)$ with
	degree one, let $\eta(v)$ be the one-vertex subgraph with vertex $v$. For each edge $e$ of $B_\ell$, with ends $u,v$ say,
	where $v$ is closer than $u$ to $r$, let $\eta(e)$ be the subpath of $\phi(e)$ with vertex set all vertices not in
	$\eta(u)\cup \eta(v)$.
	We claim that $\eta$ exhibits $B_\ell$ as a $c$-fat minor of $G$. To show this, it
	suffices to check that if $x,y$ are both vertices or edges of $B_\ell$, and it is not the case that
	one of
	$x,y$ is in $V(B_\ell)$, the other is
	in $E(B_\ell)$,
	and the edge is incident in $B_\ell$ with the vertex, then $\dist_G(\eta(x),\eta(y))>c$. Let $P$ be the minimal path of $B_\ell$
	that contains both $x,y$ (thus, $V(P)\subseteq V(B_\ell)\subseteq V(G)$), and let $v$ be the vertex of $P$ that is closest to $r$ in $B_\ell$. (See Figure \ref{fig:gettreecase1}.)
	
	Suppose first that $v=x$.
	Then we may assume that $y$ belongs to $Y_v$, and hence $\eta(y)$ is a subgraph of $\phi(Y_v)$. Since $\phi(f_v)$ has length more than $d$, and is a geodesic between $v$ and $\phi(Y_v)$, it follows that $\dist_G(v,\eta(y))\ge d+1$, and therefore
	$\dist_G(\eta(x),\eta(y))\ge d+1-2c>c$ as required.
	
	Thus we may assume that $x,y\ne v$.
	Suppose next that $x=e_v$ and $y$ belongs to $X_v$. Let the ends of $x=e_v$ in $B_\ell$ be $u,v$.  Then $y\ne u$, and so we may
	assume that either $y=f_u$ or $y$ belongs to $Y_u$. In either case, $\dist_G(u, \eta(y))\ge 2c+1$, since there is a geodesic between
	$u,\eta(y)$ that properly includes $\eta(u)\cap \phi(f_u)$, and the latter has length $2c$.
	Since $\phi(e_v)$
	is a geodesic between $v,\phi(X_v)$, and includes $\eta(x)$, it follows that for every vertex $z\in \eta(x)$,
	$\dist_G(z,\eta(y))\ge \dist_G(z,u)$. Therefore
	$$2c+1\le \dist_G(u,\eta(y))\le \dist_G(u,z) + \dist_G(z,\eta(y))\le 2\dist_G(z,\eta(y)),$$
	and consequently $\dist_G(z,\eta(y))>c$ as required.

	\begin{figure}[H]
		\centering
		\begin{tikzpicture}[scale = 1,auto=left]
			\tikzstyle{every node}=[inner sep=1.5pt, fill=black,circle,draw]
			\node (v) at (0,0) {};
			\node (u) at (-1,1) {};
			\node (w) at (1,1) {};
			\node (x) at (-1.5, 2) {};
			\node (y) at (-.5,2) {};
			\draw (v) to (u);
			\draw (v) to (w);
			\draw (u) to (x);
			\draw (u) to (y);
			\draw[dotted, thick] (v) to (0,-.5);
			
			\draw (-1.5,2.5) ellipse (.4 and .8);
			\draw (-.5,2.5) ellipse (.4 and .8);
			\draw (1,1.5) ellipse (.4 and .8);
			
			\tikzstyle{every node}=[]
			
			\draw (v) node[left] {$v$};
			\draw (u) node[left] {$u$};
			\node at (1,1.5) {$Y_v$};
			\node at (-1.5,2.5) {$X_u$};
			\node at (-.5,2.5) {$Y_u$};
			\node at (-.8,.5) {$e_v$};
			\node at (.8,.5) {$f_v$};
			\node at (-1.5,1.4) {$e_u$};
			\node at (-.5,1.4) {$f_u$};
			\node at (-1,0) {$B_\ell$};
			
			\begin{scope}[shift ={(6,0)}]
				\tikzstyle{every node}=[inner sep=1.5pt, fill=black,circle,draw]
				\node (v) at (0,0) {};
				\node (u) at (-2,2) {};
				\node (w) at (2,2) {};
				\node (x) at (-4, 4) {};
				\node (y) at (0,4) {};
				\draw (v) to (u);
				\draw (v) to (w);
				\draw (u) to (x);
				\draw (u) to (y);
				\draw[dotted, thick] (v) to (0,-1/2);
				
				\draw (-4,5) ellipse (.8 and 1.2);
				\draw (0,5) ellipse (.8 and 1.2);
				\draw (2,3) ellipse (.8 and 1.2);
				
				\tikzstyle{every node}=[]
				
				\draw (v) node[left] {$v$};
				\draw (u) node[left] {$u$};
				\node at (2,3) {$\phi(Y_v)$};
				\node at (-4,5) {$\phi(X_u)$};
				\node at (0,5) {$\phi(Y_u)$};
				\node at (0,.5) {$\eta(v)$};
				\node at (-2,2.5) {$\eta(u)$};
				\node at (-.9,1.8) {$\eta(e_v)$};
				\node at (.9,1.75) {$\eta(f_v)$};
				
				\node at (-1.1,3.75) {$\eta(f_u)$};
				\node at (-2.9,3.8) {$\eta(e_u)$};
				\node at (-2,0) {$\phi(B_\ell)\subseteq G$};
				\draw[-Stealth] (-.4,.55) to (-1.05,1.2);
				\draw[-Stealth] (.4,.55) to (1.05,1.2);
				\draw[-Stealth] (-2.4,2.55) to (-3.05,3.2);
				\draw[-Stealth] (-1.6,2.55) to (-.95,3.2);
				
				\draw [very thick, decorate,
				decoration = {calligraphic brace,
					amplitude=5pt}] (-1.95,1.95) -- (-1.1,1.1);
				\draw [very thick, decorate,
				decoration = {calligraphic brace,
					amplitude=5pt}] (-0.9,3.1) -- (-0.05,3.95);
				\draw [very thick, decorate,
				decoration = {calligraphic brace,
					amplitude=5pt}] (-3.95,3.95) -- (-3.1,3.1);
				\draw [very thick, decorate,
				decoration = {calligraphic brace,
					amplitude=5pt}] (1.1,1.1) -- (1.95,1.95);

			\end{scope}
		\end{tikzpicture}
		\caption{For the proof of \ref{gettree}.} \label{fig:gettreecase1}
	\end{figure}

	Thus,
	we may assume that either $x$ belongs to $X_v$ or $x=e_v$, and either $y$ belongs to $Y_v$
	or $y=f_v$. If $x=e_v$ and $y=f_v$, then since $\phi(e_v)\cup \phi(f_v)$ is a geodesic, and
	$$\dist_G(\eta(e_v),v)=\dist_G(\eta(f_v),v) = 2c+1,$$
	it follows that
	$\dist_G(\eta(x),\eta(y))\ge 4c+2>c$ as required. Thus, we may assume that $x$ belongs to $X_v$, and hence $\eta(x)$ is a subgraph
	of $\phi(X_v)$. If $y=f_v$, then since $\phi(e_v)\cup \phi(f_v)$ is a geodesic between $\phi(X_v), \phi(Y_v)$, there is a geodesic
	between $\eta(y), \phi(X_v)$ that contains $\phi(e_v)$, and therefore has length at least $d+1$; and so
	$\dist_G(\eta(x),\eta(y))\ge d+1>c$ as required. Finally, if $y$ belongs to $Y_v$, then $\eta(y)$ is a subgraph of $\phi(Y_v)$,
	and so
	$$\dist_G(\eta(x), \eta(y))\ge \dist_G(\phi(X_v),\phi(Y_v))\ge 2(d+1) >c,$$
	as required. Consequently, $\eta$ exhibits $B_\ell$ as a $c$-fat minor. This proves \ref{gettree}.~\bbox

	\begin{thm}\label{grouping}
		Let $\ell\ge 1$ and $d\ge 3c\ge 6$, and suppose that $G$ does not contain $H_\ell$ as a $c$-fat minor. Let $\mac A$ be a set of
		vertex-disjoint connected subgraphs of $G$. Then there is a set $\mac P$ of paths of $G$, with union $U\mac P$ say, such that:
		\begin{itemize}
			\item each member of $\mac P$ has length at most $4^{\ell+1} d$;
			\item for each $P\in \mac P$ and each $v\in V(P)$, there exists $A\in \mac A$ such that there is a path in $U\mac P$ of length
			$<4^{\ell+1} d$ between $v$ and $A$; 
			\item let $J$ be the subgraph consisting of the union of all members of $\mac A$ and all members of $\mac P$;
			then $\dist_G(u,v)>d$ for all vertices $u,v$ in distinct components of $J$.
		\end{itemize}
	\end{thm}
	\Proof
	Let $J_0$ be the union of all the subgraphs in $\mac A$. Inductively, we will define $\mac P_i$ and $J_i$ for $1\le i\le \ell+2$ as follows.
	Suppose that $i\ge 1$ and $J_{i-1}$ has been defined. Let $\mac P_i$ be the set of all paths $P$ with the following properties:
	\begin{itemize}
		\item there exist distinct components $X,Y$ of $J_{i-1}$, such that $P$ is a shortest path of $G$ between $X,Y$; and
		\item $P$ has length at most $4^{\ell+2-i}d$.
	\end{itemize}
	Let $J_i$ be the union of $J_{i-1}$ and all the members of $\mac P_i$. This completes the inductive definition of $J_i$ for $i \ge 1$.
	\\
	\\
	(1) {\em For each $i\ge 1$, let $P\in \mac P_i$, joining components $X,Y$ of $J_{i-1}$, and let the ends of $P$ be $x\in V(X)$ and
		$y\in V(Y)$. Then there is a geodesic $(\ge d)$-subdivision of $B_{i-1}$ in $G$ that is a subgraph of $X$, with root $x$, and the same for $Y,y$.}
	\\
	\\
	We proceed by induction on $i$. The result is clear for $i = 1$ (taking a ``geodesic $(\ge d)$-subdivision of $B_0$'' to mean a rooted
	tree with only one vertex), so we assume that $i\ge 2$ and the result holds for $i-1$. Let $P\in \mac P_i$, joining components
	$X,Y$ of $J_{i-1}$, and let the ends of $P$ be $x\in V(X)$ and $y\in V(Y)$. (See Figure \ref{fig:grouping}.) Thus, $x,y\in V(J_{i-1})$.
	Either $y\in V(J_{i-2})$, or $y$ lies in a path
	in $\mac P_{i-1}$ included in $Y$ with both ends in $V(J_{i-2})$, and in either case there is a path of $Y$ between
	$y,J_{i-2}$ of length at most $4^{\ell+2-i+1}d/2$. By combining this with $P$, we deduce that there is a path between $x$ and
	a vertex in $J_{i-2}\cap Y$ of length at most $4^{\ell+2-i+1}d/2+4^{\ell+2-i}d\le 3\cdot 4^{\ell+2-i}d$.
	
	\begin{figure}[H]
		\centering
		\begin{tikzpicture}[yscale=.6,auto=left]
			\draw (0,0) ellipse (2 and 3);
			\draw (7,0) ellipse (2 and 3);
			\draw (-.5,1.5) ellipse (.5 and .5);
			\draw (-.5,-1.5) ellipse (.5 and .5);

			\tikzstyle{every node}=[inner sep=1.5pt, fill=black,circle,draw]
			\node (x) at (1.7,0){};
			\node (y) at (5.3,0){};
			\node (x') at (-.35,1.1){};
			\node (y') at (-.35,-1.1){};
			
			\draw[dashed, thick] (x) to (y);
			\draw[dashed, thick] (x) to (x');
			\draw[dashed, thick] (x) to (y');
			\draw[dashed, thick] (y) to (7,0);
			\draw[] (-.5,1.5) ellipse (.5 and .5);
			\draw[] (-.5,-1.5) ellipse (.5 and .5);
			\draw[fill = white] (7.4,0) ellipse (1 and 1);
			
			\tikzstyle{every node}=[]
			\draw node at (3.5,.3) {$P$};
			\draw (x) node[above] {$x$};
			\draw (y) node[above] {$y$};
			\draw node at (1,1.5) {$X$};
			\draw node at (6,1.5) {$Y$};
			\draw node at (7.4,0) {$J_{i-2}\cap Y$};
			\draw node at (-.5,1.5) {$X'$};
			\draw node at (-.5,-1.5) {$Y'$};
			\draw (x') node[below] {$x'$};
			\draw (y') node[above] {$y'$};

		\end{tikzpicture}
		\caption{For the proof of \ref{grouping}.} \label{fig:grouping}
	\end{figure}
	
	But every path between vertices of $J_{i-2}$ that belong to distinct components of $J_{i-1}$ has length more than $4^{\ell+2-i+1}d$,
	so every path of $X$ between $x$ and $J_{i-2}$ has length
	more than $4^{\ell+2-i}d$. In particular, $x\notin V(J_{i-2})$, and so
	there is a path
	$Q\in \mac P_{i-1}$ with $x$ in its interior,
	and the two subpaths of $Q$ between $x$ and the ends of
	$Q$ both have length more than $4^{\ell+2-i}d$. Let $Q$ join components $X', Y'$ of $J_{i-2}$, and let
	the corresponding ends of $Q$ be $x', y'$.
	From the inductive hypothesis, there is a geodesic $(\ge d)$-subdivision of $B_{i-2}$ in $G$ that is a subgraph of $X'$, with root $x'$,
	and the same for $Y',y'$. But $Q$ is a geodesic between $X',Y'$, and hence the union of $Q$ and
	these two geodesic $(\ge d)$-subdivisions of $B_{i-2}$ makes a geodesic $(\ge d)$-subdivision of $B_{i-1}$ with root $x$. This proves (1).
	\bigskip
	
	Since $G$ does not contain $H_{\ell}$ as a $c$-fat minor, and therefore does not contain $B_{\ell+1}$ as a $c$-fat minor, it follows from
	\ref{gettree} that no subgraph of $G$ is a geodesic $(\ge d)$-subdivision of $B_{\ell+1}$. We deduce
	from (1) that $\mac P_{\ell+2}=\emptyset$. Consequently, every two components of $J_{\ell+1}$ have distance more than
	$4^{\ell+2-(\ell+2)}d=d$. Let $J=J_{\ell+1}$, and let $\mac P = \mac P_1\cupcup \mac P_{\ell+1}$. Then every path in $\mac P$ has
	length at most $4^{\ell+2-1}d=4^{\ell+1}d$. Moreover, for $1\le i\le \ell+1$, every vertex in $J_i$ either belongs to $J_{i-1}$ or to a path
	in $\mac P_i$, and so in either case has distance in $U\mac P$ at most $4^{\ell+2-i}d/2$ from some vertex in $J_{i-1}$; and so every 
	vertex in $J$ has distance from $J_0$ in $U\mac P$ 
	at most
	$$\sum_{1\le i\le \ell+1} 4^{\ell+2-i}d/2< 4^{\ell+1} d.$$
	This proves \ref{grouping}.~\bbox
	\section{Using the grouping lemma}

	We apply the grouping lemma to ``civilize'' a given society as follows. 
	\begin{thm}\label{newbestsoc}
		Let $k\ge 0$ be a century, and let $\tau=(d,\alpha,\beta)$ and $\tau'=(d', \alpha', \beta')$ be canons,
		where 
		\begin{align*}
			d'&\ge 3c\\
			d&\le d_0\\
			d&\ge 4^{\ell+2} d'\\
			\alpha'&= 2\alpha\\
			\beta' &\ge  \beta +2d_0+d/4.
		\end{align*}
		Let $G$ be a graph that does not contain $H_\ell$ as a $c$-fat minor, and such that there is a $k$th-century
		$\tau$-society in $G$. Then there is a civilized
		$k$th-century $\tau'$-society in $G$.
	\end{thm}
	\Proof
	Let $\mac T$ be a $k$th-century $\tau$-society in $G$, and let $W$ be the set of vertices $v\in V(G)$ such that 
	$\dist_G(v,\bd(\voro_{\mac T}(C)))\ge d_0$
	for some $\mac T$-village $C$ with $v\in \voro_{\mac T}(C)$.
	Let $Q$ be the union of $W$ and all houses of $\mac T$, and let $\mac A$ be the set of all
	vertex sets of components of $G[Q]$. Let $\mac T''$ be the union of $\mac A$ and the set of all forts of $\mac T$. It follows that all 
	members of $\mac T''$ are pairwise disjoint, and so $\mac T''$ defines Voronoi cells.
	\\
	\\
	(1) {\em For each $A\in \mac A$ there is a $\mac T$-village $C$ such that $\voro_{\mac T''}(A)\subseteq \voro_{\mac T}(C)$.}
	\\
	\\
	For each $u\in A$ there is a unique $\mac T$-village $C$ with $u\in \voro_{\mac T}(C)$. Since there are no edges between $\voro_{\mac T}(C)$ and $\voro_{\mac T}(C')$
	for distinct $\mac T$-villages $C,C'$, and $G[A]$ is connected, it follows that there is a $\mac T$-village $C$ such that
	$A\subseteq \voro_{\mac T}(C)$.
	We claim that $\voro_{\mac T''}(A)\subseteq \voro_{\mac T}(C)$. Suppose not, and choose 
	$w\in \voro_{\mac T''}(A)\setminus \voro_{\mac T}(C)$. Since $w\in \voro_{\mac T''}(A)$, the $\Lambda$-shortest path $P$ between 
	$w, V(\mac T'')$ has an end $s\in A$, and $V(P)\subseteq \voro_{\mac T''}(A)$. Since $w\notin  \voro_{\mac T}(C)$ and $s\in A\subseteq \voro_{\mac T}(C)$, 
	there is an edge $uv$ of $P$ such that $u\in \voro_{\mac T}(C)$ and $v\notin \voro_{\mac T}(C)$. The subpath of $P$ between $v,s$
	is the $\Lambda$-shortest path between $v, V(\mac T'')$. Since $V(\mac T)\subseteq V(\mac T'')$, it follows that if 
	$s\in V(C)\subseteq V(\mac T)$ then this subpath is also the $\Lambda$-shortest path between $v, V(\mac T)$, contradicting that $v\notin \voro_{\mac T}(C)$.
	Thus $s\in W$, and hence $\dist_G(s,u)\ge d_0$, and so $\dist_G(s,v)>d_0$, which is impossible since 
	$v\in V(P)\subseteq \voro_{\mac T''}(A)$, and $\dist_G(s,v) = \dist_G(A,v)$. This proves that $\voro_{\mac T''}(A)\subseteq \voro_{\mac T}(C)$, and so proves (1).
	\\
	\\
	(2) {\em If $A\in \mac A$ and $Y$ is a fort of $\mac T$ then $\dist_G(A,Y)>d$.}
	\\
	\\
	Let $P$ be a shortest path in
	$G$ between $A,Y$, with ends $x\in A$ and $y\in Y$.
	If $x\notin W$, then $x\in X$ for some house $X$ of $\mac T$, and then $|E(P)|\ge \dist_G(X,Y)>d$ as required. 
	If $x\in W$, choose a $\mac T$-village $C$ such that $x\in \voro_{\mac T}(C)$ and
	$\dist_G(x,\bd(\voro_{\mac T}(C)))\ge d_0$. There is a vertex $u\in V(P)\setminus Y$ with $u\in \bd (\voro_{\mac T}(C))$,
	and consequently 
	$$d\le d_0\le \dist_G(x,u)<|E(P)|=\dist_G(A,Y)$$
	as required. This proves (2).
	
	\bigskip
	
	We designate the members of $\mac A$ as houses of $\mac T''$, and the forts of $\mac T$ as forts of $\mac T''$. We do not claim that $\mac T''$ is a society (for any choice of canon)
	since its houses might be very close together. Nevertheless, let us use the same terminology: a {\em $\mac T''$-community} means a $\mac T''$-adjoin-connected set of its houses, and a {\em $\mac T''$-village} is a maximal $\mac T''$-community.
	\\
	\\
	(3) {\em For every $\mac T''$-village $C''$, 
		$\voro_{\mac T''}(C'')$ has quasi-line-width at most $(\alpha,\beta)$ and
		$\bigcup_{A\in C''}\bd (A)$ has quasi-size at most $(\alpha,\beta+d_0)$.}
	\\
	\\
	Let $A_1,A_2\in C''$ such that $A_1$ $\mac T''$-adjoins $A_2$.
	By (1), there are $\mac T$-villages $C_1,C_2$ such that $\voro_{\mac T''}(A_i)\subseteq \voro_{\mac T}(C_i)$ for $i = 1,2$.
	Since there is an edge between $\voro_{\mac T''}(A_1),\voro_{\mac T''}(A_2)$, and there is no edge between
	$\voro_{\mac T}(C_1),\voro_{\mac T}(C_2)$ if $C_1\ne C_2$, it follows that $C_1=C_2$; and since $C''$ is $\mac T''$-adjoin-connected, we deduce 
	that there is a $\mac T$-village $C$ such that  $\voro_{\mac T''}(C'')\subseteq \voro_{\mac T}(C)$.
	
	Now 
	$\voro_{\mac T}(C)$ has quasi-bound at most $(\alpha,\beta)$, since $\mac T$ is a $\tau$-society. 
	Since $\voro_{\mac T''}(C'')\subseteq \voro_{\mac T}(C)$, it follows that $\voro_{\mac T''}(C'')$, and therefore $V(C'')$, has quasi-line-width
	at most $(\alpha,\beta)$. It remains to check that $\bigcup_{A\in C''}\bd (A)$ has quasi-size
	at most $(\alpha,\beta+d_0)$. Let $u\in \bigcup_{A\in C''}\bd (A)$; then $u\in \bd(A)$ for some $A\in C''$. We claim that 
	$\dist_G(u,\bd(\voro_{\mac T}(C)))\le  d_0$. 
	To see this, choose
	$v\notin A$ adjacent to $u$ (this exists since $u\in \bd(A)$). Hence $v\notin Q$, since $A$ is the vertex set of a component of $G[Q]$.
	If $v\notin \voro_{\mac T}(C)$, then $u\in \bd(\voro_{\mac T}(C))$ and so $\dist_G(u,\bd(\voro_{\mac T}(C)))=0\le  d_0$ as claimed, so
	we assume that
	$v\in \voro_{\mac T}(C)$. Since $v\notin Q$, it follows that 
	$\dist_G(v,\bd(\voro_{\mac T}(C)))<d_0$. So $\dist_G(u,\bd(\voro_{\mac T}(C)))\le d_0$, as claimed. Since
	$\bd(\voro_{\mac T}(C))$ has quasi-size at most $(\alpha,\beta)$, it follows that  $\bigcup_{A\in C''}\bd (A)$ has quasi-size at most
	$(\alpha,\beta+d_0)$. This proves (3).
	
	\bigskip
	
	Let $\theta= 4^{\ell+1} d'$. Thus, $d\ge4\theta$.
	By \ref{grouping} applied to $\mac A$,
	there is a set $\mac P$ of paths of $G$, with union $U\mac P$ say, such that:
	\begin{itemize}
		\item for each $P\in \mac P$ and each $v\in V(P)$, there exists $A\in \mac A$ such that $\dist_{U\mac P}(v,A)<  \theta$;
		\item let $J$ be the subgraph of $G$ consisting of the union of all members of $\mac A$ and all members of $\mac P$;
		then $\dist_G(u,v)>d'$ for all vertices $u,v$ in distinct components of $J$.
	\end{itemize}
	Let $\mac S$ be the set of all vertex sets of components of $J$. 
	\\
	\\
	(4) {\em For each $S\in \mac S$, there is a $\mac T''$-village $C''$ such that
		every house of $\mac T''$ included in $S$ belongs to $C''$,
		and consequently $S\subseteq \voro_{\mac T''}(C'')$.}
	\\
	\\
	The distance between distinct houses in $\mac T''$ may be as small as two; but if two houses $A_1,A_2$ of $\mac T''$ are at distance 
	less than $d+1$ then they belong to the same $\mac T''$-village. To see this, suppose that $A_1,A_2$ do not belong to the same 
	$\mac T''$-village, and let $P$ be a shortest path between $A_1,A_2$.  Since there are no edges between $\voro_{\mac T''}(C_1)$ and $\voro_{\mac T''}(C_2)$
	for distinct $\mac T''$-villages $C_1,C_2$, some vertex $w$ of $P$ belongs to $\voro_{\mac T''}(Y)$ for some fort $Y$ of $\mac T''$. 
	So both subpaths of $P$ between $w$ and its ends have subpaths between 
	$\voro_{\mac T''}(Y)$ and one of $A_1, A_2$, and so both have length at least $(d+1)/2$ since $\dist_G(Y,A_i)\ge d+1$ for $i = 1,2$. 
	Hence $P$ has length at least $d+1$, a contradiction.
	
	Now let $S\in \mac S$. 
	For each $v\in S$ there exists $A\in \mac A$ such that $\dist_{G}(v, A)<\theta$; let $\mac A_v$ be the set of all such $A$. Since
	$\dist_G(A,A')\le 2(\theta-1)\le d$ for all distinct $A,A'\in \mac A_v$, it follows that they all belong to the same 
	$\mac T''$-village $C''_v$, and if $v,v'\in S$ are adjacent, then any member of $\mac A_v$ has distance at most $2(\theta-1)+1$ from each member of $\mac A_{v'}$,
	and so $C''_v=C''_{v'}$. Consequently there is a $\mac T''$-village $C''$ containing $A$ for all $A\in \mac A$ such that 
	$\dist_{G}(S, A)<\theta$.  
	
	It remains to show that $S\subseteq \voro_{\mac T''}(C'')$. Let $v\in S$, choose $A\in \mac A$ with $\dist_G(v,A)<\theta$, and choose $Y\in \mac T''$ with $v\in \voro_{\mac T''}(Y)$. Thus, $\dist_G(v,Y)\le \dist_G(v,A)$. Since $\dist_G(Y,A)\le 2(\theta-1)\le d$,
	and every
	fort of $\mac T''$ has distance more than $d$ from $A$ (by (2)), it follows that $Y$ is not a fort, and so $Y\in \mac A$. Since
	$$\dist_G(S,Y)\le \dist_G(v,Y)\le \dist_G(v,A)<\theta$$ 
	it follows that $Y\in C''$ and hence $v\in  \voro_{\mac T''}(C'')$. This proves that $S\subseteq \voro_{\mac T''}(C'')$, and so proves (4).

	\bigskip
	
	For each fort $Y\in \mac T''$, let $F(Y)$ be the set of all 
	$v\in V(G)$ such that $\dist_G(v, Y)\le \theta$. Let $\mac F$ be the set of all sets $F(Y)$ for all forts $Y$ of $\mac T''$.
	Let $\mac T'= \mac S\cup \mac F$.
	\\
	\\
	(5) {\em $\dist_G(X',Y')> d'$ for all distinct $X',Y'\in \mac T'$.}
	\\
	\\
	We have seen this already if $X',Y'\in \mac S$, from the choice of $\mac S$.
	If $X'\in \mac S$ and $Y'\in \mac F$, let $v\in X'$, choose a house $X$ of $\mac T''$ with $\dist_G(v,X)<\theta$, and let $Y'=F(Y)$
	for a fort $Y$ of $\mac T''$. Then 
	$$\dist_G(v,Y')\ge \dist_G(X,Y)-2\theta+1> d-2\theta+1>d',$$ as required. 
	Similarly, 
	since every two 
	forts of $\mac T''$ are at distance more than $d$, it follows that every two members of $\mac F$ are 
	at distance more than $d-2\theta\ge d'$. This proves (5).
	
	\bigskip
	
	In particular, the buildings in $\mac T'$ are pairwise disjoint, and we may speak of $\voro_{\mac T'}(X)$ for $X\in \mac T'$. 
	Assign the members of $\mac S$ to be houses of $\mac T'$, and the members of $\mac F$ to be forts of $\mac T'$.
	We will show that $\mac T'$ is a civilized
	$k$th-century $\tau'$-society, and here are some steps to show that:
	\\
	\\
	(6) {\em The following hold:
		\begin{itemize}
			\item For all $v\in V(G)$, there exists $X\in \mathcal{T}'$ such that $\dist_G(v,X)\le d_0$.
			\item Every fort of $\mac T'$ includes a $c$-superfat $H_k$-minor of $G$.
			\item Each fort of $\mac T'$ has quasi-bound at most $(\alpha', \beta')$.
		\end{itemize}
	}
	\noindent
	The first is true since it is true for $\mac T$ and each member of $\mac T$ is a subset of a member of $\mac T''$ and hence of $\mac T'$. The second is true
	because each fort of $\mac T'$ includes a fort of $\mac T''$ and hence one of $\mac T$.
	
	To see the third, let $Y$ be a fort of $\mac T''$, and hence of $\mac T$. Then $Y$ has quasi-bound at most $(\alpha,\beta)$. By \ref{increaserad},
	$\bd(F(Y))$ has quasi-size at most $(\alpha, \beta+\theta)$, and $F(Y)$
	has quasi-line-width at most $(2\alpha,\beta+\theta)$. 
	Since $2\alpha= \alpha'$ and $ \beta+\theta\le \beta'$, this proves the third statement, and hence proves (6).
	
	\bigskip
	
	For each $S\in \mac S$, let $C_S$ be the $\mac T''$-village that contains each house of $\mac T''$ included in $S$ (this exists, by (4)).
	\\
	\\
	(7) {\em For each $S\in \mac S$, $\voro_{\mac T'}(S)\subseteq \voro_{\mac T''}(C_S)$.}
	\\
	\\
	Suppose not; then since $S\subseteq \voro_{\mac T''}(C_S)$ by (4), there are adjacent $u,v\in \voro_{\mac T'}(S)$ such that 
	$u\in \voro_{\mac T''}(C_S)$ and $v\notin \voro_{\mac T''}(C_S)$. Choose $Y\in \mac T''$ with $v\in \voro_{\mac T''}(Y)$. Thus, 
	$\dist_G(v,C_S)\ge \dist_G(v,Y)$. 
	Since $v\notin \voro_{\mac T''}(C_S)$, it follows that $Y\notin C_S$, and since no house of $\mac T''$ $\mac T''$-adjoins $C_S$ (because $C_S$
	is a $\mac T''$-village), it follows
	that $Y$ is a fort of $\mac T''$. (See Figure \ref{fig:bestsoc}.)
	
	\begin{figure}[H]
		\centering
		
		\begin{tikzpicture}[scale=1,auto=left]
			\draw[dotted, thick] (-1.4,0) ellipse (4 and 2);
			\draw[dotted, thick] (.8,0) ellipse (3.6 and 1.7);
			\draw (0,0) ellipse (1 and 1);
			\draw (7,0) ellipse (.5 and .5);
			\draw (7,0) ellipse (1.2 and 1.2);
			\draw[dotted, thick] (6,0) ellipse (3.2 and 1.5);
			\draw (0,1/2) ellipse (.3 and .3);
			
			\tikzstyle{every node}=[inner sep=1.5pt, fill=black,circle,draw]
			\node (s) at (.8,0) {};
			\node (u) at (2.4,0) {};
			\node (v) at (3,0) {};
			
			\draw[dashed, very thick] (s) to (u);
			\draw[thick] (u) to (v);
			\draw[dashed, very thick] (v) to (7,0);
			\draw[dashed, very thick] (s) to (0,.5);
			
			\draw[fill = white] (7,0) ellipse (.5 and .5);
			\draw[fill = white] (0,1/2) ellipse (.3 and .3);
			
			\tikzstyle{every node}=[]
			\node at (0,.5) {$X$};
			\node at (7,0) {$Y$};
			\node at (0,0) {$S$};
			\node at (-4,0) {$\voro_{\mac T''}(C_S)$};
			\node at (2.8,1) {$\voro_{\mac T'}(S)$};
			\node at (5,1) {$\voro_{\mac T''}(Y)$};
			\node at (7,.8) {$F(Y)$};
			\draw (u) node[above] {$u$};
			\draw (v) node[above] {$v$};
			\draw (s) node[below left ] {$s$};
			
		\end{tikzpicture}
		\caption{For the proof of (7).} \label{fig:bestsoc}
	\end{figure}
	
	Now $v\notin F(Y)$, since $v\in \voro_{\mac T'}(S)$; and so 
	$\dist_G(v,F(Y))\le\dist_G(v,Y)-\theta$. Since $v\in \voro_{\mac T'}(S)$ and $F(Y)\in \mac T'$, it follows that 
	$$\dist_G(v,S)\le \dist_G(v,F(Y))\le\dist_G(v,Y)-\theta\le \dist_G(v,C_S)-\theta.$$
	Choose $s\in S$ such that $\dist_G(v,S) = \dist_G(v,s)$; and choose a house $X$ of $\mac T''$ with $\dist_{G[S]}(s,X)<\theta$. 
	Then $X\subseteq S$, and so $X\in C_S$ by (4). Consequently, 
	$$\dist_G(v,C_S)\le \dist_G(v,X)< \dist_G(v,s) + \theta,$$
	a contradiction. This proves (7).
	\\
	\\
	(8) {\em For every $\mac T'$-community $C'$, $\voro_{\mac T'}(C')$ has quasi-line-width at most $(\alpha,\beta)$;
		$\bigcup_{S\in C'}\bd (S)$ has quasi-size at most $(\alpha,\beta+d_0+\theta)$; and
		$\bd(\voro_{\mac T'}(C'))$ has quasi-size at most $(\alpha, \beta+2d_0+\theta)$.}
	\\
	\\
	From (7), it follows that if $S_1,S_2\in C'$ $\mac T'$-adjoin, then $C_{S_1} = C_{S_2}$, and since $C'$ is $\mac T'$-adjoin-connected, 
	all the $\mac T''$-villages $C_S \;(S\in C')$ are equal. 
	Hence there is a $\mac T''$-village $C''$ such that $\voro_{\mac T'}(C')\subseteq \voro_{\mac T''}(C'')$. By (3), 
	$\voro_{\mac T''}(C'')$ has quasi-line-width at most $(\alpha,\beta)$, and since 
	$\voro_{\mac T'}(C')\subseteq \voro_{\mac T''}(C'')$,
	it follows that $\voro_{\mac T'}(C')$ has quasi-line-width at most $(\alpha,\beta)$. 
	Now every vertex in $\bigcup_{S\in C'}\bd(S)$
	has distance at most $\theta$ from $\bigcup_{X\in C''}\bd (X)$, and since the latter has quasi-size at most $(\alpha,\beta+d_0)$
	by (3), it follows that $\bigcup_{S\in C'}\bd (S)$ has quasi-size at most $(\alpha,\beta+d_0+\theta)$. 
	Since every vertex in $\voro_{\mac T'}(C')\setminus V(C')$ has distance at most $d_0$ from $\bigcup_{S\in C'}\bd (S)$,
	\ref{increaserad} implies that
	$\bd(\voro_{\mac T'}(C'))$ has quasi-size at most $(\alpha, \beta+2d_0+\theta)$. This proves (8).
	
	\bigskip
	
	From (6) and (8), $\mac T'$ is a civilized
	$k$th-century $\tau'$-society. This proves \ref{newbestsoc}.~\bbox

	\section{Realms}

	That concludes the first part of the proof.
	Now we come to the second part of the paper, and here our goal is to show that if $G$ admits a civilized $k$th-century
	$\tau$-society then it admits a $(k+1)$st-century $\tau'$-society, for suitable $\tau, \tau'$. The proof uses ``realms'', which
	are basically civilized societies with castles added; and 
	we need a number of further definitions. 
	
	For $\gamma\ge 0$ an integer, 
	a {\em spacing with unit $\gamma$} is the function $\delta$ from $\{0\LL 2\ell\}$ to $\mathbb{R}$ defined by 
	$\delta(i)  = \gamma4^{-i}$ for $0\le i\le 2\ell$. (We will arrange that $\delta$ is integer-valued, 
	by choosing $\gamma$ carefully.)  Thus $\delta(i)= 4\delta(i+1)$ for each $i$ with $0\le i<2\ell$.
	A {\em spacing} is a spacing with unit $\gamma$ for some choice of $\gamma$.
	A {\em standard} is a triple $\sigma= (\delta, \alpha,\beta)$, where
	$\alpha\ge 12$ and $\beta\ge 1$ are integers, and
	$\delta$ is a spacing with $\delta(2\ell)\ge 5c$ and $\delta(0)\le d_0$. (The condition $\alpha\ge 12$ is just to make some arithmetic come out more smoothly.)
	
	
	Let $0\le k<\ell$ be a century and let $\sigma=(\delta, \alpha,\beta)$ be as above. We say a {\em $k$th-century realm} with {\em standard $\sigma$} in a graph $G$ is a set $\mac T$ of pairwise vertex-disjoint buildings of $G$, 
	where each member of $\mac T$ is assigned to be a house, fort or castle of $\mac T$, with no houses if $k=0$, 
	satisfying the following (where $\rk(X)=k-1,k$ or $k+1$ depending whether $X$ is a house, fort or castle):
	\begin{itemize}
		\item For all $v\in V(G)$, there exists $X\in \mathcal{T}$ such that $\dist_G(v,X)\le d_0$;
		\item If $X,Y\in \mac T$ are distinct, then $\dist_G(X,Y)>\delta(\rk(X)+\rk(Y))$; 
		\item Every fort or castle $X$ of $\mac T$ includes a $c$-superfat $H_{\rk(X)}$-minor of $G$;
		\item Every fort of $\mac T$ has quasi-bound at most $(\alpha, \beta)$; 
		\item Every castle of $\mac T$ has quasi-bound at most $(9\alpha, \beta+2d_0)$;
		\item For every $\mac T$-community $C$, $V(C)$ has quasi-line-width at most $(\alpha,\beta)$ and 
		$\bigcup_{X\in C}\bd (X)$ has quasi-size at most $(\alpha,\beta)$.
		(Again, a {\em $\mac T$-community} of $\mac T$ is an adjoin-connected set of houses of $\mac T$, and a {\em $\mac T$-village} is a maximal $\mac T$-community.)
	\end{itemize}
	If $X$ belongs to a realm $\mac T$, then it is a house, fort or castle of $\mac T$, and we call this its {\em class} under $\mac T$.
	It is easy to turn a civilized society into a realm for the same century (note, however, that we use the full strength of 
	``civilized'' -- this was the point of converting a society into a civilized one):
	\begin{thm}\label{bestsoctorealm}
		Let $0\le k<\ell$ be a century, and let $\sigma=(\delta,\alpha,\beta)$ be a standard. 
		If $k=0$, every civilized $k$th-century $(\delta(0),\alpha,\beta)$-society $\mac T$ in $G$ is also a $k$th-century realm with standard $\sigma$.
		If $k>0$, every civilized $k$th-century $(\delta(2k-2),\alpha,\beta)$-society $\mac T$ in $G$ is also a $k$th-century realm with standard $\sigma$.
	\end{thm}
	\Proof
	If $k=0$, and  $X,Y\in \mac T$ are distinct, then $\dist_G(X,Y)>\delta(0)$ since $\mac T$ is a $0$th-century 
	$(\delta(0),\alpha,\beta)$-society; and $\delta(0) \ge \delta(\rk(X)+\rk(Y))$ because $0$th-century societies
	have no houses.  If $k>0$, then similarly $\dist_G(X,Y)>\delta(2k-2)\ge \delta(\rk(X)+\rk(Y))$. So in either case, 
	$\dist_G(X,Y)> \delta(\rk(X)+\rk(Y))$.
	This proves \ref{bestsoctorealm}.~\bbox

	\section{A sketch of the rest of the proof}\label{sec:sketch}
	
	Let us give an overall view of the remainder of the proof, before we embark on its details. To reach the goal expressed at the start of the previous section, it suffices (in view of \ref{bestsoctorealm}) to be able  
	to convert a $k$th-century realm (of a reasonably high standard) to a $(k+1)$st-century $\tau$-society, where $\tau$
	is some suitable canon. 
	Suppose then that we have some $k$th-century realm $\mac T$.
	
	\subsection{Making it optimal}
	
	The first step is to ``optimize'' $\mac T$.
	This means, roughly, choosing $\mac T$ with the same standard and with as many castles as possible, combining houses and forts to make castles wherever we can.
	If some village adjoins three forts, 
	then we could try
	to combine the Voronoi closures of these houses and forts into a castle; or if there are three villages joining a given fort to three other forts,
	again we could hope to combine their closures into a castle. 
	(With some care: we need that the final approach to  each of the three important forts is via a geodesic of length $c+1$, so we can apply \ref{claw}). This is the content of \ref{newaddclaw}.
	If we succeed, we make a new realm, still in the same century and with the same standard.
	
	Suppose that, after a sequence of such promotions, we have converted $\mac T$ to a new realm $\mac T'$ (with the same century), and we want to show that 
	$\mac T'$ has the same standard.  
	Then in particular, we need that
	the boundary of the union of the Voronoi cells of each village 
	of $\mac T'$ has bounded  
	quasi-size. (Let us call this ``fact F''.) This is awkward to guarantee, because the Voronoi cells might change as we move to new realms, and therefore which sets are villages 
	might also change. But for any building of $\mac T$
	that remains in $\mac T'$, its Voronoi closure in $\mac T'$ is a subset of its closure in $\mac T$, so the villages of $\mac T'$ are 
	communities (not necessarily maximal) of $\mac T$. We don't know in advance which communities of $\mac T$ might end up as 
	villages in some $\mac T'$; so to guarantee that fact F will hold for $\mac T'$, 
	we will arrange that fact F holds for $\mac T$ in a strengthened, hereditary form. This is the reason for the very strong realm axiom, 
	that for every community $X$, the                                  
	union of the 
	boundaries of the members of $X$ has bounded quasi-size. (The point is that the union of the boundaries may be very different from 
	the boundary of the union.)
	
	When $\mac T$ is optimal, the set of its forts and villages has a sort of linear structure; each village adjoins 
	at most two forts, and for each fort, there are at most two other forts to which it is connected via a village. So each 
	component of the
	graph of villages and forts defined by adjoining pairs has a ``spine'', a path or cycle that contains all forts in the component,
	together
	with some villages each attaching to one or two consecutive members of the spine. (This is \ref{contactgraph}.) By exploiting this linear structure,
	we prove (in \ref{composeline} and \ref{smallcoarsepw}) that the 
	union of the Voronoi closures of all houses and forts has bounded quasi-line-width. This is a key step in our attempt to move to a new century; 
	we hope to promote all the current houses and forts to be houses of a $(k+1)$st-century realm, and at least we have control of the quasi-line-width of the union of their Voronoi cells. The problem is to bound the quasi-size of their boundaries.
	
	But we can get more than this from optimality of the realm. We might be able to combine some triple of forts into a castle even if they are 
	{\em not} adjoin-connected via villages. If we can connect together three forts in a claw, moving between them either 
	via villages or
	by paths (``passages'') that are not very long and do not go very close to other forts or castles, then we want to do so. This step is critical later
	when we introduce ``governments''.
	
	\subsection{Governments and revolution}
	
	Castles have rank $k+1$. 
	In the same way that we combined the Voronoi cells of three forts into a castle, we could try to combine the cells of three castles  into some 
	``super-castle'' of rank $k+2$, and so on; these new higher-level objects are called ``palaces'', and the set of them,
	together with the members of $\mac T$ not included in any palace, is called a ``government''. 
	But in a way this is simpler than combining forts into castles. There were several different scenarios for the latter; if some village adjoins three forts, or some fort
	can reach three others via villages, or if we can connect them via passages. For $j\ge k+1$, promoting three palaces of rank 
	$j$ 
	to a palace of rank $j+1$ is simpler: we look for an adjoin-connected set of houses and forts  that adjoins all three 
	palaces, and if we find one with the appropriate properties (a ``cabal'') we promote it all to make a new palace with rank $j+1$, and thereby change the government (a 
	``revolution''). (Except ``adjoin'' is the wrong term now, see below.)
	The more complex methods we used to make castles are still available at the higher levels, but seem not to help.
	
	If a house or fort has been used to glue some three palaces into a higher-rank palace, then it is no longer available
	for glueing other things together. We call the houses and forts that have not yet been used ``rebels''. 
	But we need to keep track of the buildings of the realm that have been combined into higher structures: the realm is not changing,
	but some of its buildings 
	are {\em also} being used as parts of palaces. 
	
	There is an extra complication: we have a partition into Voronoi cells defined by 
	the realm, and also a partition into Voronoi cells defined by the current government, and they are not the same. We use 
	``$\mac T$-adjoin'' and ``$\mac G$-adjoin'', to show which set of Voronoi cells we are referring to.
	
	There is a delicate issue here, one that gave us a lot of trouble. Suppose we have found a set $X$ of rebels, 
	$\mac G$-adjoin-connected, 
	and $X$ $\mac G$-adjoins three palaces of rank $j$ in the current government, and we want to promote all this to make a palace of rank 
	$j+1$. We already saw that the union of the $\mac T$-Voronoi cells of $X$ has bounded quasi-line-width, and therefore so does 
	the union of the $\mac G$-cells (which are subsets); but we also need that the boundary of the latter
	has bounded quasi-size.
	That was not a problem when we built castles, because the boundaries of villages and forts have bounded 
	quasi-size by definition,  and we were only using at most three villages and at most four forts; but now it becomes a big problem.
	
	We can get partway around it by making $X$ minimal subject to its $\mac G$-adjoining three palaces of the same rank, and therefore 
	the total number of palaces that $X$ $\mac G$-adjoins is bounded (unless $|X|=1$). So we would be happy if we could get a bound on the quasi-size of the 
	interface
	between the union of the $\mac G$-cells of $X$ and $\voro_{\mac G}(S)$ of each palace $S$ that $X$ $\mac G$-adjoins. 
	(We can get a bound
	on the quasi-size of the boundary of $\voro_{\mac G}(S)$, from the definition of a palace; but we 
	cannot use this to bound the interface, because 
	there is a vicious circle among the constants.) 
	
	This is where we use ``passages''. 
	Any term of $X$ that $\mac G$-adjoins $S$ is joined by a passage to one of the members of $\mac T$
	within $S$ (or it gives us no problem). But $X$ is $\mac G$-adjoin-connected, and so $\mac T$-adjoin-connected, and this gives a 
	linear structure as we saw earlier.  
	Interior terms of this linear structure already $\mac T$-adjoin (possibly via villages) two forts in $X$, and they cannot reach 
	three forts, even via a passage, from the optimality of $\mac T$. So, if an interior term $\mac G$-adjoins $S$, then its passage into $S$ 
	is very restricted. That gives us much more control of 
	the interface between $X$ and $S$, and eliminates the vicious circle. (This is \ref{palacebound} and \ref{cabalsaresmall}.)
	
	\subsection{Societies and the next century}
	
	If some $\mac G$-adjoin-connected set of rebels has $\mac G$-Voronoi closure with a boundary of large quasi-size, then it must 
	$\mac G$-adjoin many different palaces, and so three of the same rank,  and we
	can have a revolution. So when we finally reach a stable government, each $\mac G$-adjoin-connected set of rebels has closure with boundary of 
	small quasi-size. (This is \ref{notmuchleft}.)
	Consequently, we can make a $(k+1)$st-century society $\mac T'$ by declaring that each rebel is a house of $\mac T'$, and each 
	palace of the government is a fort of $\mac T'$; and that completes the proof.

	\section{Growing $c$-fat trees}
	
	We will try to grow a $c$-fat copy of $H_{t+1}$ by taking three $c$-fat copies of
	$H_{t}$, sufficiently far apart, and connecting them together appropriately. This goes more smoothly if we work with ``superfat'' rather than ``fat'', as the next result shows.
	
	\begin{thm}\label{claw}
		Let $c\ge 1$, let $t\ge 0$, and for $i = 1,2,3$, let $\eta_i$ exhibit $H_t$ as a $c$-superfat minor of a graph $G$,
		such that $\dist_G(\eta_i(H_t), \eta_j(H_t))> 5c+2$ for all distinct $i,j\in \{1,2,3\}$.
		Let $W\subseteq V(G)$ such that
		$G[W]$ is connected, and
		$\dist_G(\eta_i(H_t),W)= c+1$ for $i = 1,2,3$; and for $i = 1,2,3$, let $P_i$ be a geodesic from $W$ to $\eta_i(H_t)$.
		Then there is a mapping $\eta$ that exhibits $H_{t+1}$ as a $c$-superfat minor of $G$,
		such that
		$$\eta(H_{t+1})\subseteq G[W]\cup \eta_1(H_t)\cup \eta_2(H_t)\cup \eta_3(H_t)\cup P_1\cup P_2\cup P_3.$$
	\end{thm}
	\Proof (See Figure \ref{fig:superfat}.) For $1\le i<j\le 3$, $\dist_G(P_i,P_j)> (5c+2) -2(c+1)= 3c$, since $P_i,P_j$ both have length $c+1$, and they have ends at
	distance
	$>5c+2$. If $t=0$, let $\eta$ map the root of $H_{t+1}$ to $W$, its leaves to the three subgraphs $\eta_i(H_t)$ for
	$i = 1,2,3$, and each edge of $H_{t+1}$ to the interior of the corresponding geodesic $P_i$; then $\eta$ satisfies the
	theorem. Thus, we may assume that $t\ge 1$.
	
	\begin{figure}[H]
		\centering
		
		\begin{tikzpicture}[scale=1,auto=left]
			\draw (0,5) ellipse (5 and 1);
			
			\tikzstyle{every node}=[inner sep=1.5pt, fill=black,circle,draw]
			\node (w1) at (-4,4.7) {};
			\node (w2) at (0,4.4) {};
			\node (w3) at (4,4.7) {};

			\draw (-4.5,3) to (-6.2,3);
			\draw (-4,3) to (-5,1);
			\draw (-4,3) to (-3,1);
			\draw (0,3) to (-2.2,3);
			\draw (0,3) to (-1,1);
			\draw (0,3) to (1,1);
			\draw (4,3) to (1.8,3);
			\draw (4,3) to (3,1);
			\draw (4,3) to (5,1);
			
			\draw[fill = white]  (-4,3) ellipse (.6 and .5);
			\draw[fill = white]  (0,3) ellipse (.6 and .5);
			\draw[fill = white]  (4,3) ellipse (.6 and .5);
			
			\draw[fill = white]  (-6.2,3) ellipse (.6 and .5);
			\draw[fill = white]  (-2.2,3) ellipse (.6 and .5);
			\draw[fill = white]  (1.8,3) ellipse (.6 and .5);

			\draw[fill = white] (-5,1) ellipse (.6 and .5);
			\draw[fill = white]  (-3,1) ellipse (.6 and .5);
			\draw[fill = white]  (-1,1) ellipse (.6 and .5);
			\draw[fill = white]  (1,1) ellipse (.6 and .5);
			\draw[fill = white]  (3,1) ellipse (.6 and .5);
			\draw[fill = white]  (5,1) ellipse (.6 and .5);

			\node (u1) at (-5.9,3.25) {};
			\node (u2) at (-1.1,3) {};
			\node (u3) at (4,3.35) {};
			
			\draw (w1) -- (u1);
			
			\draw (w2) -- (u2);
			\draw (w3) -- (u3);
			\tikzstyle{every node}=[]
			
			\node at (0,5) {$W$};
			\node at (-5.3,4) {$P_1$};
			\node at (-1,3.6) {$P_2$};
			\node at (3.7,4) {$P_3$};
			\node at (-4,3) {$\eta_1(r)$};
			\node at (0,3) {$\eta_2(r)$};
			\node at (4,3) {$\eta_3(r)$};
			
			\node at (-6.2,3) {$\eta_1(B_3)$};
			\node at (-2.2,3) {$\eta_2(B_3)$};
			\node at (1.8,3) {$\eta_3(B_3)$};
			
			\node at (-5,1) {$\eta_1(B_1)$};
			\node at (-3,1) {$\eta_1(B_2)$};
			\node at (-1,1) {$\eta_2(B_1)$};
			\node at (1,1) {$\eta_2(B_2)$};
			\node at (3,1) {$\eta_3(B_1)$};
			\node at (5,1) {$\eta_3(B_2)$};
			
			\node at (-5.1,2.8) {$\eta_1(e_3)$};
			\node at (-1.1,2.8) {$\eta_2(e_3)$};
			\node at (2.9,2.8) {$\eta_3(e_3)$};
			
			\node at (-5.1,2) {$\eta_1(e_1)$};
			\node at (-2.9,2) {$\eta_1(e_2)$};
			\node at (-1.1,2) {$\eta_2(e_1)$};
			
			\node at (1.1,2) {$\eta_2(e_2)$};
			\node at (2.9,2) {$\eta_3(e_1)$};
			\node at (5.1,2) {$\eta_3(e_2)$};
		\end{tikzpicture}
		\caption{Growing a superfat minor. The figure shows the three ways that $P_h$ might attach.} \label{fig:superfat}
	\end{figure}
	
	Let $r$ be the root of $H_t$, let $B_1,B_2,B_3$ be the three components of $H_t\setminus \{r\}$, and for $h = 1,2,3$, let $e_h$ be the
	edge of $H_t$ between $r$ and $V(B_h)$. For $h,i\in \{1,2,3\}$, let $\eta_h(B_i^+)$ denote the subgraph
	$$\eta_h(e_i) \cup \bigcup_{x\in U(B_i)}\eta_h(x).$$
	We know that for $h = 1,2,3$, and all distinct $i,j\in \{1,2,3\}$,
	$$\dist_G(\eta_h(B_i^+), \eta_h(B_j^+))\ge 3c+1.$$
	Consequently, there is at most one value of $i\in \{1,2,3\}$ such that $\dist_G(P_h, \eta_h(B_i^+))\le c$, since
	$P_h$ has length $c+1$ and $\eta_h$ is $c$-superfat. By relabelling, we may assume that $\dist_G(P_h, \eta_h(B_i^+))>c$ for $i = 1,2$.
	
	Since $\dist_G(\eta_i(H_t), \eta_j(H_t))> 5c+2$ and $P_i,P_j$ have length $c+1$, it follows that
	$\dist_G(P_i, \eta_j(H_t))> 4c+1$, and $\dist_G(P_i, P_j)> 3c$, for all distinct $i,j\in \{1,2,3\}$. Moreover,
	$\dist_G(P_h, \eta_h(B_i^+))> c$ for $i = 1,2$ and $h = 1,2,3$.
	
	Let $r'$ be the root of $H_{t+1}$, and let $B_1',B_2',B_3', e_1', e_2', e_3'$
	be defined as usual. Thus, for $h = 1,2,3$,  $B_h'$ is isomorphic to $H_t\setminus V(B_3)$; let $\phi_h$ be such an isomorphism.
	Define $\eta(r') = G[W]$,
	and for $h = 1,2,3$, let $\eta(x) = \eta_h(\phi_h(x))$
	for each $x\in U(B_h')$. For $h = 1,2,3$, we define $\eta(e_h')$ as follows. If there is an end $y_h$ of $P_h$ in $\eta_h(r)$,
	let $\eta(e_h') = (P_h\setminus W)\setminus \{y_h\}$ (since $c\ge 1$, this subgraph is non-null).
	If not, then there is an end of $P_h$ in $\eta_h(B_3^+)$; let
	$\eta(e_h') = (P_h\setminus W)\cup \eta_h(B_3^+)$. Then
	$\eta$ exhibits $H_{t+1}$ as a $c$-superfat minor of $G$. This proves \ref{claw}.~\bbox
	
	
	\section{Optimizing within a century}
	
	We have a $k$th-century realm with some given standard, and we first want to choose an ``optimal'' one,
	by which we mean, roughly, one with
	the set of castles maximal;
	but since $G$ and the realm might be infinite, we need to formulate this
	carefully. If $0\le k<\ell$, we say a $k$th-century realm $\mac T_2$
	is an {\em extension} of a $k$th-century realm $\mac T_1$ if for each member $X_1\in \mac T_1$ there exists $X_2\in \mac T_2$ such that:
	\begin{itemize}
		\item $X_1\subseteq X_2$; and
		\item either $X_2$ is a building of higher class than $X_1$ (that is, either $X_1$ is a house of $\mac T_1$ and $X_2$ is a fort or castle of
		$\mac T_2$, or $X_1$ is a fort of $\mac T_1$ and $X_2$ is a castle of $\mac T_2$), or $X_1=X_2$ and $X_1$ has the same class under $\mac T_1$ and under $\mac T_2$.
	\end{itemize}
	Fix some standard $\sigma$.
	Let us say a  $k$th-century realm $\mac T$ with standard $\sigma$ is {\em optimal} (for $\sigma$) if no other  $k$th-century realm with the same
	standard
	is an extension of $\mac T$.
	It follows from Zorn's lemma that if there is a $k$th-century realm with standard $\sigma$, then one of its extensions is optimal.
	
	Next, an easy lemma. Let $F$ be a subset of vertices of a graph $H$. Let $P_1,P_2,P_3$ be paths of $H$ with a common end $z$ that are 
	otherwise pairwise vertex-disjoint. Let $P_i$ have ends $y_i, z$ for $i = 1,2,3$. We say the subgraph $P_1\cup P_2\cup P_3$
	is an {\em $F$-claw} of $H$ if for $i = 1,2,3$, $y_i\in F$, $y_i\ne z$, and no internal vertex of $P_i$ is in $F$. The vertex $z$ might belong to $F$. (See Figure \ref{fig:Fclaw}.)
	
	\begin{figure}[H]
		\begin{center}
			\begin{tikzpicture}[]
				\tikzstyle{every node}=[inner sep=1.5pt, fill=black,circle,draw]
				\node (y1) at (-1,0) {};
				\node (y2) at (0,0) {};
				\node (y3) at (1,0) {};
				\node (z) at (0,-2) {};
				\draw[dashed, thick] (y1) -- (z);
				\draw[thick, dashed] (y2) -- (z);
				\draw[thick, dashed] (y3) -- (z);
				\tikzstyle{every node}=[]
				\draw (y1) node[above] {$F$};
				\draw (y2) node[above] {$F$};
				\draw (y3) node[above] {$F$};
				\draw (z) node[below] {$F$?};

			\end{tikzpicture}
		\end{center}
		\caption{An $F$-claw. The dashed lines are paths.} \label{fig:Fclaw}
	\end{figure}
	
	\begin{thm}\label{noFclaw}
		Let $H$ be a connected graph and let $F\subseteq V(H)$, such that $H\setminus F$ has no edges.
		If no subgraph of $H$ is an $F$-claw,
		then there is an induced path or cycle of $H$ that contains all members of $F$, and every vertex of $H$ not in $F$ is adjacent only to one or two consecutive members of $F$. (See Figure \ref{fig:noFclaw}.)
	\end{thm}
	We omit the proof, which is easy.
	
	\begin{figure}[H]
		\begin{center}
			\begin{tikzpicture}[]
				\tikzstyle{every node}=[inner sep=1.5pt, fill=black,circle,draw]
				\node (y0) at (0,0) {};
				\node (y1) at (1,0) {};
				\node (y2) at (2,0) {};
				\node (y3) at (3,0) {};
				\node (y4) at (4,0) {};
				\node (y5) at (5,0) {};
				\node (y6) at (6,0) {};
				\node (y7) at (7,0) {};
				\node (y8) at (8,0) {};
				\node (y9) at (9,0) {};
				\node (y10) at (10,0) {};
				
				\tikzstyle{every node}=[inner sep=1.5pt, fill=white,circle,draw]
				\node (m1) at (1/2,.5) {};
				\node (m2) at (1/2,1) {};
				\node (m3) at (1.5,.5) {};
				\node (m4) at (2.5,0) {};
				\node (m5) at (2.5, .5) {};
				\node (m6) at (2.5, 1) {};
				\node (m7) at (4.5,.5) {};
				\node (m8) at (6.5,0) {};
				\node (m9) at (6.5,.5) {};
				\node (m10) at (6.5,1) {};
				\node (m11) at (7.5,.5) {};
				\node (m12) at (9.5,.5) {};
				
				\node (a1) at (.9,.7) {};
				\node (a2) at (1.3,.7) {};
				\node (a3) at (2,.7) {};
				\node (a4) at (4,.7) {};
				\node (a5) at (4.7,.7) {};
				\node (a6) at (5.3,.7) {};
				\node (a7) at (7,.7) {};
				\node (a8) at (8.7,.7) {};
				\node (a9) at (9,.7) {};
				\node (a10) at (9.3,.7) {};
				
				\foreach \from/\to in {y0/y1,y1/y2,y3/y4,y4/y5,y5/y6,y7/y8,y8/y9,y9/y10,
					m1/y0,m1/y1,m2/y0,m2/y1,m3/y1,m3/y2,m4/y2,m4/y3,m5/y2,m5/y3,m6/y2,m6/y3,m7/y4,m7/y5, m8/y6,m8/y7, m9/y6,m9/y7, m10/y6,m10/y7,m11/y7,m11/y8,m12/y9,m12/y10,
					a1/y1,a2/y1,a3/y2,a4/y4, a5/y5, a6/y5, a7/y7, a8/y9, a9/y9,a10/y9}
				\draw [-] (\from) -- (\to);

				\draw[dashed, thick] (-.7,0) -- (y0);
				\draw[dashed, thick] (y10) -- (10.7,0);
			\end{tikzpicture}
		\end{center}
		\caption{When there is no $F$-claw. The black vertices are in $F$. The horizontal path shown may be part of a cycle, or part of a finite or infinite path.} \label{fig:noFclaw}
	\end{figure}

	With $H,F$ as before, suppose that no two vertices in $V(H)\setminus F$ are adjacent. Then each of $P_1,P_2,P_3$ has length one or two, and they all 
	have length one unless their common end is in $F$. That gives us five cases (see Figure \ref{fig:stable}.)
	
	\begin{figure}[H]
		\begin{center}
			\begin{tikzpicture}[scale = .8]
				\tikzstyle{every node}=[inner sep=1.5pt, fill=black,circle,draw]
				
				\node (y1) at (-1,0) {};
				\node (y2) at (0,0) {};
				\node (y3) at (1,0) {};
				\tikzstyle{every node}=[inner sep=1.5pt, fill=white,circle,draw]
				\node (z) at (0,-2) {};
				\draw[] (y1) -- (z);
				\draw[] (y2) -- (z);
				\draw[] (y3) -- (z);
				
				\begin{scope}[shift ={(4,0)}]
					\tikzstyle{every node}=[inner sep=1.5pt, fill=black,circle,draw]
					
					\node (y1) at (-1,0) {};
					\node (y2) at (0,0) {};
					\node (y3) at (1,0) {};
					\node (z) at (0,-2) {};
					\draw[] (y1) -- (z);
					\draw[] (y2) -- (z);
					\draw[] (y3) -- (z);
				\end{scope}
				\begin{scope}[shift ={(8,0)}]
					\tikzstyle{every node}=[inner sep=1.5pt, fill=black,circle,draw]
					
					\node (y1) at (-1,0) {};
					\node (y2) at (0,0) {};
					\node (y3) at (1,0) {};
					\node (z) at (0,-2) {};
					\tikzstyle{every node}=[inner sep=1.5pt, fill=white,circle,draw]
					\node (w1) at (-1/2,-1) {};
					\draw (y1) -- (w1);
					\draw (w1) -- (z);
					\draw (y2) -- (z);
					\draw (y3) -- (z);
				\end{scope}
				\begin{scope}[shift ={(12,0)}]
					\tikzstyle{every node}=[inner sep=1.5pt, fill=black,circle,draw]
					
					\node (y1) at (-1,0) {};
					\node (y2) at (0,0) {};
					\node (y3) at (1,0) {};
					\node (z) at (0,-2) {};
					\tikzstyle{every node}=[inner sep=1.5pt, fill=white,circle,draw]
					\node (w1) at (-1/2,-1) {};
					\node (w2) at (0,-1) {};
					\draw (y1) -- (w1);
					\draw (w1) -- (z);
					\draw (y2) -- (w2);
					\draw (w2) -- (z);
					\draw (y3) -- (z);
					
				\end{scope}
				\begin{scope}[shift ={(16,0)}]
					\tikzstyle{every node}=[inner sep=1.5pt, fill=black,circle,draw]
					
					\node (y1) at (-1,0) {};
					\node (y2) at (0,0) {};
					\node (y3) at (1,0) {};
					\node (z) at (0,-2) {};
					\tikzstyle{every node}=[inner sep=1.5pt, fill=white,circle,draw]
					\node (w1) at (-1/2,-1) {};
					\node (w2) at (0,-1) {};
					\node (w3) at (1/2,-1) {};
					\draw (y1) -- (w1);
					\draw (w1) -- (z);
					\draw (y2) -- (w2);
					\draw (w2) -- (z);
					\draw (y3) -- (w3);
					\draw (w3) -- (z);
					
				\end{scope}

			\end{tikzpicture}
		\end{center}
		\caption{$F$-claws when $F$ hits all edges. The black vertices are in $F$.} \label{fig:stable}
	\end{figure}
	
	Suppose again that $F\subseteq V(H)$, but there is one special edge that may join two vertices in $V(H)\setminus F$ (we call it 
	the ``long edge''). If we enumerate the possible $F$-claws now, we still have the five cases of Figure \ref{fig:stable} (and any of their edges might be the long edge),
	and in addition there are three more shown in Figure \ref{fig:longedgeclaw}. (We omit the proof, which is easy case-checking.)
	\begin{figure}[H]
		\begin{center}
			\begin{tikzpicture}[scale = .8]
				
				\begin{scope}[shift ={(4,0)}]
					\tikzstyle{every node}=[inner sep=1.5pt, fill=black,circle,draw]
					
					\node (y1) at (-1,0) {};
					\node (y2) at (0,0) {};
					\node (y3) at (1,0) {};
					\node (z) at (0,-2) {};
					\tikzstyle{every node}=[inner sep=1.5pt, fill=white,circle,draw]
					\node (x1) at (1/3,-4/3) {};
					\node (x2) at (2/3,-2/3) {};

					\draw[] (y1) -- (z);
					\draw[] (y2) -- (z);
					\draw[] (y3) -- (x2);
					\draw[] (z) -- (x1);
					\draw[ultra thick] (x2) -- (x1);
				\end{scope}
				
				\begin{scope}[shift ={(8,0)}]
					\tikzstyle{every node}=[inner sep=1.5pt, fill=black,circle,draw]
					
					\node (y1) at (-1,0) {};
					\node (y2) at (0,0) {};
					\node (y3) at (1,0) {};
					\node (z) at (0,-2) {};
					\tikzstyle{every node}=[inner sep=1.5pt, fill=white,circle,draw]
					\node (w1) at (-1/2,-1) {};
					\node (x1) at (1/3,-4/3) {};
					\node (x2) at (2/3,-2/3) {};
					
					\draw (y1) -- (w1);
					\draw (w1) -- (z);
					\draw (y2) -- (z);
					\draw[] (y3) -- (x2);
					\draw[] (z) -- (x1);
					\draw[ultra thick] (x2) -- (x1);
					
				\end{scope}
				
				\begin{scope}[shift ={(12,0)}]
					\tikzstyle{every node}=[inner sep=1.5pt, fill=black,circle,draw]
					
					\node (y1) at (-1,0) {};
					\node (y2) at (0,0) {};
					\node (y3) at (1,0) {};
					\node (z) at (0,-2) {};
					\tikzstyle{every node}=[inner sep=1.5pt, fill=white,circle,draw]
					\node (w1) at (-1/2,-1) {};
					\node (w2) at (0,-1) {};
					\node (x1) at (1/3,-4/3) {};
					\node (x2) at (2/3,-2/3) {};
					
					\draw (y1) -- (w1);
					\draw (w1) -- (z);
					\draw (y2) -- (w2);
					\draw (w2) -- (z);
					\draw[] (y3) -- (x2);
					\draw[] (z) -- (x1);
					\draw[ultra thick] (x2) -- (x1);

				\end{scope}

			\end{tikzpicture}
		\end{center}
		\caption{$F$-claws containing the long edge (drawn thick).} \label{fig:longedgeclaw}
	\end{figure}
	
	Let $\mac T$ be a $k$th-century realm in $G$, with standard $\sigma$, in the usual notation. 
	A $\mac T$-community {\em adjoins} a fort $X$ if one of its
	members adjoins $X$.
	A {\em passage} is a path $P$ of $G$
	with the following properties:
	\begin{itemize}
		\item $P$ has length at most $3d_0+1$;
		\item there exist distinct $X_1,X_2\in \mac T$, each either a house or fort of $\mac T$, such that one end of $P$ is in $X_1$ and the other in $X_2$;
		\item no internal vertex of $P$ belongs to $X_1\cup X_2$; and
		\item $\dist_G(P,Y)>  \delta(k+1+\rk(Y))$ for each $Y\in \mac T$ with $Y\ne X_1,X_2$.
	\end{itemize}
	We say that $P$ {\em joins} $X_1,X_2$, and $P$ is {\em incident} with $X_1,X_2$. If $X_1$ is a house in a $\mac T$-community $C$, we also speak of 
	$P$ {\em joining} $C,X_2$. 
	\begin{thm}\label{adjoinpassage}
		If $X_1,X_2$ are houses or forts of $\mac T$, and $X_1$ adjoins $X_2$, there is a passage of length at most $2d_0+1$ joining $X_1,X_2$.
	\end{thm}
	\Proof
	Choose $v_i\in \voro_{\mac T}(X_i)$ for $i = 1,2$ such that $v_1,v_2$ are adjacent. Let $Q_i$ be the $\Lambda$-shortest path between 
	$v_i, X_i$; so $Q_i$
	is a path of $G[\voro_{\mac T}(X_i)]$ for $i = 1,2$. Let $P$ be the union of $Q_1,Q_2$ and the edge $v_1v_2$.
	Then $Q_1,Q_2$ both have length at most $d_0$ and so $P$ has length at most $2d_0+1$. Since $v_2\notin \voro_{\mac T}(X_1)$ and vice versa,
	it follows
	that the lengths of $Q_1,Q_2$ differ by at most one; and since the lengths sum to at least $\delta(\rk(X_1)+\rk(X_2))\ge 2c$, both lengths
	are at least $c$. 
	
	To show that $P$ is a passage, it remains to check that $\dist_G(P,Y)>  \delta(k+1+\rk(Y))$ for each $Y\in \mac T$ with $Y\ne X_1,X_2$. To see this, let $R$ be a shortest path from $Y$ to $V(P)$, with an end $p\in V(P)$. From the symmetry we may assume that $p\in V(Q_1)$.
	Let $P'$ be the subpath of $P$ between $X_1,p$. Since $p\in \voro_{\mac T}(X_1) \setminus \voro_{\mac T}(Y)$, it follows that $|E(R)|\ge |E(P')|$.
	But the sum of their lengths is more than $\delta(\rk(X_1)+\rk(Y))$, and so $R$ has length more than
	$\delta(\rk(X_1)+\rk(Y))/2\ge \delta(k+1+\rk(Y))$ since $\rk(X_1)\le k$. This proves \ref{adjoinpassage}.~\bbox

	We call passages as in \ref{adjoinpassage} {\em adjoinment passages}. (We will never need an adjoinment passage between two houses.)
	The eight graphs drawn in Figures \ref{fig:stable} and \ref{fig:longedgeclaw} all 
	represent configurations of forts, $\mac T$-communities and passages that we can prove do not appear if $\mac T$ is optimal; the 
	vertices in $F$ (drawn solid)
	represent forts, the other vertices (drawn hollow) represent pairwise disjoint $\mac T$-communities, and each edge represents a passage joining the forts or 
	$\mac T$-communities represented by its ends. All edges represent 
	adjoinment passages 
	except possibly one, the ``long edge'' when it appears. Thus, for example, the first graph in Figure \ref{fig:stable} shows a $\mac T$-community adjoining three forts. The first graph of Figure \ref{fig:longedgeclaw} shows a fort adjoining two forts and a $\mac T$-community $C$ say;
	$C$ is joined to another $\mac T$-community $D$ by a passage, and $D$ adjoins a fourth fort. 
	If $\sigma$ is a standard, we will often use the notation $\sigma=(\delta, \alpha,\beta)$ for the terms of $\sigma$
	without explicitly saying so.
	
	\begin{thm}\label{newaddclaw}
		Let $\sigma$ be a standard, let $G$ be a graph, and let 
		$\mathcal{T}$ be a $k$th-century realm in $G$, optimal for $\sigma$.
		There do not exist forts, pairwise disjoint $\mac T$-communities, and passages forming a configuration represented in any of the drawings in 
		Figures \ref{fig:stable} or \ref{fig:longedgeclaw}.
	\end{thm}
	\Proof
	Suppose that these things do exist; so, three or four forts $X_1\LL X_{m_1}$, at most four $\mac T$-communities $C_1\LL C_{m_2}$, 
	and at most seven passages $P_1\LL P_{m_3}$, making a configuration as in the drawing. 
	Choose the numbering such that $X_1,X_2,X_3$ are the three forts represented by the three leaves in the drawing, and $P_i$ is 
	the passage with an end in $X_i$ for $i = 1,2,3$. Again for $i = 1,2,3$, let $p_i$ be the end of $P_i$ in $X_i$, let $\eta_i$ exhibit $H_k$
	as a $c$-superfat minor of $G[X_i]$, let $J_i = \eta_i(H_k)$, and let $Q_i$ be a path of $G[X_i]$ between $p_i$ and $V(J_i)$. Choose a vertex $r_i$
	of the path $P_i\cup Q_i$ with $\dist_G(r_i, J_i)\le c+1$ such that the subpath of $P_i\cup Q_i$ between $r_i$ and $J_i$
	is maximal (thus, $r_i$ might lie in $V(P_i)$ or
	in $X_i$). Let $R_i$ be the $\Lambda$-shortest path in $G$ between $r_i$ and $J_i$. (Thus, $R_i$ is not necessarily a path of 
	$G[X_i]\cup P_i$.) Let $S_i$ be the subpath of $P_i\cup Q_i$ between $r_i$ and the end of $P_i\cup Q_i$ not in $X_i$. (See Figure \ref{fig:shortcut}.)
	From the choice of $r_i$, 
	$$\dist_G(S_i,J_i) = \dist_G(r_i, J_i) =|E(R_i)|=c+1.$$ 
	Let
	$$W=V(S_1\cup S_2\cup S_3)\cup \bigcup_{4\le i\le m_1} X_i \cup \bigcup_{1\le i\le m_2}\voro_{\mac T}(C_i) \cup \bigcup_{4\le i\le m_3} V(P_i).$$
	\begin{figure}[h!]
		\begin{center}
			\begin{tikzpicture}[]
				
				\tikzstyle{every node}=[inner sep=1.5pt, fill=black,circle,draw]
				
				
				\draw (0,0) ellipse (.7 and .5);
				\draw (0,0) ellipse (1.8 and 2.2);
				
				\node (p) at (-1.2,-1.5) {};
				\draw[very thick, dashed] (-3, 2) to [out= 270, in=180] (p);
				
				\draw[very thick, dashed] (p) to [out= 80, in=270] (.5,-.2);
				\node (r) at (-2.65,-.5) {};
				\draw[very thick, dashed] (r) to (-.6,-.1);
				
				\tikzstyle{every node}=[]
				
				\node at (0,0) {$J_i$};
				\node at (0,1.5) {$X_i$};
				\node[above right] at (r) {$r_i$};
				\node[below, right] at (p) {$p_i$};
				\node at (0,-1) {$Q_i$};
				\node at (-1.3,0) {$R_i$};
				\node at (-2.7,.6) {$S_i$};
				
				\node at (-4,-.5) {$P_i$};
				\draw[
				gray, ultra thin,decoration={markings,mark=at position 1 with {\arrow[black,scale=2]{>}}},
				postaction={decorate}
				]
				(-3.9,-.3) to[bend left=30] (-3,.4);
				\draw[
				gray, ultra thin,decoration={markings,mark=at position 1 with {\arrow[black,scale=2]{>}}},
				postaction={decorate}
				]
				
				(-3.9,-.7) to[bend right=30] (-2.4,-1);
			\end{tikzpicture}
		\end{center}
		\caption{For the proof of \ref{newaddclaw}. The vertex $r_i$ might belong to $Q_i$. } \label{fig:shortcut}
	\end{figure}
	\\
	\\
	(1) {\em $G[W]$ is connected, and $\dist_G(W, J_i) = c+1$ for $i = 1,2,3$.}
	\\
	\\
	$G[W]$ is clearly connected, since the corresponding configuration in Figures \ref{fig:stable} or \ref{fig:longedgeclaw} remains 
	connected when its leaves are deleted. Since $V(S_i)\subseteq W$ and $\dist_G(S_i, J_i)=c+1$, it suffices to show that $\dist_G(W, J_i) \ge c+1$ for $i = 1,2,3$. We claim:
	\begin{itemize}
		\item $\dist_G(S_i,J_i)\ge c+1$ from the choice of $r_i$;
		\item $\dist_G(P_j, J_i) \ge \dist_G(P_j, X_i)> \delta(k+1+\rk(X_i))=\delta(2k+1)\ge c+1$
		for all $j\in \{1\LL m_3\}$ with $j\ne i$, from the definition of a passage since $P_j$ is not incident with $X_i$;
		\item $\dist_G(S_j, J_i)\ge c+1$ for $j=1,2,3$ with $j\ne i$, since every vertex of $S_j$ belongs either to $P_j$
		(and the claim follows from the previous bullet) or to $X_j$ (and $\dist_G(X_j,X_i)> \delta(2k)\ge c+1$);
		\item $\dist_G(X_j, J_i)\ge \dist_G(X_j,X_i)> \delta(2k)\ge c+1$ for all $j\in \{1\LL m_1\}$ with $j\ne i$;
		\item $\dist_G(\voro_{\mac T}(C_j), J_i)\ge \dist_G(\voro_{\mac T}(C_j), X_i)\ge \dist_G(C_j, X_i)/2\ge c+1$, since
		if $v\in \voro_{\mac T}(C_j)$, then $\dist_G(v,C_j)\le \dist_G(v, X_i)$ and 
		$$\dist_G(v,C_j) + \dist_G(v, X_i)\ge \dist_G(C_j, X_i)\ge 2c+2.$$
	\end{itemize}
	This proves (1).

	\bigskip
	
	Let $Z=W\cup X_1\cup X_2\cup X_3\cup V(R_1\cup R_2\cup R_3)$.
	\\
	\\
	(2) {\em $Z$ has quasi-bound at most $(9\alpha,\beta+2d_0)$.}
	\\
	\\
	Let 
	$$A=\bigcup_{1\le i\le m_1} X_i \cup \bigcup_{1\le i\le m_2}V(C_i).$$
	Since each of
	$X_1\LL X_{m_1}, V(C_1)\LL V(C_{m_2})$
	has quasi-line-width at most $(\alpha,\beta)$, so does $A$ (because there are no edges between these sets).
	Since the boundary of each $X_i$ and of each $V(C_i)$ has quasi-size at most $(\alpha, \beta)$, and $m_1+m_2\le 8$, it follows that $\bd(A)$ 
	has quasi-size at most 
	$(8\alpha, \beta)$. Each vertex in $Z\setminus A$ has distance at most $(3d_0+1)/2\le 2d_0$ from some vertex in $\bd(A)$, and hence by \ref{increaserad}, 
	$Z$ has quasi-line-width at most $(9\alpha, \beta+2d_0)$, and $\bd(Z)$ has quasi-size at most $(8\alpha, \beta+2d_0)$. This proves (2).
	\\
	\\
	(3) {\em For each $Y\in \mathcal{T}$, either $\dist_G(Z,Y)> \delta(k+1+\rk(Y))$, or $Y\subseteq Z$.}
	\\
	\\
	We assume that $Y\not\subseteq Z$, and so $Y\ne X_1\LL X_{m_1}$. Let $M$ be a shortest path from $Y$ to $Z$, and let $z$ be its end in $Z$.
	Suppose that $z\in V(P_i)$ for some $i\in \{1\LL m_3\}$. Since $P_i$ is a passage joining two houses or forts included in $Z$, and 
	$Y\not\subseteq Z$, it follows that $P_i$ is not incident with $Y$, and so $\dist_G(P_i,Y)>  \delta(k+1+\rk(Y))$ from the definition of a passage, and the claim holds. 
	So we assume that either $z$ belongs to one of the forts
	$X_i$, or to one of the paths $R_i$, or to $\voro_{\mac T}(X)$ for some house $X$ in one of the $\mac T$-communities $C_i$. 
	In either case, $z\in \voro_{\mac T}(X)$ for some 
	house or fort $X$ included in $Z$ (because for $i = 1,2,3$, every vertex of $R_i$ has distance at most $c+1$ from $X_i$, and so
	$V(R_i)\subseteq  \voro_{\mac T}(X_i)$).
	Consequently 
	$$|E(M)|\ge \dist_G(X,Y)/2>\delta(\rk(X)+\rk(Y))/2\ge \delta(k+1+\rk(Y))$$
	as required. This proves (3).
	
	\bigskip
	
	Let
	$\mathcal{T}'$ consist of the set of members of $\mac T$ that are not included in $Z$, together with $Z$, where $Z$ is designated as a castle of $\mac T'$, and each other member
	of $\mac T'$ has the same class in $\mac T'$ that it has in $\mac T$.
	We claim that
	$\mac T'$ is a $k$th-century realm with standard $\sigma$.
	
	To show this, we must check that
	the members of $\mathcal{T}'$ are pairwise vertex-disjoint buildings (which is true, since any member of $\mac T'$ that intersects $Z$ is included in $Z$), and:
	\begin{itemize}
		\item For all $v\in V(G)$, there exists $X\in \mathcal{T}'$ such that $\dist_G(v,X)\le d_0$;
		\item If $X,Y\in \mac T'$ are distinct, then $\dist_G(X,Y)>\delta(\rk(X)+\rk(Y))$;
		\item Every fort or castle $X$ of $\mac T'$ includes a $c$-superfat $H_{\rk(X)}$-minor of $G$; and
		\item Every castle of $\mac T'$ has quasi-bound at most $(9\alpha, \beta+2d_0)$.
	\end{itemize}
	(The remaining conditions follow directly from the fact that $\mac T$ is a $k$th-century realm with standard~$\sigma$.)
	The first is clear. For the second, we may assume that $X=Z$, so $\rk(X)=k+1$. Since $Y\in \mac T'$ and therefore 
	$Y\not\subseteq Z$, the second bullet follows from (3).
	
	For the third, we only need check this when $X=Z$, and in that case the claim follows from \ref{claw} applied to $W, J_1,J_2,J_3$ and the geodesics $R_1,R_2,R_3$.
	The fourth bullet holds by (2). 
	This proves that $\mac T'$ is a $k$th-century realm with standard $\sigma$,
	contrary to the optimality of $\mac T$. 
	So there are no such
	$X_1\LL X_{m_1}, C_1\LL C_{m_2},P_1\LL P_{m_3}$. This proves \ref{newaddclaw}.~\bbox

	Let $C\subseteq \mac T$. We need to talk about $\mac T$-communities that are subsets of $C$ and maximal with this property.
	Let us call such sets {\em $C$-villages}. 
	If $X_1,X_2\in C$ are both forts, we say they {\em $C$-semiadjoin} if either $X_1$ adjoins $X_2$, or there is a $C$-village  
	that adjoins them both. If $X\in C$ is a fort, we say $X$ is {\em $C$-peripheral} if 
	$X$ $C$-semiadjoins at most one other fort in $C$. If $D$ is a $C$-village, we say $D$ is {\em $C$-peripheral} if $D$ adjoins at most one fort in $C$, and any such fort is $C$-peripheral.
	(See Figure \ref{fig:peripheral}.)
	If $A\subseteq V(G)$, $\mac T[A]$ denotes the set of members of $\mac T$ included in $A$. 
	\begin{figure}[H]
		\begin{center}
			\begin{tikzpicture}[]
				\node at (0,1) {\text{peripheral}};
				\tikzstyle{every node}=[inner sep=1.5pt, fill=black,circle,draw]
				\node (y1) at (1,0) {};
				\node (y2) at (2,0) {};
				\node (y3) at (3,0) {};
				\node (y4) at (4,0) {};
				\node (y5) at (5,0) {};
				\node (y6) at (6,0) {};
				\node (y7) at (7,0) {};
				\node (y8) at (8,0) {};
				\node (y9) at (9,0) {};
				\node (y10) at (10,0) {};
				
				\tikzstyle{every node}=[inner sep=1.5pt, fill=white,circle,draw]
				\node (m1) at (.4,.5) {};
				\node (m3) at (1.5,.5) {};
				\node (m4) at (2.5,0) {};
				\node (m5) at (2.5, .5) {};
				\node (m6) at (2.5, 1) {};
				\node (m7) at (4.5,.5) {};
				\node (m8) at (6.5,0) {};
				\node (m9) at (6.5,.5) {};
				\node (m10) at (6.5,1) {};
				\node (m11) at (7.5,.5) {};
				\node (m12) at (9.5,.5) {};
				
				\node (a1) at (.75,.7) {};
				\node (a2) at (1.1,.7) {};
				\node (a3) at (2,.7) {};
				\node (a4) at (4,.7) {};
				\node (a5) at (4.7,.7) {};
				\node (a6) at (5.3,.7) {};
				\node (a7) at (7,.7) {};
				\node (a8) at (8.7,.7) {};
				\node (a9) at (9,.7) {};
				\node (a10) at (9.3,.7) {};
				
				\foreach \from/\to in {y1/y2,y3/y4,y4/y5,y5/y6,y7/y8,y8/y9,y9/y10,
					m1/y1,m3/y1,m3/y2,m4/y2,m4/y3,m5/y2,m5/y3,m6/y2,m6/y3,m7/y4,m7/y5, m8/y6,m8/y7, m9/y6,m9/y7, m10/y6,m10/y7,m11/y7,m11/y8,m12/y9,m12/y10,
					a1/y1,a2/y1,a3/y2,a4/y4, a5/y5, a6/y5, a7/y7, a8/y9, a9/y9,a10/y9}
				\draw [-] (\from) -- (\to);

				\draw[dashed, thick] (y10) -- (10.7,0);
				\draw[dotted, thick] (-.1,.5) ellipse (1.4 and 1.1);
				
			\end{tikzpicture}
		\end{center}
		\caption{The black vertices represent forts in $C$, the others are $C$-villages, and edges are adjoinments.} \label{fig:peripheral}
	\end{figure}

	An {\em integer interval} is a set $I$ of integers such that if $a,b\in I$ then $I$ contains all integers between $a,b$. We say $(a,b)$
	is an {\em end-set} of an integer interval $I$ if $a,b\in I$, and $a-1,b+1\notin I$ (and hence $I$ is finite).
	If $\mac T$ is a realm, we denote the set of members of $\mac T$ that are houses or forts by $\mac T^0$.
	\begin{thm}\label{contactgraph}
		Let $\sigma$ be a standard, let
		$G$ be a graph, and
		let $\mathcal{T}$ be a $k$th-century realm in $G$, optimal for $\sigma$. 
		Each $\mac T$-community adjoins at most two forts, and 
		each fort $\mac T$-semiadjoins at most two other forts. Consequently, if 
		$C$ is an adjoin-connected subset of $\mac T^0$ that contains a fort, then  
		the forts in $C$ can be numbered as $X_i\;(i\in I)$, where $I$ is an integer interval, with the following properties:
		\begin{itemize}
			\item for each $i\in I$, if $i+1\in I$ then $X_i$ $C$-semiadjoins $X_{i+1}$, and no other pair of forts in $C$  $\mac T$-semiadjoin  each other 
			except possibly $X_a,X_b$, where $(a,b)$ is an end-set of $I$ and $b\ge a+2$;
			\item each $\mac T$-community  $S$  with $S\subseteq C$ adjoins at most two forts in $C$ (which therefore $C$-semiadjoin if there are two).
		\end{itemize}
		Moreover, if $C,D$ are disjoint adjoin-connected subsets of $\mac T^0$, and $P$ is a passage joining some $X\in C$ and some  $Y\in D$,
		then:
		\begin{itemize}
			\item  either $X$ is a $C$-peripheral fort, or $X$ belongs to a $C$-peripheral $C$-village, or $D$ contains no forts; and 
			\item either $Y$ 
			is a $D$-peripheral fort, or $Y$ belongs to a $D$-peripheral $D$-village,
			or $C$ contains no forts.
		\end{itemize}
	\end{thm}
	\Proof
	Each $\mac T$-community adjoins at most two forts, because otherwise we have the configuration represented by the first drawing in Figure \ref{fig:stable}. 
	Similarly, each fort $\mac T$-semiadjoins at most two forts, because otherwise we have a configuration represented by one of the other
	drawings in the same figure.
	
	Let $H$ be the graph with vertex set the set of forts in $C$, together with the set of all $C$-villages;
	two forts, or a fort and a $C$-village, are adjacent in $H$ if they adjoin. (No two $C$-villages adjoin.) 
	By \ref{noFclaw}, 
	the forts in $C$ can be numbered as in the theorem.
	
	Finally, suppose that $C,D$ are disjoint adjoin-connected subsets of $\mac T^0$, and $P$ is a passage joining some $X\in C$ and some  $Y\in D$.
	Suppose that $X$ is not a $C$-peripheral fort and does not belong to a $C$-peripheral $C$-village, and there is a fort in $D$. Thus either 
	\begin{itemize}
		\item $X$ is a fort, and $X$ $C$-semiadjoins two other forts in $C$; or
		\item $X$ is a house, and the $C$-village containing $X$ adjoins two forts in  $C$; or
		\item $X$ is a house, and the $C$-village containing $X$ adjoins a fort  in $C$ that is not $C$-peripheral.
	\end{itemize}
	Moreover, the end of $P$ in $V(D)$ either belongs to a fort of $D$, or it belongs to a $D$-village that adjoins 
	a fort of $D$. This again gives one of the configurations represented in Figures \ref{fig:stable} and \ref{fig:longedgeclaw}
	(the latter if $P$ joins a house with a house), in contradiction to \ref{newaddclaw}.
	This proves 
	\ref{contactgraph}.~\bbox
	
	We will need the following lemma about composing line-decompositions.
	\begin{thm}\label{composeline}
		Let $\mathcal{A}$ be a set of disjoint nonempty subsets of $V(G)$, with union $W$ say, such that each $A\in \mathcal{A}$
		has quasi-bound at most
		$(a,b)$. Suppose also that there is a line-decomposition of
		$G[W]$ in which each bag is the union of at most $k$ members of $\mathcal{A}$. Then $G[W]$ has quasi-line-width at most $((k+1)a,b)$.
	\end{thm}
	\Proof
	By hypothesis, there is a line-decomposition $(B_t':t\in T)$
	of $G[W]$
	such that each bag is the union of at most $k$ members of $\mac A$. 
	Let $Z$ be the union of the sets $\bd(A)\;(A\in \mathcal{A})$, and define $B_t=B'_t\cap Z$ for each $t\in T$. Then $(B_t:t\in T)$ 
	is a line-decomposition 
	of $G[Z]$
	such that each bag is the union of at most $k$ sets of the form $\bd(A)$, and so has quasi-size at most $(ka,b)$. For each $A\in \mathcal{A}$,
	since $A$ is nonempty, there exists $r(A)\in T$ such that $B'_{r(A)}\cap A\ne \emptyset$, and hence $\bd(A)\subseteq B_{r(A)}$ (since  $B_{r(A)}$ is a union of
	boundaries of members of $\mac A$ and the members of $\mac A$ are pairwise disjoint).
	By duplicating points of $T$, we may assume that the elements
	$r(A)\;(A\in \mathcal{A})$ are all distinct; and also by duplicating all the points of $T$, we may assume that each $r(A)$ has a successor 
	$s(A)\in T$ (that is,
	an element $s(A)\in T$ different from $r(A)$, such that $r(A)<s(A)$, and there is no $t\in T$ with $r(A)<t<s(A)$); and moreover, $B_{s(A)}=B_{r(A)}$. For each $A\in \mathcal{A}$,
	let $(C_t^A:t\in T^A)$ be a line-decomposition of $G[A]$ with quasi-width at most $(a,b)$. By adding two new elements to $T^A$, we may assume that
	$T^A$ has a maximum and minimum element; and so, by inserting $T^A$ into $T$, we may assume that $T^A$ equals
	$\{t\in T:r(A)\le t\le s(A)\}$, where $B_t = B_{r(A)}$ for each $t\in T^A$. For each $t\in T$, if $t\in T^A$ for some (necessarily unique) $A\in \mathcal{A}$, let
	$D_t=B_t\cup C^A_t$. If there is no such $A$, let $D_t=B_t$.
	
	We claim that $(D_t:t\in T)$ is a line-decomposition of $G[W]$. If $v\in W$, choose $A\in \mathcal{A}$ with $v\in A$; then $v\in C^A_t$
	for some $t\in T^A$, and hence $v\in D_t$. Next, suppose that $uv\in E(G[W])$.  If there exists $A\in \mathcal{A}$ with $u,v\in A$,
	then $u,v\in C^A_t$
	for some $t\in T^A$, and hence $u,v\in D_t$. If there is no such $A$, choose $A,A'\in \mathcal{A}$ with $u\in A$ and $v\in A'$.
	Then $u\in \bd(A)$ and $v\in \bd(A')$, and so $uv$ is an edge of $G[Z]$. Choose $t\in T$ such that $u,v\in B_t$; then $u,v\in D_t$.
	
	Finally, suppose that $r<s<t$ are elements of $T$, and $v\in D_r\cap D_t$. We need to show that $v\in D_s$. Choose $A\in \mathcal{A}$
	with $v\in A$. If $v\notin Z$, then each of $r,s,t$ belong to $T^A$, since no other bags contain $v$; and then, since
	$v\in C^A_r\cap C^A_t$, it follows that $v\in C^A_s\subseteq D_s$ as required. Now suppose that $v\in Z$. Then, for $p\in T$,
	$v$ belongs to $D_p$ if and only if $v\in B_p$; and so $v\in B_r\cap B_t\subseteq B_s\subset D_s$. This proves that
	$(D_t:t\in T)$ is a line-decomposition of $G[W]$. It is easy to check that each of its bags has 
	quasi-size at most $((k+1)a,b)$. This proves \ref{composeline}.~\bbox

	\begin{thm}\label{smallcoarsepw}
		Let $\sigma$ be a standard, let
		$G$ be a graph,
		let $\mathcal{T}$ be a $k$th-century realm in $G$, optimal for $\sigma$.
		Then $\voro_{\mac T} (\mac T^0)$
		has quasi-line-width at most
		$(10\alpha,\beta+d_0)$.
	\end{thm}
	\Proof
	Let $S=\voro_{\mac T} (\mac T^0)$.
	Let $J$ be the graph with vertex set the set of forts of $\mac T$  together with the set of all $\mac T$-villages, where we say $X,Y\in V(J)$
	are adjacent if $X$ adjoins $Y$ (and consequently at least one of $X,Y$ is a fort).
	From \ref{contactgraph}, each $\mac T$-village has degree at most two in $J$, and its neighbours in $J$ are forts. By \ref{contactgraph}, 
	$J$ has line-width at most three. For each $X\in V(J)$, let $A(X)=\voro_{\mac T}(X)$.
	Thus the sets $A(X)\;(X\in V(J))$ are pairwise disjoint and have union $S$. 
	Each $X\in V(J)$ has quasi-bound at most $(\alpha,\beta)$. Every vertex in $\voro_{\mac T}(X)$ has distance 
	at most $d_0$ from $X$, so by \ref{increaserad}, 
	it follows that $\bd(\voro_{\mac T}(X))$ has quasi-size at most $(\alpha,\beta+d_0)$, and
	$\voro_{\mac T}(X)$ has quasi-line-width at most $(2\alpha,\beta+d_0)$. 
	By \ref{composeline},
	$S$ has quasi-line-width at most $(10\alpha,\beta+d_0)$. This proves \ref{smallcoarsepw}.~\bbox

	\section{Governments}
	
	Let $0\le k\le \ell-1$, let $\sigma$ be a standard, let $G$ be a graph, and let 
	$\mathcal{T}$ be a $k$th-century realm in $G$, optimal for $\sigma$. 
	In the same way that we built a castle from three forts in \ref{newaddclaw}, now we want to build higher-level objects called ``palaces'', but the construction
	regulations are more complicated. 
	For $k+1\le t\le \ell$, a {\em palace (over $\mac T$)} of {\em rank} $t$ is a building $A$ in $G$
	with the following properties. If $t=k+1$ then $A$ is  a castle of $\mac T$, and $\mac P_A=\emptyset$. If $t\ge k+2$, then
	there is a set $\mac P_A$ of paths of $G[A]$, each of length at most $d_0+1$, and with $|\mac P_A|\le 3^{t-k}-3$, such that:
	\begin{itemize}
		\item every member of $\mac T$ is either a subset of $A$ or disjoint from $A$;
		\item $A\subseteq \left(\bigcup_{X\in \mac T[A]} \voro_{\mac T}(X)\right) \cup \left(\bigcup_{P\in \mac P_A}V(P)\right)$;
		\item $A$ includes a $c$-superfat $H_t$-minor of $G$; 
		\item there are at most $3^{t-k-1}-1$ forts of $\mac T$ that are included in $A$ and are $\mac T[A]$-peripheral;
		\item there are at most $(3^{t-k}-3)/2$ $\mac T[A]$-peripheral $\mac T[A]$-villages; and
		\item there are exactly $3^{t-k-1}$ castles of $\mac T$ included in $A$.
	\end{itemize}
	
	
	Define $\gamma = 8(\ell-k-1)3^{\ell-k-1}+1$.
	For $k+1\le t\le \ell$, define $\pi_t=3^{t-k-1}(10+\gamma) -\gamma$. 
	So $\pi_{k+1} = 10$, and $\pi_t=3\pi_{t-1}+2\gamma$ for $t\ge k+2$.
	We say a {\em government} for $\mac T$ is a set $\mac G$ of 
	pairwise disjoint buildings, each either a palace over $\mac T$, or a house or fort of $\mac T$, 
	satisfying:
	\begin{itemize}
		\item If $A_1,A_2\in \mac G$
		are distinct, of ranks $i,j$ respectively, then $\dist_G(A_1, A_2)> \delta(i+j)$.
		\item Each castle of $\mac T$ is a subset of some palace in $\mac G$.
		\item Every house or fort of $\mac T$ either belongs to $\mac G$ or is a subset of a palace in $\mac G$; we call the houses and forts that are not subsets of
		palaces {\em rebels}. (Thus, no castles are rebels.)
		\item Each palace $A\in \mac G$ has quasi-bound at most $(\pi_t\alpha,\beta+4d_0+1)$, 
		where $t$ is the rank of $A$.
	\end{itemize}
	If $A\in \mac G$, we denote its rank by $\gtype(A)$. (Rebel houses and forts have the same rank in $\mac G$ as they do in $\mac T$.)
	
	Thus, adding palaces is much like adding castles, but there is one important difference: the realm does not change. 
	When a house or fort is used to make 
	a castle, it is never going to be seen again, and we can forget it. But when a house, fort or castle is used to make a palace,
	it is not lost: it remains in the realm, and we may need it in the future, to help analyze properties of the palace.
	
	Keeping $k, \sigma, G, \mac T$ fixed, we will consider different governments for $\mac T$. There is at least one, because 
	$\mac T$ itself is a government for $\mac T$,
	since all castles have quasi-bound at most 
	$$(9\alpha,\beta+2d_0)\le (\pi_{k+1}\alpha,\beta+4d_0+1).$$
	We call this
	the {\em self-government} of $\mac T$.
	
	Let $\mac G$ be a government for $\mac T$.
	Thus, $\mac G$ is a sort of generalization of a realm, and once again we can use it to define Voronoi cells, using the same tie-breaker $\Lambda$ as before. 
	If $X\in \mac G$, let $\voro_{\mac G} (X)$ be the 
	corresponding Voronoi cell. 
	Now we have two different Voronoi partitions, one defined by $\mac T$ and one defined by $\mac G$. We need to work with them both, so from now on, 
	we write ``$\mac T$-adjoin'' or ``$\mac G$-adjoin'' to indicate which partition we are using.
	Since every member of $\mac T$ is a subset of a member of $\mac G$, it follows that if $X$ is a rebel, then
	$\voro_{\mac G}(X)\subseteq \voro_{\mac T}(X)$; and consequently every $\mac G$-community is a $\mac T$-community. The converse is false: sets of rebels that are $\mac T$-communities need not be $\mac G$-communities.
	
	\begin{thm}\label{palacebound}
		In the same notation, suppose that $G$ does not contain $H_\ell$ as a $c$-fat minor. Let $C$ be a $\mac T$-adjoin-connected set of rebels that contains at least
		one fort, and let $A\in \mac G$ be a palace.
		Let $U$ be the set of vertices in $\voro_{\mac G}(C)$ with a neighbour in  $\voro_{\mac G}(A)$.
		Then $U$ has quasi-size at most
		$$\left(4\cdot 3^{\ell-k-1}\alpha-2\alpha, \beta+4d_0+1\right).$$
	\end{thm}
	\Proof Let $A$ have rank $t$; so $k+1\le t\le \ell-1$ (since $G$ does not contain $H_\ell$ as a $c$-fat minor). 
	We say a subset of $A$ is {\em special} if either it is a castle, or it is a $\mac T[A]$-peripheral fort, or it is a house in 
	a $\mac T[A]$-peripheral $\mac T[A]$-village, or it is the vertex set of a member of $\mac P_A$.
	\\
	\\
	(1) {\em For each $u\in U$, 
		there is a path of length at most $2d_0+1$ between $u$ and the boundary of some special subset of $A$.}
	\\
	\\
	Since $u\notin A$, it suffices to show that there is a path of length at most $2d_0+1$ between $u$ and some special subset of $A$.
	Choose $X\in C$ with  $u\in \voro_{\mac G}(X)$, and $v \in \voro_{\mac G}(A)$ adjacent to $u$.
	Since $u\in \voro_{\mac G}(X)$, the $\Lambda$-shortest path $P$ between $u,X$ is in $G[\voro_{\mac G}(X)]$. 
	Let $Q$ be the $\Lambda$-shortest path between $v,A$ (which is therefore a path of $G[\voro_{\mac G}(A)]$). 
	Let the ends of $Q$ be $v,q$ where $q\in A$. If $q$ belongs to the vertex set of a member of $\mac P_A$, then
	the claim is true, so we may assume that there exists 
	$X'\in \mac T[A]$ with $q\in \voro_{\mac T}(X')$.
	If $X'$ is a castle of $\mac T$, then $X'$ is special and again the claim is true, so we assume that $\rk(X')\le k$.
	Let $R$ be the $\Lambda$-shortest path from $q$ to $X'$. Thus,
	$R$ has length 
	at most $d_0$. (But $R$ might not be contained in $G[A]$, and might contain vertices of $Q$ different from $q$.)
	Let $P'$ be a path between $X,X'$ in the union of $P$, the edge 
	$uv$, $Q$ and $R$. Thus, $P'$ has length at most $3d_0+1$.
	(See Figure \ref{fig:passage}.)
	\begin{figure}[H]
		\begin{center}
			\begin{tikzpicture}[]
				
				\tikzstyle{every node}=[inner sep=1.5pt, fill=black,circle,draw]
				
				\node (u) at (3,0) {};
				\node (v) at (4,0) {};
				\node (s) at (4.6,0) {};
				\node (q) at (6,0) {};
				\draw[dashed, thick]  (1,0) -- (u);
				\draw[dashed, thick] (v) -- (s) -- (q) -- (8,0);
				\draw[dashed, thick] (3.5,2) -- (s);
				\draw (u) -- (v);
				\draw[fill = white] (1,0) ellipse (1/2 and 1/2);
				\draw[fill = white] (8,0) ellipse (1/2 and 1/2);
				\draw[fill = white] (3.5,2) ellipse (1/2 and 1/2);
				\draw[thick, dotted] (3,0) arc [start angle=0, end angle=30, radius=3];
				\draw[thick, dotted] (3,0) arc [start angle=0, end angle=-30, radius=3];
				\draw[thick, dotted] (4,0) arc [start angle=180, end angle=150, radius=4];
				\draw[thick, dotted] (4,0) arc [start angle=180, end angle=200, radius=4];
				\draw[thick, dotted] (6,0) arc [start angle=180, end angle=135, radius=2];
				\draw[thick, dotted] (6,0) arc [start angle=180, end angle=215, radius=2];
				\tikzstyle{every node}=[]
				
				\node at (1,0) {$X$};
				\node at (8,0) {$X'$};
				\node at (3.5,2) {$Y$};
				\node at (2.2,-.3) {$P$};
				\node at (5,-.3) {$Q$};
				\node at (6.8,-.3) {$R$};
				\node at (2,1) {$\voro_{\mac G}(X)$};
				\node at (5.2,1) {$\voro_{\mac G}(A)$};
				\node at (7.2,1) {$\voro_{\mac T}(X')$};
				\draw (u) node[above right] {$u$};
				\draw (v) node[above left] {$v$};
				\draw (s) node[above right] {$s$};
				\draw (q) node[above left] {$q$};
				
			\end{tikzpicture}
			\caption{The passage. The picture is over-simplified, because $R$ might contain more vertices of $Q$.} \label{fig:passage}
		\end{center}
		
	\end{figure}

	We claim that $P'$ is a passage joining $X,X'$; and it 
	remains to check that 
	$$\dist_G(P',Y)>  \delta(k+1+\rk(Y))$$ 
	for each $Y\in \mac T$ with $Y\ne X,X'$. We will show the stronger result that $\dist_G(Y,s) >  \delta(k+1+\rk(Y))$
	for each $s\in V(P\cup Q\cup R)$ and each such $Y$. 
	
	Suppose first that $s\in V(P)$. 
	Since $s\in \voro_{\mac T}(X)$, $\dist_G(s,Y)\ge \dist_G(s,X)$ and the two sum to at least 
	$$\dist_G(X,Y)>\delta(\rk(X)+\rk(Y))\ge \delta(k+\rk(Y)).$$ 
	Hence 
	$$\dist_G(Y,s)> \delta(k+\rk(Y))/2\ge \delta(k+1+\rk(Y)),$$
	as required. 
	Similarly, if $s\in V(R)$, then 
	$$\dist_G(Y,s)> \delta(\rk(X')+\rk(Y))/2\ge \delta(k+1+\rk(Y))$$
	since $\rk(X')\le k$,
	so we assume that $s\in V(Q)$. 
	Write $\ell_1 = \dist_G(Y,s)$, and $\ell_2= \dist_G(s,q)$, and $\ell_3 = \dist_G(X', q)$ for brevity.
	Since $q\in  \voro_{\mac T}(X')$, it follows that 
	$$\dist_G(Y,s)+\dist_G(s,q)\ge \dist_G(Y,q)\ge \dist_G(X', q),$$
	and so $\ell_1+\ell_2\ge \ell_3$. 
	
	We recall that $Q$ is the $\Lambda$-shortest path between $v,A$, and $q$ is its end in $A$, and $s\in V(Q)$.
	Now there are three cases: $Y$ is a rebel; $Y$ is a subset of a palace in $\mac G$ different from $A$; and
	$Y\subseteq A$. But in each case, $\dist_G(Y,s)\ge \dist_G(s,q)$: if $Y$ is a rebel or a subset of a   
	different palace then because $s\in \voro_{\mac G}(A)$, and if $Y\subseteq A$ then trivially.
	Thus, $\ell_1\ge \ell_2$. 
	It follows that 
	$$\ell_1 + (2\ell_1)+ (\ell_1+\ell_2) \ge \ell_1+ 2\ell_2+ \ell_3,$$
	that is, $4\ell_1\ge \ell_1+\ell_2+\ell_3$. Consequently, 
	$$\dist_G(Y,s)\ge \dist_G(Y,X')/4> \delta(k+\rk(Y))/4\ge \delta(k+1+\rk(Y)),$$
	as required. (This is why we need the definition of ``spacing'' to be powers of 4, not of 2.)
	This proves that $P'$ is a passage joining $X,X'$.
	
	Now $X\in C$, and $X'$ is a house or fort in $\mac T[A]$. Let $D$ be the maximal subset of $\mac T[A]$ that contains $X'$ and is 
	$\mac T$-adjoin-connected. Since $C$ is $\mac T$-adjoin-connected, and contains a fort by hypothesis, 
	it follows from 
	\ref{contactgraph} that
	either $X'$
	is a $D$-peripheral fort, or $X'$ belongs to a $D$-peripheral $D$-village.
	But every $D$-peripheral fort is $\mac T[A]$-peripheral,
	and every $D$-peripheral $D$-village is a $\mac T[A]$-peripheral $\mac T[A]$-village, from the maximality of $D$. Consequently
	$X'$ is special.  
	This proves (1).

	\bigskip
	
	The special houses in $A$ can be partitioned into
	at most $(3^{t-k}-3)/2$ $\mac T$-villages, and there are at most $3^{t-k-1}-1$ special forts in $A$, and $3^{t-k-1}$ special castles,
	and $|\mac P_A|\le 3^{t-k}-3$.
	For each of these $\mac T$-villages,
	the union of the boundaries of its members has quasi-size at most $(\alpha,\beta)$; the boundary of each of the special forts has quasi-size at most 
	$(\alpha,\beta)$;
	the boundary of each castle has quasi-size at most $(9\alpha, \beta+2d_0)$, and the vertex set of each member of $\mac P_A$ 
	has quasi-size at most $(1,d_0)$. 
	Consequently, the union of the boundaries of the
	special subsets of $A$ has quasi-size at most
	\begin{align*}
		&\left(\left(3^{t-k}-3\right)(\alpha+4)/2+ \left(3^{t-k-1}-1\right)\alpha+ 3^{t-k-1}(9\alpha), \beta+2d_0\right)\\
		\le &
		\left(\left(4\cdot 3^{t-k}-2\right)\alpha,\beta+2d_0\right)
	\end{align*}
	(since $\alpha\ge 12$).
	Hence, by (1), we deduce that $U$ has quasi-size at most
	$$\left(\left(4\cdot 3^{t-k}-2\right)\alpha,\beta+4d_0+1\right).$$
	Since $t<\ell$, this proves \ref{palacebound}.~\bbox
	
	\section{Revolution}
	
	Again, let $\sigma$ be a standard, let $G$ be a graph that does not contain $H_{\ell}$ as a $c$-fat minor, let $\Lambda$ be a tie-breaker, and let
	$\mathcal{T}$ be a $k$th-century realm in $G$, optimal for $\sigma$.
	Let $\mac G$ be a government
	for $\mac T$. It follows that all members of $\mac G$ have rank $<\ell$. 
	For each rebel $X$, $\voro_{\mac G}(X)\subseteq \voro_{\mac T}(X)$,
	and so if two rebels $\mac G$-adjoin then they $\mac T$-adjoin each other; but whether two $\mac T$-adjoining rebels are $\mac G$-adjoining depends on the government.
	
	Suppose that there is a palace $A$ of $\mac T$ that does not belong to $\mac G$, that either includes or is disjoint from every member of $\mac G$. 
	Let $\mac G'$ be obtained from $\mac G$ by removing all members of $\mac G$ that are included in $A$, and adding $A$ itself. If $\mac G'$ is another government
	for $\mac T$, we say it is obtained from $\mac G$ by a {\em revolution}. Let us try to assemble the ingredients for a revolution.

	We say a set $C$ of rebels is {\em $\mac G$-village-closed} if the set of houses in $C$ is a union of $\mac G$-villages; in other words, if no rebel house outside of $C$ $\mac G$-adjoins a rebel house in $C$.
	A set $C$ of rebels is a {\em cabal} if:
	\begin{itemize}
		\item $C$ is $\mac G$-adjoin-connected and $\mac G$-village-closed;
		\item zero, one or two forts in $C$ are designated as ``leader forts''; and zero, one, two or three $\mac G$-villages included in $C$ are designated as
		``leader $\mac G$-villages'';
		\item every $C$-peripheral fort is a leader fort, and every $C$-peripheral $C$-village is a leader $\mac G$-village;
		\item if $X\in C$ $\mac G$-adjoins some rebel not in $C$, then either $X$ is a leader fort, or $X$ belongs to a leader $\mac G$-village;
		\item for some $j\in \{k+1\LL \ell-1\}$ there are three palaces $A_1,A_2,A_3\in \mac G$ of rank $j$, such that for $1\le i\le 3$,
		some member of $C$ $\mac G$-adjoins $A_i$; and
		\item either $|C|=1$, or $C$ is a $\mac G$-village, or for each $j'\in \{k+1\LL \ell-1\}$ there are at most four palaces $A\in \mac G$ of rank $j'$ such that
		some member of $C$ $\mac G$-adjoins $A$.
	\end{itemize}
	We will show that if we can find a cabal, we can use it to make a revolution by fusing together some three palaces of $\mac G$
	of the same rank to make a new palace of rank one bigger. The next result is used to control the quasi-bound of the new palace.
	
	\begin{thm}\label{cabalsaresmall}
		With notation as above, suppose that $C$ is a cabal.
		Then $\bd(\voro_{\mac G}(C))$ has quasi-size at most
		$$(16(\ell-k-1)3^{\ell-k-1}\alpha, \beta+4d_0+1).$$
	\end{thm}
	\Proof
	Since $C$ is a cabal, there exists $j\in \{k+1\LL \ell-1\}$ and three palaces $A_1,A_2,A_3\in \mac G$ of rank $j$, such that for $1\le i\le 3$,
	some member of $C$ $\mac G$-adjoins $A_i$. It follows that $k<j\le \ell-1$, and in particular $k\le \ell-2$.
	If $|C|=1$, say $C=\{X\}$, then $\bd(X)$ has quasi-size at most $(\alpha,\beta)$, and so $\bd(\voro_{\mac G}(C))$ has quasi-size at most $(\alpha,\beta+d_0)$, since every
	vertex in $\bd(\voro_{\mac G}(X))$ has distance at most $d_0$ from some vertex in $\bd(X)$. Similarly, if $C$ is a $\mac G$-village, and hence a $\mac T$-community, then $\bd(V(C))$
	has quasi-size at most
	$(\alpha,\beta)$, and hence again $\bd(\voro_{\mac G}(C))$ has quasi-size at most $(\alpha,\beta+d_0)$. Thus, we may assume that neither of these
	is true, and so for each $j\in \{k+1\LL \ell-1\}$ there are at most four palaces $A\in \mac G$ of rank $j$ such that some member of 
	$C$ $\mac G$-adjoins $A$.
	Let $\mac Q$ be the set of all such palaces in $\mac G$; so $|\mac Q|\le 4(\ell-k-1)$.  If some $X\in C$ $\mac G$-adjoins some rebel $Y$
	not in $C$, then either $X$ is a leader fort, or $X$ belongs to a leader $\mac G$-village; and
	the boundary of each leader fort has quasi-size at most $(\alpha, \beta)$, and so does the boundary of the union
	of the members of each leader $\mac G$-village. Consequently, $\bd(\voro_{\mac G}(C))$ is the union of at most five sets of quasi-size at most $(\alpha, \beta+d_0)$, and at most
	$4(\ell-k-1)$ further sets each with quasi-size at most
	$$\left(4\cdot 3^{\ell-k-1}\alpha-2\alpha, \beta+4d_0+1\right),$$
	by \ref{palacebound}.
	It follows that $\bd(\voro_{\mac G}(C))$ has quasi-size at most
	$$\left(5\alpha+4(\ell-k-1)\left(4\cdot 3^{\ell-k-1}\alpha-2\alpha\right),\beta+4d_0+1\right).$$
	This proves \ref{cabalsaresmall}.~\bbox

	\begin{figure}[h!]
		\begin{center}
			\begin{tikzpicture}[]
				
				\tikzstyle{every node}=[inner sep=1.5pt, fill=black,circle,draw]
				\def\q{1.4}
				\def\p{1.3}
				\def\m{3.5}
				\draw (0,0) ellipse ({\q} and {\q});
				\def\r{1.2}
				\node (y1) at (0,\r) {};
				\node (y2) at ({\r*cos(210)}, {\r*sin(210)}) {};
				\node (y3) at ({\r*cos(330)}, {\r*sin(330)}) {};
				\def\s{3}
				
				\draw[dashed, thick] (y1) -- (0,\s);
				\draw[thick, dashed] (y2) -- ({\s*cos(210)}, {\s*sin(210)});
				\draw[thick, dashed] (y3) -- ({\s*cos(330)}, {\s*sin(330)});
				
				\begin{scope}[shift ={(0,\m)}]
					\def\qq{.4}
					\def\pp{.3}
					\def\mm{1}
					\draw (0,0) ellipse ({\qq} and {\qq});
					\def\rr{.32}
					\def\ss{.85}
					\draw[dashed, thick] (0,\rr) -- (0,\ss);
					\draw[thick, dashed] ({\rr*cos(210)}, {\rr*sin(210)}) -- ({\ss*cos(210)}, {\ss*sin(210)});
					\draw[thick, dashed] ({\rr*cos(330)}, {\rr*sin(330)}) -- ({\ss*cos(330)}, {\ss*sin(330)});
					
					\begin{scope}[shift ={(0,1)}, scale= .23]
						\def\qq{.4}
						\def\pp{.3}
						\def\mm{1}
						\draw (0,0) ellipse ({\qq} and {\qq});
						\draw (0,\mm) ellipse ({\pp} and {\pp});
						\draw ({\mm*cos(210)}, {\mm*sin(210)}) ellipse ({\pp} and {\pp});
						\draw ({\mm*cos(330)}, {\mm*sin(330)}) ellipse ({\pp} and {\pp});
						\def\rr{.32}
						\def\ss{.85}
						\draw[dashed, thick] (0,\rr) -- (0,\ss);
						\draw[thick, dashed] ({\rr*cos(210)}, {\rr*sin(210)}) -- ({\ss*cos(210)}, {\ss*sin(210)});
						\draw[thick, dashed] ({\rr*cos(330)}, {\rr*sin(330)}) -- ({\ss*cos(330)}, {\ss*sin(330)});
					\end{scope}
					
					\begin{scope}[shift ={({cos(210)}, {sin(210)})}, scale= .23]
						\def\qq{.4}
						\def\pp{.3}
						\def\mm{1}
						\draw (0,0) ellipse ({\qq} and {\qq});
						\draw (0,\mm) ellipse ({\pp} and {\pp});
						\draw ({\mm*cos(210)}, {\mm*sin(210)}) ellipse ({\pp} and {\pp});
						\draw ({\mm*cos(330)}, {\mm*sin(330)}) ellipse ({\pp} and {\pp});
						\def\rr{.32}
						\def\ss{.85}
						\draw[dashed, thick] (0,\rr) -- (0,\ss);
						\draw[thick, dashed] ({\rr*cos(210)}, {\rr*sin(210)}) -- ({\ss*cos(210)}, {\ss*sin(210)});
						\draw[thick, dashed] ({\rr*cos(330)}, {\rr*sin(330)}) -- ({\ss*cos(330)}, {\ss*sin(330)});
					\end{scope}
					
					\begin{scope}[shift ={({cos(330)}, {sin(330)})}, scale= .23]
						\def\qq{.4}
						\def\pp{.3}
						\def\mm{1}
						\draw (0,0) ellipse ({\qq} and {\qq});
						\draw (0,\mm) ellipse ({\pp} and {\pp});
						\draw ({\mm*cos(210)}, {\mm*sin(210)}) ellipse ({\pp} and {\pp});
						\draw ({\mm*cos(330)}, {\mm*sin(330)}) ellipse ({\pp} and {\pp});
						\def\rr{.32}
						\def\ss{.85}
						\draw[dashed, thick] (0,\rr) -- (0,\ss);
						\draw[thick, dashed] ({\rr*cos(210)}, {\rr*sin(210)}) -- ({\ss*cos(210)}, {\ss*sin(210)});
						\draw[thick, dashed] ({\rr*cos(330)}, {\rr*sin(330)}) -- ({\ss*cos(330)}, {\ss*sin(330)});
					\end{scope}

				\end{scope}
				\begin{scope}[shift ={({\m*cos(210)}, {\m*sin(210)})}]
					\def\qq{.4}
					\def\pp{.3}
					\def\mm{1}
					\draw (0,0) ellipse ({\qq} and {\qq});
					\def\rr{.32}
					\def\ss{.85}
					\draw[dashed, thick] (0,\rr) -- (0,\ss);
					\draw[thick, dashed] ({\rr*cos(210)}, {\rr*sin(210)}) -- ({\ss*cos(210)}, {\ss*sin(210)});
					\draw[thick, dashed] ({\rr*cos(330)}, {\rr*sin(330)}) -- ({\ss*cos(330)}, {\ss*sin(330)});
					
					\begin{scope}[shift ={(0,1)}, scale= .23]
						\def\qq{.4}
						\def\pp{.3}
						\def\mm{1}
						\draw (0,0) ellipse ({\qq} and {\qq});
						\draw (0,\mm) ellipse ({\pp} and {\pp});
						\draw ({\mm*cos(210)}, {\mm*sin(210)}) ellipse ({\pp} and {\pp});
						\draw ({\mm*cos(330)}, {\mm*sin(330)}) ellipse ({\pp} and {\pp});
						\def\rr{.32}
						\def\ss{.85}
						\draw[dashed, thick] (0,\rr) -- (0,\ss);
						\draw[thick, dashed] ({\rr*cos(210)}, {\rr*sin(210)}) -- ({\ss*cos(210)}, {\ss*sin(210)});
						\draw[thick, dashed] ({\rr*cos(330)}, {\rr*sin(330)}) -- ({\ss*cos(330)}, {\ss*sin(330)});
					\end{scope}
					
					\begin{scope}[shift ={({cos(210)}, {sin(210)})}, scale= .23]
						\def\qq{.4}
						\def\pp{.3}
						\def\mm{1}
						\draw (0,0) ellipse ({\qq} and {\qq});
						\draw (0,\mm) ellipse ({\pp} and {\pp});
						\draw ({\mm*cos(210)}, {\mm*sin(210)}) ellipse ({\pp} and {\pp});
						\draw ({\mm*cos(330)}, {\mm*sin(330)}) ellipse ({\pp} and {\pp});
						\def\rr{.32}
						\def\ss{.85}
						\draw[dashed, thick] (0,\rr) -- (0,\ss);
						\draw[thick, dashed] ({\rr*cos(210)}, {\rr*sin(210)}) -- ({\ss*cos(210)}, {\ss*sin(210)});
						\draw[thick, dashed] ({\rr*cos(330)}, {\rr*sin(330)}) -- ({\ss*cos(330)}, {\ss*sin(330)});
					\end{scope}
					
					\begin{scope}[shift ={({cos(330)}, {sin(330)})}, scale= .23]
						\def\qq{.4}
						\def\pp{.3}
						\def\mm{1}
						\draw (0,0) ellipse ({\qq} and {\qq});
						\draw (0,\mm) ellipse ({\pp} and {\pp});
						\draw ({\mm*cos(210)}, {\mm*sin(210)}) ellipse ({\pp} and {\pp});
						\draw ({\mm*cos(330)}, {\mm*sin(330)}) ellipse ({\pp} and {\pp});
						\def\rr{.32}
						\def\ss{.85}
						\draw[dashed, thick] (0,\rr) -- (0,\ss);
						\draw[thick, dashed] ({\rr*cos(210)}, {\rr*sin(210)}) -- ({\ss*cos(210)}, {\ss*sin(210)});
						\draw[thick, dashed] ({\rr*cos(330)}, {\rr*sin(330)}) -- ({\ss*cos(330)}, {\ss*sin(330)});
					\end{scope}
					
				\end{scope}
				\begin{scope}[shift ={({\m*cos(330)}, {\m*sin(330)})}]
					\def\qq{.4}
					\def\pp{.3}
					\def\mm{1}
					\draw (0,0) ellipse ({\qq} and {\qq});
					\def\rr{.32}
					\def\ss{.85}
					\draw[dashed, thick] (0,\rr) -- (0,\ss);
					\draw[thick, dashed] ({\rr*cos(210)}, {\rr*sin(210)}) -- ({\ss*cos(210)}, {\ss*sin(210)});
					\draw[thick, dashed] ({\rr*cos(330)}, {\rr*sin(330)}) -- ({\ss*cos(330)}, {\ss*sin(330)});
					
					\begin{scope}[shift ={(0,1)}, scale= .23]
						\def\qq{.4}
						\def\pp{.3}
						\def\mm{1}
						\draw (0,0) ellipse ({\qq} and {\qq});
						\draw (0,\mm) ellipse ({\pp} and {\pp});
						\draw ({\mm*cos(210)}, {\mm*sin(210)}) ellipse ({\pp} and {\pp});
						\draw ({\mm*cos(330)}, {\mm*sin(330)}) ellipse ({\pp} and {\pp});
						\def\rr{.32}
						\def\ss{.85}
						\draw[dashed, thick] (0,\rr) -- (0,\ss);
						\draw[thick, dashed] ({\rr*cos(210)}, {\rr*sin(210)}) -- ({\ss*cos(210)}, {\ss*sin(210)});
						\draw[thick, dashed] ({\rr*cos(330)}, {\rr*sin(330)}) -- ({\ss*cos(330)}, {\ss*sin(330)});
					\end{scope}
					
					\begin{scope}[shift ={({cos(210)}, {sin(210)})}, scale= .23]
						\def\qq{.4}
						\def\pp{.3}
						\def\mm{1}
						\draw (0,0) ellipse ({\qq} and {\qq});
						\draw (0,\mm) ellipse ({\pp} and {\pp});
						\draw ({\mm*cos(210)}, {\mm*sin(210)}) ellipse ({\pp} and {\pp});
						\draw ({\mm*cos(330)}, {\mm*sin(330)}) ellipse ({\pp} and {\pp});
						\def\rr{.32}
						\def\ss{.85}
						\draw[dashed, thick] (0,\rr) -- (0,\ss);
						\draw[thick, dashed] ({\rr*cos(210)}, {\rr*sin(210)}) -- ({\ss*cos(210)}, {\ss*sin(210)});
						\draw[thick, dashed] ({\rr*cos(330)}, {\rr*sin(330)}) -- ({\ss*cos(330)}, {\ss*sin(330)});
					\end{scope}
					
					\begin{scope}[shift ={({cos(330)}, {sin(330)})}, scale= .23]
						\def\qq{.4}
						\def\pp{.3}
						\def\mm{1}
						\draw (0,0) ellipse ({\qq} and {\qq});
						\draw (0,\mm) ellipse ({\pp} and {\pp});
						\draw ({\mm*cos(210)}, {\mm*sin(210)}) ellipse ({\pp} and {\pp});
						\draw ({\mm*cos(330)}, {\mm*sin(330)}) ellipse ({\pp} and {\pp});
						\def\rr{.32}
						\def\ss{.85}
						\draw[dashed, thick] (0,\rr) -- (0,\ss);
						\draw[thick, dashed] ({\rr*cos(210)}, {\rr*sin(210)}) -- ({\ss*cos(210)}, {\ss*sin(210)});
						\draw[thick, dashed] ({\rr*cos(330)}, {\rr*sin(330)}) -- ({\ss*cos(330)}, {\ss*sin(330)});
					\end{scope}
					
				\end{scope}

				\tikzstyle{every node}=[]
				\node at (-.7,2.7+1) {$A_1$};
				\node at ({2.7*cos(210)-1.4}, {2.7*sin(210)}) {$A_2$};
				\node at ({2.7*cos(330)+1.3}, {2.7*sin(330)}) {$A_3$};
				\node at ({2.2*cos(83)+.1}, {2.2*sin(83)}) {$S_1'$};
				\node at ({2.2*cos(203)+.4}, {2.2*sin(203)-.5}) {$S_2'$};
				\node at ({2.2*cos(337)-.4}, {2.2*sin(337)-.6}) {$S_3'$};
				
				\node at (0,0) {$\voro_{\mac G}(C)$};
				
			\end{tikzpicture}
		\end{center}
		
		\caption{The headquarters of the cabal $C$. 
			The three objects labelled $A_i$ are palaces of rank $j$.} \label{fig:palace}
	\end{figure}
	
	With notation as before, suppose that $C$ is a cabal, and let
	$j, A_1,A_2,A_3$ be as in the fifth bullet in the definition of a cabal.
	For $i = 1,2,3$, choose $u_i\in \voro_{\mac G}(C)$ with a neighbour $v_i\in \voro_{\mac G}(A_i)$, and 
	define $P_i$ to be the path formed by the union of the edge $u_iv_i$ and 
	the 
	$\Lambda$-shortest path in $G$ between $v_i, A_i$. Hence, $P_i$ has length at most $d_0+1$.
	Again for $i = 1,2,3$, let $p_i$ be the end of $P_i$ in $A_i$, let $\eta_i$ exhibit $H_j$
	as a $c$-superfat minor of $G[A_i]$, let $J_i = \eta_i(H_j)$, and let $Q_i$ be a path of $G[A_i]$ between $p_i$ and $V(J_i)$. Choose a vertex $r_i$
	of the path $P_i\cup Q_i$ with $\dist_G(r_i, J_i)\le c+1$ such that the subpath of $P_i\cup Q_i$ between $r_i$ and $J_i$
	is maximal (thus, $r_i$ might lie in $V(P_i)$ or
	in $A_i$). Let $R_i$ be the $\Lambda$-shortest path in $G$ between $r_i$ and $J_i$. (Thus, $R_i$ is not necessarily a path of
	$G[A_i]\cup P_i$; but, since it has length $c+1$ and has an end in $J_i$, all its vertices belong to $\voro_{\mac G}(A_i)$.) 
	Let $S_i$ be the subpath of $P_i\cup Q_i$ between $u_i,r_i$. 
	From the choice of $r_i$,
	$$\dist_G(S_i,J_i) = \dist_G(r_i, J_i) =|E(R_i)|=c+1.$$
	Since $V(P_i)\subseteq \voro_{\mac G}(A_i)\cup \{u_i\}$ and $V(Q_i)\subseteq A_i$, and $V(R_i)\subseteq \voro_{\mac G}(A_i)$,
	and $S_i$ is a subpath of $P_i\cup Q_i$, it follows that $V(P_i\cup Q_i\cup R_i\cup S_i)\subseteq \voro_{\mac G}(A_i)\cup \{u_i\}$ for $i = 1,2,3$.
	Let 
	$$A=\voro_{\mac G}(C)\cup A_1\cup A_2\cup A_3 \cup V(P_1\cup R_1)\cup V(P_2\cup R_2)\cup V(P_3\cup R_3),$$
	and define $\mac P_A$ to be the 
	union of $\{P_1,P_2,P_3,R_1,R_2,R_3\}$ and the three sets $\mac P_{A_i}$ for $i = 1,2,3$ (taking $\mac P_{A_i}=\emptyset$
	if $A_i$ is a castle). It follows that 
	$$A\subseteq \voro_{\mac G}(C)\cup \voro_{\mac G}(A_1)\cup \voro_{\mac G}(A_2)\cup \voro_{\mac G}(A_3).$$
	We call $A$ the {\em headquarters} of the cabal $C$. (See Figure \ref{fig:palace}.)
	We will show that $A$ is a new palace, and causes 
	a revolution. 
	
	\begin{thm}\label{usecabal1}
		With notation as above, let $A$ be the headquarters of the cabal $C$.
		Then $A$ is a palace over $\mac T$ of rank $j+1$.
	\end{thm}
	\Proof 
	Certainly $G[\voro_{\mac G}(C)]$ is connected, because 
	$C$ is $\mac G$-adjoin-connected, and $G[\voro_{\mac G}(X)]$ is connected for each $X\in C$. 
	Moreover, $G[A_i]$ is connected, and the path $P_i$ joins
	$G[\voro_{\mac G}(C)]$ and $G[A_i]$ for $i = 1,2,3$, and so $G[A]$ is connected.
	
	The members of $\mac P_A$ have length at most $d_0+1$, since this is true for the paths in $\mac P_{A_i}$ for $i = 1,2,3$, and true for 
	$P_i, R_i$ for $i = 1,2,3$. Moreover, 
	$$|\mac P_A|=6+\sum_{1\le i\le 3}|\mac P_{A_i}|\le 6+3(3^{j-k}-3)=3^{j+1-k}-3.$$
	We claim next that:
	\\
	\\
	(1) {\em $A\subseteq \left(\bigcup_{X\in \mac T[A]} \voro_{\mac T}(X)\right) \cup \left(\bigcup_{P\in \mac P_A}V(P)\right)$.}
	\\
	\\
	For $i =1,2,3$, since $A_i$ is a palace and $\mac P_{A_i}\subseteq \mac P_A$, it follows that 
	$$A_i\subseteq \bigcup_{X\in \mac T[A_i]} \voro_{\mac T}(X)\cup \bigcup_{P\in \mac P_{A_i}}V(P) \subseteq \bigcup_{X\in \mac T[A]} \voro_{\mac T}(X)\cup \bigcup_{P\in \mac P_A}V(P).$$
	Also, 
	$\voro_{\mac G}(C)\subseteq \bigcup_{X\in \mac T[A]} \voro_{\mac T}(X)$ since $C\subseteq \mac T[A]$. And trivially
	$V(P_i\cup R_i)\subseteq \bigcup_{P\in \mac P_A}V(P)$ since $P_i, R_i \in \mac P_A$. This proves (1).
	\\
	\\
	(2) {\em Every member of $\mac T$ is a subset of $A$ or disjoint from $A$.}
	\\
	\\
	Let $Y\in \mac T$ with $Y\cap A\ne \emptyset$, and choose $y\in Y\cap A$. Thus, by (1), either $y\in \voro_{\mac T}(X)$ for some 
	$X\in \mac T[A]$, or $y\in V(P)$ for some $P\in \mac P_A$. In the first case, $Y\cap \voro_{\mac T}(X)\ne \emptyset$, and so $Y=X$
	since $X,Y\in \mac T$, and hence $Y\subseteq A$. In the second case, either $y$ belongs to some member of
	$\mac P_{A_i}$ for some $i\in \{1,2,3\}$, or $y\in V(P_i\cup R_i)$ for some $i\in \{1,2,3\}$. If $y$ belongs to some member 
	of $\mac P_{A_i}$ for some $i\in \{1,2,3\}$, then $Y\cap A_i\ne \emptyset$, and so $Y\subseteq A_i$ (because $A_i$ is a palace)
	and so $Y\subseteq A$. Thus we may assume that 
	$$y\in V(P_i\cup R_i)\subseteq \voro_{\mac G}(A_i)\cup \{u_i\}.$$
	Hence either $Y\cap \voro_{\mac G}(A_i)\ne \emptyset$, or $Y\cap  \voro_{\mac G}(C)\ne \emptyset$, and in either case it follows that 
	$Y\subseteq A$. This proves (2).
	
	\bigskip
	
	To show that $A$ is a palace, it remains to check that:
	\begin{itemize}
		\item $A$ includes a $c$-superfat $H_{j+1}$-minor of $G$;
		\item there are at most $3^{j-k}-1$ forts of $\mac T$ that are included in $A$ and are $\mac T[A]$-peripheral;
		\item there are at most $(3^{j+1-k}-3)/2$ $\mac T[A]$-peripheral $\mac T[A]$-villages; and
		\item there are exactly $3^{j-k}$ castles of $\mac T$ included in $A$.
	\end{itemize}
	The first follows from \ref{claw}, applied with $W=\voro_{\mac G}(C)\cup V(S_1\cup S_2\cup S_3)$ and the geodesics $R_1,R_2,R_3$.
	For the second, every $\mac T[A]$-peripheral fort included in $A$ either belongs to $C$ or is a $\mac T[A_i]$-peripheral fort
	included in some $A_i$. Since there are at most two in $C$ (from the definition of a cabal) and at most $3^{j-1-k}-1$ in $A_i$
	for $i = 1,2,3$, this proves the second bullet. The third and fourth are proved similarly. This proves \ref{usecabal1}.~\bbox
	
	\begin{thm}\label{usecabal2}
		In the same notation, let $A$ be the headquarters of the cabal $C$, and $\mac G'$ be the union of $\{A\}$ and the set of members of $\mac G$
		that are disjoint from $A$. Then $\mac G'$ is a government for $\mac T$.
	\end{thm}
	\Proof
	We need to check that:
	\begin{itemize}
		\item If $A'\in \mac G'$ with $A'\ne A$
		then $\dist_G(A, A')> \delta(\rk(A')+j+1)$.
		\item Each castle of $\mac T$ is a subset of some palace in $\mac G'$.
		\item Every house or fort of $\mac T$ either belongs to $\mac G'$ or is a subset of a palace in $\mac G'$. 
		\item Each palace $A'\in \mac G'$ has quasi-bound at most
		$\left(\pi_t\alpha, \beta+4d_0+1 \right)$
		where $t$ is the rank of $A'$.
	\end{itemize}
	For the first bullet, let $A'\in \mac G'$ with $A'\ne A$. Let $v\in A$; we need to show that $\dist_G(v,A')>\delta(\gtype(A')+j+1)$.
	Since 
	$$v\in A\subseteq \voro_{\mac G}(C)\cup \voro_{\mac G}(A_1)\cup \voro_{\mac G}(A_2)\cup \voro_{\mac G}(A_3),$$
	there exists $X\in \mac G$ with rank at most $j$, with $X\subseteq A$, and with $v\in  \voro_{\mac G}(X)$. Since $A',X\in \mac G$ and 
	$v\in  \voro_{\mac G}(X)$, it follows that 
	$\dist_G(v,A')\ge \dist_G(v,X)$, but 
	$$\dist_G(v,A')+\dist_G(v,X)>\delta(\gtype(A')+\rk(X)),$$
	so
	$$\dist_G(v,A')> \delta(\gtype(A')+\rk(X))/2> \delta(\gtype(A')+j+1).$$
	This proves the first bullet.
	The second holds since it holds for $\mac G$, and similarly so does the third.
	Finally, for the fourth bullet, it is only necessary to check the claim when $A' = A$.
	We know:
	\begin{itemize}
		\item $A_1, A_2, A_3$ all have quasi-bound at most
		$(\pi_{j}\alpha, \beta+4d_0+1)$; 
		\item $V(P_i), V(R_i)$ have quasi-bound at most $(1,d_0)$ for $i = 1,2,3$;
		\item 
		$\voro_{\mac G}(C)$ has quasi-line-width at most $(10\alpha,\beta+d_0)$, by \ref{smallcoarsepw}; and 
		\item  $\bd(\voro_{\mac G}(C))$ has quasi-size at most
		$(16(\ell-k-1)3^{\ell-k-1}\alpha, \beta+4d_0+1)$ by \ref{cabalsaresmall}.
	\end{itemize}
	Since there are no edges between any of the sets $A_1, A_2, A_3, \voro_{\mac G}(C)\setminus \bd(\voro_{\mac G}(C))$, 
	and $\pi_{j}\ge 10$, 
	it follows that the union of these four sets also has quasi-line-width at most $(\pi_{j}\alpha, \beta+4d_0+1)$.
	By adding 
	$$\bd(\voro_{\mac G}(C))\cup V(P_1\cup P_2\cup P_3\cup R_1\cup R_2\cup R_3)$$
	to each of the bags of a corresponding line-decomposition, we deduce that $A$ has quasi-line-width at most 
	$$(\pi_{j}\alpha+ 16(\ell-k-1)3^{\ell-k-1}\alpha+6, \beta+4d_0+1).$$
	Its boundary $\bd(A)$ is a subset of the union of the boundaries of $A_1, A_2, A_3, \voro_{\mac G}(C)$ and the set
	$V(P_1\cup P_2\cup P_3\cup R_1\cup R_2\cup R_3)$, and so has quasi-size at most
	$$(3\pi_{j}\alpha+ 16(\ell-k-1)3^{\ell-k-1}\alpha+6, \beta+4d_0+1).$$
	Since 
	$$\pi_{j+1}= 3\pi_{j}+ 2\gamma \ge 3\pi_{j}+ 16(\ell-k-1)3^{\ell-k-1}+1$$ 
	and $\alpha\ge 6$, it follows that
	$A$ has quasi-bound at most $(\pi_{j+1}\alpha, \beta+4d_0+1)$. 
	This proves the fourth bullet, and consequently
	$\mac G'$ is a government. This proves \ref{usecabal2}.~\bbox
	\section{Stable government}
	
	With the century $k$ fixed, and given a $k$th-century realm $\mac T$ that is optimal for some standard $\sigma$,
	we want to modify the self-government to make the best possible government. To do so, we say a government $\mac G_2$
	{\em extends} a government $\mac G_1$ if:
	\begin{itemize}
		\item every rebel of $\mac G_1$ either is a rebel of $\mac G_2$, or is a subset of a palace of $\mac G_2$;
		\item every palace $X_1$ of $\mac G_1$ is a subset of a palace $X_2$ of $\mac G_2$, where $X_1\subseteq X_2$, and either
		the rank of $X_2$ in $\mac G_2$ is greater than the
		rank of $X_1$ in $\mac G_1$, or $X_1=X_2$ and the two ranks are equal.
	\end{itemize}
	Let us say a government $\mac G$ is {\em stable} if no other government  extends $\mac G$.  Since the  self-government 
	exists, it follows from Zorn's lemma that there is a stable government. If a government is stable, then
	by \ref{usecabal1} and \ref{usecabal2} there are no cabals.
	
	\begin{thm}\label{notmuchleft}
		With notation as before, 
		let $\mac G$ be a stable government for $\mac T$, and
		let $C$ be a maximal $\mac G$-adjoin-connected set of rebels. Then
		$\voro_{\mac G}(C)$ has quasi-bound at most 
		$$\left(\ell\,3^{\ell+1}\alpha, \beta+4d_0+1\right).$$
	\end{thm}
	\Proof
	By \ref{smallcoarsepw}, $\voro_{\mac G}(C)\subseteq \voro_{\mac T}(C)$ has quasi-line-width at most $(10\alpha,\beta+d_0)$, so it remains to bound the quasi-size of $\bd(\voro_{\mac G}(C))$. If $k=\ell-1$, then there are no castles in $\mac T$ and no palaces in $\mac G$ 
	(because $G$ does not contain $H_{\ell}$ as a $c$-fat minor), and so $\bd(\voro_{\mac G}(C))=\emptyset$
	and the claim is true. So we may assume that $k\le \ell-2$.
	
	If $C$ is a $\mac T$-community, then, since $\bd(V(C))$ has quasi-size at most $(\alpha,\beta)$ and every vertex in $\bd(\voro_{\mac G}(C))$ has distance at most
	$d_0$ from $\bd(V(C))$, it follows that $\bd(\voro_{\mac G}(C))$ has quasi-size at most $(\alpha,\beta+d_0)$.
	So we may assume that $C$ is not a $\mac T$-community.
	
	We say that $C'\subseteq C$ is {\em dangerous} if $C'$ is $\mac G$-adjoin-connected, and 
	there exist $j>k$ and three members $A_1,A_2,A_3\in \mac G$, all of rank $j$, such that for $1\le i\le 3$, some member of $C'$ $\mac G$-adjoins $A_i$.
	If $C$ is dangerous, we want to choose a minimal subset $D'$ of $C$ that is still dangerous, and $\mac G$-village-closed, and $\mac G$-adjoin-connected, 
	\\
	\\
	(1) {\em $C$ is not dangerous.}
	\\
	\\
	Suppose that $C$ is dangerous, but let us digress for a moment and consider what we need. What we would really like is to obtain 
	a ``good'' dangerous subset $C'$, where ``good'' means 
	\begin{itemize}
		\item $C'$ is $\mac G$-village-closed and  $\mac G$-adjoin-connected, 
		\item only a bounded number of palaces in the government $\mac G$-adjoin members of $C'$;
		\item no rebel house or fort that is not in $C'$ $\mac G$-adjoins a non-$C'$-peripheral fort in $C$, or a member of a non-$C'$-peripheral $C'$-village; and
		\item the number of $C'$-peripheral $C'$-villages is bounded. (There are automatically at most two $C'$-peripheral forts.)
	\end{itemize}
	The third condition is needed for \ref{cabalsaresmall}, and the fourth condition is needed to prove \ref{palacebound}.
	
	$C$ itself satisfies the first and third conditions. By simply choosing $C'$ to be a minimal subset of $C$ that
	is $\mac G$-village-closed,  $\mac G$-adjoin-connected and dangerous, we can arrange the second bullet and fourth bullets, but that might wreck the third,
	so we need to be more cautious. Indeed, there are some cases where no good subset exists. For instance, if $C$ consists of a 
	single fort or a single $\mac G$-village, and $\mac G$-adjoins many palaces in $\mac G$, there is nothing we can do; or if $|C|=6$,
	and consists of three forts and three $\mac G$-villages making a six-cycle under $\mac G$-adjoinment, and $C$ is minimally dangerous,
	then again there is no good subset. This explains why cabals sometimes have ``leader'' terms that are not peripheral. 
	
	Let us continue the proof of (1).
	If some $\mac G$-village included in $C$ is dangerous, then it is a cabal (designating itself as the only leader $\mac G$-village), a contradiction. So
	no $\mac G$-village included in $C$ is dangerous.
	Thus we assume (for a contradiction) that $C$ is dangerous and not a $\mac T$-community, and no $\mac G$-village in $C$ is dangerous.
	Since $C$ is $\mac G$-adjoin-connected, it is $\mac T$-adjoin-connected, and so we may number the forts in
	$C$ as $X_i\;(i\in I)$ where $I$ is an integer interval, as in \ref{contactgraph}.
	Since $C$ is a maximal $\mac G$-adjoin-connected set of rebels, it follows that $C$ is $\mac G$-village-closed.
	
	For $i_1,i_2\in I$ with $i_1\le i_2$, let $C(i_1,i_2)$ be the union of $\{X_{i_1}\LL X_{i_2}\}$ and all $\mac G$-villages that $\mac G$-adjoin
	one of $X_{i_1}\LL X_{i_2}$.
	It follows that
	$C(i_1,i_2)$ is $\mac G$-village-closed.
	Since $C$ is dangerous, we may choose $i_1,i_2$ with $i_2-i_1$ minimal such that $C(i_1,i_2)$ is dangerous. (Possibly $i_1=i_2$.)
	
	Define $D$ as follows:
	\begin{itemize}
		\item if $i_2=i_1$, let $D=\{X_{i_1}\}$;
		\item if $i_2=i_1+1$, let $D$ be the union of $\{X_{i_1}, X_{i_2}\}$ and all $\mac G$-villages that $\mac G$-adjoin  both of $X_{i_1}, X_{i_2}$;
		\item if $i_2\ge i_1+2$, let 
		$D$ be the union of $\{X_{i_1}\LL X_{i_2}\}$ and all $\mac G$-villages that $\mac G$-adjoin one of $X_{i_1+1}\LL X_{i_2-1}$.
	\end{itemize}
	In each case, 
	$D$ is $\mac G$-village-closed, and $\mac G$-adjoin-connected.
	Each member of $C(i_1,i_2)\setminus D$ is a house, and belongs to a $\mac G$-village that $\mac G$-adjoins one or both of $X_{i_1},X_{i_2}$ and none of
	$X_{i_1+1}\LL X_{i_2-1}$.
	
	Since $C(i_1,i_2)$ is dangerous, there exist $j>k$ and three members $A_1,A_2,A_3\in \mac G$, all of rank $j$, such that for $1\le i\le 3$, some member of $C(i_1,i_2)$
	$\mac G$-adjoins  $A_i$.
	For $i = 1,2,3$, $A_i$ $\mac G$-adjoins either some member of $D$, or some $\mac G$-village
	included in $C(i_1,i_2)\setminus D$.
	So by adding to $D$ at most three of the $\mac G$-villages included in $C(i_1,i_2)\setminus D$, we can construct  a set $D'$ that is dangerous,
	$\mac G$-adjoin-connected, and $\mac G$-village-closed. Let us construct such a set $D'$ by adding to $D$ as few $\mac G$-villages as possible;  so in particular, if $D$ is dangerous then $D'=D$.
	We designate $X_{i_1}, X_{i_2}$ as leader forts of $D'$, and the (at most three) $\mac G$-villages with union $D'\setminus D$ as leader $\mac G$-villages of $D'$.
	We claim that $D'$ is a cabal (which will be a contradiction). To show this, we must check that
	\begin{itemize}
		\item $D'$ is $\mac G$-adjoin-connected, and $\mac G$-village-closed;
		\item at most two forts of $D'$ are designated as leader forts, and at most three $\mac G$-villages of $D'$ are designated as
		leader $\mac G$-villages;
		\item every $D'$-peripheral fort is a leader fort, and every $D'$-peripheral $D'$-village is a leader $\mac G$-village;
		\item if $X\in D'$ $\mac G$-adjoins some rebel not in $D'$, then either $X$ is a leader fort, or $X$ belongs to a leader $\mac G$-village;
		\item for some $j\in \{k+1\LL \ell-1\}$ there are three palaces $A_1,A_2,A_3\in \mac G$ of rank $j$, such that for $1\le i\le 3$,
		$A_i$ $\mac G$-adjoins some member of $D'$; and
		\item either $|D'|=1$, or $D'$ is a $\mac G$-village, or for each $j'\in \{k+1\LL \ell-1\}$ there are at most four palaces $A\in \mac G$ of rank $j'$ such that
		$A$ $\mac G$-adjoins a member of $D'$.
	\end{itemize}
	The first two bullets are clear.
	The third holds since no member of $D$ is $D$-peripheral except possibly the forts $X_{i_1}, X_{i_2}$, and no $D$-village is $D$-peripheral. Hence the only $D'$-peripheral $D'$-villages are those (at most three) that we added.
	
	For the fourth bullet, suppose that $X\in D'$ $\mac G$-adjoins some rebel $Y\notin D'$.
	If $X$ is a fort, then $X=X_i$ for some $i\in \{i_1\LL i_2\}$; and since $Y\notin D$, it follows that $i\in\{i_1,i_2\}$  and so $X$ is a leader fort.
	Now we assume that $X$ is a house;  let $B$ be the $\mac G$-village of $\mac G$ that contains $X$. Since $D'$
	is $\mac G$-village-closed, it follows that $B\subseteq D'$. Since $Y\notin D'$, it follows that $Y\notin B$, and so $B$ $\mac G$-adjoins $Y$ and hence $Y$
	is a fort. Choose $i\in \{i_1\LL i_2\}$ such that $B$ $\mac G$-adjoins $X_i$, with $i\notin \{i_1,i_2\}$ if possible.
	If $i\ne i_1,i_2$, then $X_i$ $\mac T$-semiadjoins $X_{i-1}, X_{i+1}$ and $Y$, contrary to \ref{contactgraph}.
	Hence we cannot choose $i\notin \{i_1,i_2\}$; so $B\not\subseteq D$, and therefore $B$ is a leader $\mac G$-village.
	This proves the fourth bullet.
	
	The fifth bullet holds since $D'$ is dangerous. Finally, for the sixth bullet, 
	suppose first that $D'\ne D$, and let $B$ be a $\mac G$-village in $D'\setminus D$. Let $j'\in \{k+1\LL \ell-1\}$.
	Since we added as few $\mac G$-villages to $D$ as possible to make a dangerous set, $D'\setminus B$ $\mac G$-adjoins at most two  members of $\mac G$ of rank $j'$; and since
	$B$ is not dangerous, $B$ also $\mac G$-adjoins at most two  members of $\mac G$ of rank $j'$. Hence $D'$ $\mac G$-adjoins at most four 
	members of $\mac G$ of rank $j'$ and the sixth bullet holds. So we may assume that $D'=D$, and so $D$ is dangerous.
	If $i_1\ne i_2$, then for every palace $A\in \mac G$, if $A$ $\mac G$-adjoins some member
	of $D$
	then it also $\mac G$-adjoins a member of one of $C(i_1+1,i_2),C(i_1,i_2-1)$, and since $C(i_1+1,i_2)$ and $C(i_1,i_2-1)$ are not dangerous (by the minimality of $i_2-i_1$), it follows that
	the sixth bullet holds. So we may assume that $i_1=i_2$. Since $D=D'$, it follows that $D=\{X_{i_1}\}$ and so $|D|=1$, and again the sixth bullet holds.
	
	Hence $D'$ is a cabal, a contradiction since $\mac G$ is a stable government.
	This proves (1).
	
	\bigskip
	
	Thus, $C$ is not dangerous, and so 
	for $k+1\le j\le \ell-1$, there are at most two palaces $A\in \mac G$ with rank $j$ such that $A$ $\mac G$-adjoins some member of $C$.
	Let $\mac Q$ be the set of all such palaces in $\mac G$; so $|\mac Q|\le 2(\ell-k-1)$. But since $C$ is a maximal
	$\mac G$-adjoin-connected set of rebels, it follows that for each $u\in \bd(\voro_{\mac G}(C))$ there exists $A\in \mac Q$ such that $u$ has a neighbour in $\voro_{\mac G}(A)$;
	and so by \ref{palacebound}, $\bd(\voro_{\mac G}(C))$ has quasi-size at most
	$$\left(8(\ell-k-1) 3^{\ell-k-1}\alpha, \beta+4d_0+1\right),$$
	and the theorem holds.
	This proves \ref{notmuchleft}.~\bbox
	
	\section{Into a new century}\label{sec:newcent}
	
	Finally, we can move into a new century, and deduce the main theorem.
	\begin{thm}\label{realmtosoc}
		Let $k<\ell$ be a century, let $\sigma= (\delta, \alpha,\beta)$ be a 
		standard, and let $G$ be a graph admitting a $k$th-century realm with standard $\sigma$, and not containing $H_\ell$ as a $c$-fat minor.
		Let $\tau$ be the canon
		$$(\delta(2\ell), \ell\,3^{2\ell}\alpha, \beta+4d_0+1).$$
		Then $G$ admits a $(k+1)$st-century $\tau$-society.
	\end{thm}
	\Proof
	Since $G$ admits a $k$th-century realm with standard $\sigma$, 
	it admits a $k$th-century realm, $\mac T$ say, that is optimal for $\sigma$. Since there is the self-government for $\mac T$, there is a stable
	government $\mac G$ for $\mac T$. Let $\mac T'=\mac G$, where every house or fort of $\mac G$ is designated as a house of $\mac T'$, and every palace of $\mac G$ is designated as a fort of $\mac T'$.
	If $X\in \mac T'$, we write $\rk'(X)=k$ or $k+1$ depending whether $X$ is a house or fort of $\mac T'$.
	We claim that $\mac T'$ is a $(k+1)$st-century $\tau$-society.
	
	Let $\alpha'=\ell\,3^{2\ell}\alpha$ and $\beta'=\beta+4d_0+1$.
	We must check that
	\begin{enumerate}
		\item The members of $\mac T'$ are pairwise vertex-disjoint and induce connected subgraphs of $G$.
		\item For all $v\in V(G)$, there exists $X\in \mathcal{T}'$ such that $\dist_G(v,X)\le d_0$.
		\item If $X,Y\in \mac T'$ are distinct then $\dist_G(X,Y)$ is more than $\delta(2\ell)$.
		\item Every fort of $\mac T'$ includes a $c$-superfat $H_{k+1}$-minor of $G$.
		\item Each fort of $\mac T'$ has quasi-bound at most $(\alpha', \beta')$.
		\item For every $\mac T'$-village $C$,
		$\voro_{\mac T'}(C)$ has quasi-bound at most $(\alpha',\beta')$.
	\end{enumerate}
	Statement 1 is clear, and statement 2 holds since each member of $\mac T$ is a subset of some member of $\mac T'$.
	The third statement is clear.
	
	The fourth statement holds since if $X$ is a fort of $\mac T'$, then $X=A$ for some $A\in \mac G$; and since $A$ has
	rank some $t>k$, $A$ includes a $c$-superfat $H_t$-minor of $G$, and so also includes a $c$-superfat $H_{k+1}$-minor of $G$.
	
	For statement 5, let $X$ be a fort of $\mac T'$, and hence $X=A$ for some palace $A\in \mac G$.  Since $\mac G$ is a
	government, $A$ has quasi-bound at most
	$(\pi_t\alpha, \beta')$
	(where $t$ is the rank of $A$ in $\mac G$). But $\pi_t$ is at most
	$$3^{\ell-k-2}(10+8(\ell-k-1)3^{\ell-k-1}+1)=8(\ell-k-1)3^{2\ell-2k-3} + 11\cdot 3^{\ell-k-2}\le \ell\,3^{2\ell}.$$
	We deduce that statement 5 holds.
	
	Finally, for statement 6, let $C$ be a $\mac T'$-village. Thus, $C$ is a maximal $\mac G$-adjoin-connected set of rebels of 
	$\mac G$,
	so the statement follows from \ref{notmuchleft}. This proves \ref{realmtosoc}.~\bbox
	
	By combining \ref{realmtosoc} with \ref{bestsoctorealm} we deduce:
	
	\begin{thm}\label{bestsoctosoc}
		Let $k\ge 0$ be a century with $k<\ell$, and let $\tau=(d,\alpha,\beta)$ and $\tau'=(d', \alpha', \beta')$ be canons, 
		satisfying:
		\begin{align*}
			\alpha&\ge 12\\
			d'&\ge 5c\\
			d_0&\ge 4^{2\ell}d'\\
			d &= 4^{2\ell -2k+2} d'\\
			\alpha'&=\ell\,3^{2\ell}\alpha\\
			\beta' &= \beta+4d_0+1
		\end{align*}
		Let $G$ be a connected graph that does not contain $H_\ell$ as a $c$-fat minor. If $G$ admits a civilized $k$th-century $\tau$-society
		then it admits a $(k+1)$st-century $\tau'$-society. 
	\end{thm}
	\Proof
	Define $\delta(t) = 4^{2\ell-t}d'$ for $0\le t\le 2\ell$. Then $\delta$ is a spacing, since $\delta(2\ell)=d'\ge 5c$ and 
	$\delta(0)\le d_0$. Hence $\sigma=(\delta,\alpha,\beta)$ is a standard, because $\alpha\ge 12$.
	Since $d=\delta(2k-2)$ if $k\ge 1$, and $d\ge \delta(0)$ if $k=0$, \ref{bestsoctorealm} implies that $G$ admits a $k$th-century realm with standard $\sigma$.
	By \ref{realmtosoc}, $G$ admits a $(k+1)$st-century $\tau'$-society. This proves \ref{bestsoctosoc}.~\bbox
	
	In turn, by combining the previous result with \ref{newbestsoc}, we obtain:
	
	\begin{thm}\label{soctosoc}
		Let $k\ge 0$ be a century with $k<\ell$, and let $\tau=(d,\alpha,\beta)$ and $\tau'=(d', \alpha', \beta')$ be canons,
		satisfying:
		\begin{align*}
			\alpha&\ge 6\\
			d&\le d_0\\
			d'&\ge 5c\\
			d &\ge 4^{4\ell -2k+4} d'\\
			\alpha'&=2\ell\,3^{2\ell}\alpha\\
			\beta' &\ge  \beta+7d_0
		\end{align*}
		Let $G$ be a connected graph that does not contain $H_\ell$ as a $c$-fat minor. If $G$ admits a $k$th-century $\tau$-society
		then it admits a $(k+1)$st-century $\tau'$-society.
	\end{thm}
	\Proof It follows that $4^{2\ell}d'\le d_0$. Let $d''= 4^{2\ell -2k+2} d'$ (and so $d\ge 4^{2\ell+2} d''$). Let $\alpha''=2\alpha$ (and so $\alpha'=\ell\,3^{2\ell}\alpha''$).
	Let $\beta'' = \beta +2d_0+4^{\ell+1} d''$. Since
	$$\beta'\ge \beta+7d_0= \beta'' -4^{\ell+1} d''+5d_0$$
	and 
	$$4^{\ell+1} d''=4^{3\ell -2k+3} d'\le 4^{3\ell -2k+3} 4^{-4\ell +2k-4}d\le  4^{-\ell -1}d_0\le d_0-1.$$
	it follows that $\beta' \ge  \beta''+4d_0+1$.
	Since $G$ admits a $k$th-century $\tau$-society and $d\ge 4^{\ell+2}d''$, 
	it also admits a $k$th-century $(4^{\ell+2}d'',\alpha,\beta)$-society. 
	By \ref{newbestsoc}, $G$ admits a civilized $k$th-century $(d'', \alpha'', \beta'')$-society, since $\beta''\ge \beta+2d_0+4^{\ell+1} d''$. 
	By \ref{bestsoctosoc},
	$G$ admits a $(k+1)$st-century $\tau'$-society. This proves \ref{soctosoc}.~\bbox

	Now we can deduce the main result:
	\begin{thm}\label{mainthm4}
		For all integers $\ell\ge 1$ and $c\ge 0$
		there exist $\alpha,\beta$ such that every graph
		that does not contain $H_\ell$ as a $c$-fat minor has quasi-line-width at most $(\alpha,\beta)$.
	\end{thm}
	\Proof
	It suffices to prove the result when $\ell\ge 4$ and $c\ge 2$. 
	For $0\le k\le \ell$, let us define $d_k = 5\cdot 4^{5\ell(\ell-k)}c$ (note that this at last defines $d_0=5\cdot 4^{5\ell^2}c$), 
	and $\alpha_k= 6(2\ell\,3^{2\ell})^k$, and $\beta_k= 1+7kd_0$, and $\tau_k = (d_k, \alpha_k, \beta_k)$.
	Let $\alpha = \alpha_\ell$ and $\beta = \beta_\ell$; we claim that $\alpha, \beta$ satisfy the theorem. 
	By working with each component of $G$ separately, we may assume that $G$ is connected.
	By Zorn's lemma, there exists $S\subseteq V(G)$ maximal such that $\dist_G(u,v)> d_0$ for all distinct $u,v\in S$.
	Let $\mac T_0 = \{\{v\}:v\in S\}$ where each member of $\mac T_0$ is assigned to be a fort of $\mac T_0$.
	Then $\mathcal{T}_0$ is trivially a $0$th-century $\tau_0$-society. By $\ell$ applications of \ref{soctosoc} we deduce that 
	$G$ admits an $\ell$th-century $\tau_\ell$-society $\mac T$. Since $G$ does not contain $H_\ell$ as a $c$-fat minor, it follows that 
	$\mac T$ has no forts, and, since $G$ is connected, $\mac T$ has only one village $C$ say, with $\voro_{\mac T}(C)=V(G)$.
	Moreover, $\voro_{\mac T}(C)$ has quasi-bound at most $(\alpha_\ell,\beta_\ell)$, from the definition of a society.
	This proves \ref{mainthm4}.~\bbox
	
	Consequently, this proves \ref{mainthm3}, and hence \ref{mainthm2}, and together with \ref{quasitw} proves \ref{mainthm}.
	
	\section{Graph searching}
	
	For finite graphs $G$, these results are related to graph searching. 
	Intuitively, imagine that there is an infection loose in the graph; any vertex of the graph instantly becomes infected if it is
	adjacent to an infected vertex, unless a doctor is positioned on it. A search can be thought of as a sequence of
	positions of doctors in a procedure to eliminate the infection.  (Or, equally, a sequence of moves by a band of cops to capture an invisible robber.) 
	
	Let $G$ be a finite graph, let $X\subseteq V(G)$, and let $F$ be a subset of $V(G)$ with $\bd(F)\subseteq X\subseteq F$.
	We say that $(X,F)$ is a {\em split} of $G$. If $(X,F)$ and $(X',F')$ are splits,
	we say that $(X,F)$ {\em justifies} $(X',F')$ if $F\setminus X\subseteq F'$, and either $X\subseteq X'$ or $X'\subseteq X$.
	A {\em search} of $G$ is a finite sequence $(X_1,F_1)\LL (X_n,F_n)$ of splits of $G$ such that
	\begin{itemize}
		\item $(X_1,F_1) = (\emptyset,V(G))$, and $(X_n,F_n) = (\emptyset, \emptyset)$; and
		\item $(X_i,F_i)$ justifies $(X_{i+1},F_{i+1})$ for $1\le i\le n-1$.
	\end{itemize}
	We call the sets $X_i$ the {\em bags} of the search. (In terms of infections, at the $i$th stage, $X_i$ is the set of vertices
	being treated, and all infected vertices are in $F_i$. In cops and robbers language, at time $i$ the cops occupy the vertices in $X_i$,
	and the robber is known to be within $F_i$.)

	For instance, let $(B_1\LL B_n)$ be a path-decomposition of a finite graph $G$. Define $R_i = B_i\cup B_{i+1}\cupcup B_n$
	for $1\le i\le n$. Then the sequence
	$$(\emptyset,R_1),(B_1, R_1), (B_1\cap B_2, R_2), (B_2,R_2), (B_2\cap B_3, R_3),(B_3,R_3)\LL (B_{n-1}\cap B_n, R_{n}), (B_n,R_n), (\emptyset,\emptyset)$$
	is a search.
	
	A search $(X_1,F_1)\LL (X_n,F_n)$ is {\em monotone} if $F_{i+1}\subseteq F_i$ and $F_{i+1}\cap X_i\subseteq X_{i+1}$ for $1\le i<n$.
	The search above, derived from a path-decomposition, is monotone, and it is easy to see that every monotone search arises from a
	path-decomposition in this way, so monotone searches are essentially the same as path-decompositions.
	
	One important parameter of a search $(X_1,F_1)\LL (X_n,F_n)$ is its {\em bag-size} $\max(|X_1|\LL |X_n|)$. If we fix some number $k$, when is there a
	search with bag-size at most $k$? It was proved by Kirousis and
	Papadimitriou~\cite{kirousis1,kirousis2}, building on work of LaPaugh~\cite{lapaugh} (who proved the same thing for a slightly different kind of search), that:
	\begin{thm}\label{getmonotone}
		For each integer $k$ and finite graph $G$, if there is a search with bag-size at most $k$,
		then there is a monotone search with bag-size at most $k$.
	\end{thm}
	Consequently, there is a search with bag-size at most $k$ if and only if $G$ has path-width at most $k-1$.
	
	Here we are concerned with searches in which each bag has quasi-size at most $(a,b)$, for some
	fixed $(a,b)$. We say that a finite graph $G$ is {\em $(a,b)$-searchable} if there is a search $(X_1,F_1)\LL (X_n,F_n)$
	of $G$ such that $X_i$ has quasi-size at most $(a,b)$ for $1\le i\le n$.
	It is easy to prove that for all $a,b$, there exists $c,\ell\ge 2$ such that if $G$ is $(a,b)$-searchable then $G$ does not contain
	$H_{\ell}$ as a $c$-fat minor. (Proved in \ref{easyhalf}, below).
	
	Our result \ref{mainthm3} gives a converse:
	\begin{thm}\label{searching}
		For all $c,\ell\ge 2$, there exist $a,b$ such that if a finite graph $G$ does not contain
		$H_{\ell}$ as a $c$-fat minor, then $G$ is $(a,b)$-searchable (because it has bounded quasi-line-width). 
	\end{thm}
	
	We want to make two comments. First, we also have a direct proof of \ref{searching} (for finite $G$), that we found before we found the proof of 
	\ref{mainthm3}, and it is considerably simpler. The idea is, to replace the assertions that certain sets have bounded quasi-bound
	by the weaker assertions that these sets are {\em subsets} of sets with bounded quasi-bound. Then all our problems with communities go 
	away, and so we don't need the sections about societies. Except the proof of \ref{composeline} does not work any more; but it works if 
	we replace
	quasi-line-width with being $(a,b)$-searchable with $a,b$ bounded, and this is why we have to introduce graph searching. So, it would 
	be really nice if we could find a direct proof of:
	
	\begin{thm}\label{eqnce}
		For all $a,b\ge 1$ there exist $a',b'\ge 1$ such that if a finite graph is $(a,b)$-searchable then it has quasi-line-width at most $(a',b')$.
	\end{thm}
	We know that this is true, because of the main theorem of this paper, but if a simple direct proof could be found, it would give a 
	simpler proof of \ref{mainthm3}.
	
	Second, how do $a',b'$ depend on $a,b$ in \ref{eqnce}? In particular, can we take $a=a'$? Yes we can if $b=0$, by \ref{getmonotone}, 
	and in that case we can also take $b'=0$. But what happens if $b>0$? For instance, if $b=1$, can we take $a'=a$?
	
	On this topic, let us mention that there is an interesting conjecture of Chudnovsky,  Gollin, Krnc, and Milani\v{c}~\cite{balanced} 
	with a special case that is closely related to our results, the following:
	\begin{thm}\label{mariaconj}
		{\bf Conjecture:} For every tree $T$ there exists $k$ such that if a finite graph $G$ has no induced minor isomorphic to $T$, then there is  a balanced separator in
		$G$ that is the union of at most $k$ balls of radius one.
	\end{thm}
	An ``induced minor'' is much the same as a 1-fat minor: $G$ contains $H$ as an {\em induced minor} if $H$ can be obtained from $G$ by 
	vertex-deletion and edge-contraction, and the deletion of all resultant loops and parallel edges. A {\em balanced separator} in $G$ is a cutset 
	$X$ such that each component of $G\setminus X$ contains at most $|V(G)|/2$ vertices. It is easy to see that if $G$ is $(a,1)$-searchable 
	for some bounded $a$ then the desired balanced separator exists. 
	
	Finally, here is a proof of the claim we made earlier in this section:
	
	\begin{thm}\label{easyhalf}
		For all $a,b\ge 1$, no finite graph that contains $H_{4a}$ as a $2b$-fat minor is $(a,b)$-searchable.
	\end{thm}
	\Proof Given integers $a,b$, let $c=2b$, and let $H=H_{4a}$,
	with root $r$. 
	Let $G$ be a finite graph that contains $H$
	as a $c$-fat minor, and choose $\eta$ that exhibits $H$ as a $c$-fat minor of $G$.
	
	Suppose that there is a search $(X_1,F_1)\LL (X_n, F_n)$ of $G$ such that each $X_i$ has quasi-size at most $(a,b)$.
	For $1\le i\le n$, let $Y_i$ be the set of $h\in V(H)$ such that either
	\begin{itemize}
		\item $V(\eta(h))\cap X_i\ne \emptyset$; or
		\item $h$ is an end of some edge $f\in E(H)$ such that $V(\eta(f))\cap X_i\ne \emptyset$.
	\end{itemize}
	There exists $Z_i\subseteq V(G)$ with $|Z_i|\le a$ such that every vertex $x\in X_i$ satisfies $\dist_G(x,Z_i)\le b$.
	Since $\eta$ exhibits $H$ as a $c$-fat minor of $G$, and $c=2b$, it follows that if $x_1,x_2\in V(H)$ are distinct vertices of $H$
	then $\dist_G(\eta(x_1),\eta(x_2))> c$; and so for each $z\in Z_i$, there is at most one vertex $x\in V(H)$ with $\dist_G(z,\eta(x))\le b$.
	Similarly, there is at most one edge $f\in E(H)$ such that $\dist_G(z,\eta(f))\le b$, and if there is one of each type, then the edge and vertex are incident in $H$. Consequently $|Y_i|\le 2a$.
	\\
	\\
	(1) {\em For $1\le i\le n$, if $Q$ is a component of $H\setminus Y_i$, and $F_i\cap V(\eta(x))\ne \emptyset$
		for some vertex or edge $x$ of $Q$, then
		$V(\eta(x))\subseteq F_i\setminus X_i$ for every vertex or edge $x$ of $Q$.}
	\\
	\\
	Let $R$ be the subgraph of $G$ induced on the union of $\eta(x)$ over all vertices and edges $x$ of $Q$.
	Thus $R$ is
	connected since $Q$ is connected. For each $v\in V(R)$, if $v\in X_i$, choose a vertex or edge $x$ of $Q$ with $v\in V(\eta(x))$; then
	$x\in Y_i$ (if $x$ is a vertex) or one end of $x$ is in $Y_i$ (if $x$ is an edge), and in either case some vertex of $Q$ is in $Y_i$,
	contradicting that $Q$ is a component of $H\setminus Y_i$. Thus $V(R)\cap X_i=\emptyset$.
	By hypothesis,  $F_i\cap V(\eta(x))\ne \emptyset$
	for some vertex or edge $x$ of $Q$. But $V(\eta(x))\subseteq V(R)$, and so $V(R)\cap F_i\ne \emptyset$.
	Hence, $V(R)\subseteq F_i\setminus X_i$, since $(X_i,F_i)$ is a split of $G$ and $R$ is connected and $V(R)\cap X_i=\emptyset$. This proves (1).
	
	\bigskip
	
	It follows from (1) that if $Q$ is a component of $H\setminus Y_i$ and $F_i\cap V(\eta(x))\ne \emptyset$ for some vertex or edge $x$ 
	of $Q$, then $V(\eta(x))\subseteq F_i$ for every vertex or edge $x$ of $Q$. 
	For $1\le i\le n$, let $P_i$ be the union of $Y_i$ and the vertex sets of all components $Q$ of $H\setminus Y_i$ with the property that
	$F_i\cap V(\eta(x))\ne \emptyset$
	for some vertex or edge $x$ of $Q$. 
	
	Now each $(Y_i, P_i)$ is a split of $H$.
	We claim that:
	\\
	\\
	(2) {\em $P_i\subseteq P_{i+1}\cup Y_i$ for $1\le i<n$.}
	\\
	\\
	To see this, let $h\in P_i\setminus Y_i$; we need to show that $h\in P_{i+1}$, and may therefore assume that $h\notin Y_{i+1}$.
	Let
	$Q_i, Q_{i+1}$ be the components of $H\setminus Y_i, H\setminus Y_{i+1}$ respectively that contain $h$.
	Since $h\in P_i$, it follows that
	$V(\eta(x))\subseteq F_i\setminus X_i\subseteq F_{i+1}$ for every vertex or edge $x$ of $Q_i$. Hence,
	$F_{i+1}\cap V(\eta(x))\ne \emptyset$
	for some vertex or edge $x$ of $Q_{i+1}$ (namely, $x=h$), and so $V(Q_{i+1})\subseteq P_{i+1}$. This proves (2).
	
	\bigskip
	
	From (2), and since $(Y_1,P_1) = (\emptyset, V(H))$ and $(Y_n,P_n) = (\emptyset, \emptyset)$,
	the sequence $(Y_1,P_1)\LL (Y_n,P_n)$ is a search of $H$. But $H=H_{4a}$, and therefore has
	path-width at least $2a+1$ (see \cite{GM1}), and so is not $(2a+1,0)$-searchable, by a result of Kirousis and
	Papadimitriou~\cite{kirousis1, kirousis2}. Yet $Y_1\LL Y_n$ each have cardinality at most $2a$, a contradiction.
	This proves \ref{easyhalf}.~\bbox
	
	\section*{Declaration of generative AI use}
	We acknowledge the use of ChatGPT in proofreading and checking the manuscript. All the mathematical ideas and arguments presented, and the writing, are due to the authors.
	
	\section*{Remarks on open access}
	For the purpose of open access, the authors have applied a CC BY public copyright licence to
	any author accepted manuscript arising from this submission.

\end{document}